\pdfoutput=1
\documentclass[twoside,12pt]{report}
\usepackage{tikz}
\usepackage{cite}
\usepackage{url}
\usepackage{amsmath}
\usepackage{amssymb}
\usepackage{amsthm}
\usepackage{latexsym}
\usepackage{graphicx}
\usepackage{float}
\usepackage{epstopdf}
\usepackage{verbatim}
\usepackage[section]{placeins}
\usepackage{amsfonts}
\usepackage{color}
\usepackage{calligra}
\usepackage[T1]{fontenc}
\usepackage{lmodern}
\usepackage[utf8]{inputenc}
\usepackage[]{xcolor}
\usepackage{bm}
\usepackage{subcaption}
\usepackage{enumerate}
\usepackage[protrusion=true,expansion=true]{microtype}
\usepackage{fancyhdr}
\usepackage[hmarginratio=1:1,left=1.5in,right=1.5in]{geometry}
\usepackage[bottom]{footmisc}
\usepackage[pdfborder={0,0,0},pdfauthor={Alexander Giles}, colorlinks=true,
linkcolor=blue, citecolor=blue]{hyperref}
\usepackage{setspace}

\usepackage{graphicx,amsmath,amsfonts,amssymb,amscd,bezier,amstext,makeidx,mathtools,xcolor,siunitx}
\usepackage{amsthm}
\usepackage{bm}
\usepackage{color}
\usepackage{latexsym}
\usepackage{graphicx}
\usepackage{amsthm,amssymb,amscd}
\usepackage{longtable}

\usepackage{rotating}
\usepackage{tikz}

\usepackage{pifont}

\usepackage[toc,page]{appendix}

\definecolor{apricot}{RGB}{253, 213, 177}

\DeclarePairedDelimiter\norm{\lVert}{\rVert}
\DeclarePairedDelimiter\singlenorm{|}{|}

\providecommand{\floor}[1]{\left \lfloor #1 \right \rfloor }
\providecommand{\ceil}[1]{\left \lceil #1 \right \rceil }

\DeclareMathOperator{\arccosh}{arcosh}

\newtheorem{lemma}{Lemma}
\newtheorem{theorem}[lemma]{Theorem}

\newtheorem{conjecture}[lemma]{Conjecture}
\numberwithin{equation}{section} \numberwithin{lemma}{section}

\makeatletter
\newcommand{\contraction}[5][1ex]{%
  \mathchoice
      {\contraction@\displaystyle{#2}{#3}{#4}{#5}{#1}}%
      {\contraction@\textstyle{#2}{#3}{#4}{#5}{#1}}%
      {\contraction@\scriptstyle{#2}{#3}{#4}{#5}{#1}}%
      {\contraction@\scriptscriptstyle{#2}{#3}{#4}{#5}{#1}}}%
\newcommand{\contraction@}[6]{%
  \setbox0=\hbox{$#1#2$}%
  \setbox2=\hbox{$#1#3$}%
  \setbox4=\hbox{$#1#4$}%
  \setbox6=\hbox{$#1#5$}%
  \dimen0=\wd2%
  \advance\dimen0 by \wd6%
  \divide\dimen0 by 2%
  \advance\dimen0 by \wd4%
  \vbox{%
    \hbox to 0pt{%
      \kern \wd0%
      \kern 0.5\wd2%
      \contraction@@{\dimen0}{#6}%
      \hss}%
    \vskip 0.2ex%
    \vskip\ht2}}

\newcommand{\contraction@@}[3][0.06em]{%
  \hbox{%
    \vrule width #1 height 0pt depth #3%
    \vrule width #2 height 0pt depth #1%
    \vrule width #1 height 0pt depth #3%
    \relax}}

\definecolor{deepred}{RGB}{209,34,80}
\definecolor{darkgreen}{RGB}{18,194,15}
\definecolor{amber}{RGB}{222,194,13}
\definecolor{darkergreen}{RGB}{0,133,0}
\definecolor{betagreen}{RGB}{128,255,128}
\definecolor{hotpink}{RGB}{85,26,139}
\begin{document}
\fancyhead[ER]{\fancyplain{}\leftmark{}}
\fancyhead[OR]{\fancyplain{}{\thepage}}
\begin{titlepage}
  \begin{center}
{\large \textbf{Connectivity and Centrality in \protect\\ Dense Random Geometric
Graphs}}
    \\[0.5cm]
      {\large {Alexander P. Kartun-Giles}}\\[1.5cm]
      \includegraphics[width=4cm]{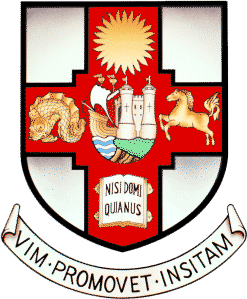}\\[0.75cm]
      \textsc{University of Bristol}\\
      \textsc{EPSRC Center for Doctoral Training in Communications}\\[1cm]
\small{\textit{A dissertation submitted to the University of Bristol in
accordance with the requirements for the degree of Doctor of Philosophy in
the Faculty of Engineering}}\\[0.2cm]
      March 2017.
  \end{center}
\end{titlepage}
\pagenumbering{roman}
\newpage
\begin{center}$\dagger$\end{center}
\newpage

\begin{abstract}
Due to shorter range communication becoming more prevalent with the development of multiple-input, multiple-output antennas (MIMO) and millimeter wave communications, multi-hop, intra-cell communication is anticipated to play a major role in 5G. This is developed in this thesis. Our analysis involves a stochastic spatial network model called a
\textit{random geometric graph}, which we use to model a network of interconnected devices communicating wirelessly without any separate, pre-established infrastructure.

Extending recent work on the \textit{random connection model} by relaxing the
convexity restriction on the bounding geometry $\mathcal{V}$, we study the
asymptotic connectivity of a model of randomly distributed, randomly linked vertices embedded within an
annulus. Vertices, resembling communicating devices, are linked according to a Rayleigh fading model of signal attenuation in built-up urban environments. We notice that vertices near obstructions in the
domain have an exceptionally high betweenness centrality, and highlight the
need to quantify this sort of measure analytically. Following up
on this, in a theoretical \textit{continuum limit}, where vertex density goes
to infinity and the connection range goes to zero compared with the size of
the domain, we approximate betweenness centrality in a random network
with a known special function called the elliptic integral of the second
kind. This provides analytic support to various delay tolerant networking
protocols currently used in 5G scenarios, wireless sensor networks, and
intelligent vehicle networks, which avoid what are often intractably complex algorithmic computations, even when distributed over many mobile processors.

Due to the obvious break with our model at lower densities, we then relax the
continuum limit, and investigate the expected number of shortest paths which
run between two spatially separated vertices in a unit disk graph in order to develop an understanding of betweenness centrality at finite density, which is based on counting shortest paths in networks. It also develops so called range-free localisation beyond estimating distances to
location-aware anchors via hop counts, instead using the more sensitive
\textit{number} of geodesic paths. We conclude with the difficulties of
combinatorial enumeration in spatial networks, and highlight the need to
develop probabilistic techniques to deal with the deep issues of spatial dependence.
\end{abstract}

\newpage
\begin{center}$\dagger$\end{center}
\newpage

\section*{Acknowledgements}
I acknowledge with great thanks and gratitude the help and support of Professor Carl P.
Dettmann, Dr Orestis Georgiou of Toshiba Research Europe, and Professor Justin P. Coon of the Oxford University Communications Research Group, all of whom have provided the supervision of this thesis. Also, thanks to everyone at the CDT, including Bagots, Brett, Dave, Leo, Suzanne, Mark, Simon, and everyone else, and finally, thanks to all those at 26 Cadogan Square.
\newpage
\begin{center}$\dagger$\end{center}
\newpage
\section*{Author's Declaration}
\vspace{5mm}
I declare that the work in this dissertation was carried out in accordance
with the requirements of the University’s Regulations and Code of Practice
for Taught Programmes and that it has not been submitted for any other
academic award. Except where indicated by specific reference in the text,
this work is my own work. Work done in collaboration with, or with the
assistance of others, is indicated as such. I have identified all material in
this dissertation which is not my own work through appropriate referencing
and acknowledgement. Where I have quoted or otherwise incorporated material
which is the work of others, I have included the source in the references.
Any views expressed in the dissertation, other than referenced material, are
those of the author.\\

\noindent Alexander Paul Kartun-Giles\\
4th of March 2017
\newpage
\begin{center}$\dagger$\end{center}
\newpage
\tableofcontents

\listoffigures

\listoftables

\newpage
\begin{center}$\dagger$\end{center}
\newpage
\pagenumbering{arabic}
\chapter{Introduction}\label{c:intro}
Ultra-dense spatial deployment of cellular base stations is among the most promising ways in which the data capacity of large scale wireless networks is to be enhanced \cite{ge2016}. The density of this deployment is highly anticipated, as of 2016, to be around 40-50 base stations per km\textsuperscript{2}. New technology introduced into these units will result in their communication range falling by an order of magnitude compared to the 4G base stations which are currently deployed in an ad hoc fashion around our cities. This is partly due to the appearance of
multiple-input multiple-output (MIMO) antennas, which offer faster data rates (Gbit/s), but utilise over a hundred antennas simultaneously, each at reduced power \cite{marzetta2010}. Compounding this, the anticipated high-frequency millimeter wave technologies in the 30-300 Ghz range will attenuate very quickly, allowing only short range communication between devices \cite{rappaport2013}. Therefore, in order to increase the range of special gateway-enabled cells,
networks will have to employ a form of multi-hop, cell-to-cell communication. According to this much
anticipated theory, 5G is therefore likely to consist, at least in part, of densely deployed,
millimeter wave connected `small' cells routing data from mobile devices between themselves in a sort of daisy chain toward gateway cells, and then finally into the backhaul network, enabling much higher throughput networks \cite{gupta2015}.

Now, these small cells will form a spatially embedded network.
Combinatorial analysis of this sort of system has long since been a topic of
great research interest. In this thesis, we contribute to this growing field
by analysing a combinatorial object called a
\textit{random geometric graph}, an example of which is shown in Fig.
\ref{fig:rgg}, in order to perform statistical network analysis developing future intra-cell communications. We will focus on two models:
\begin{itemize} 
\item \textbf{Soft random geometric graph}, also known as the \textbf{Random Connection Model} when embedded in the entire real plane rather than a bounding geometry such as a disk, is a type of random network formed by distributing
vertices in a region $\mathcal{V} \subseteq \mathbb{R}^{d}$ according
to a Poisson point process $\mathcal{Y}$ of density $\rho$, and then adding
an edge between points $\{x,y\} \in \mathcal{Y}$ with probability
$H(\norm{x-y})$, where $H: \mathbb{R}^{+} \to \left[0,1\right]$ is the
connection function, and $\norm{x-y}$ is the distance between vertices given
by some metric \cite{penrose2016}. We often use the Rayleigh fading connection function $H(\norm{x-y}) =\exp{\left(-\beta\norm{x-y}^2\right)}$ defined in Section \ref{sec:rayleigh}. This is parametrised by $\beta \in \mathbb{R}^{+}$, and we then write $r_0=1/\sqrt{\beta}$ to denote the typical distances over which vertices connect, when the context is clear.

\item The \textbf{Unit disk model} is the limit of this model where the
connection function $H(\norm{x-y}) \to \mathbf{1}\{\norm{x-y} < r_0\}$, with
$r_0 \in \mathbb{R}^{+}$ the critical range over which vertices connect \cite{penrosebook}.
\end{itemize}

These random graphs have been used to model the following two types of wireless network:

\begin{figure}[t]
\begin{center}
\includegraphics[scale=0.7]{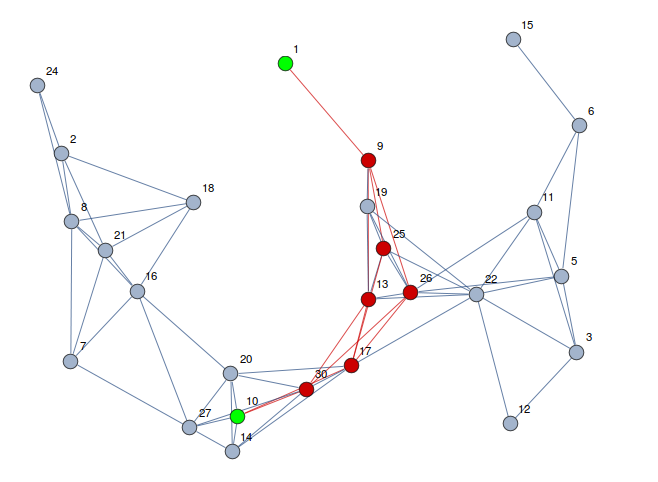}
\caption{An example random geometric graph. These combinatorial objects can
be used to model networks of millimeter wave-connected 5G base stations in
ultra-dense urban deployment \cite{ge2016}. It has been formed from a Poisson point process
by adding an edge between vertices whenever they fall within a Euclidean
distance of 0.25 of each other. All six shortest paths from vertex 1 to 10
have been highlighted in red as a structure of intersecting paths.}
\label{fig:rgg}
\end{center}
\end{figure}

\begin{itemize}
\item \textbf{Ad hoc wireless networks}, where devices communicate, as
discussed, between themselves without utilising separate, pre-established
infrastructure. This is achieved by routing data in a `multi-hop\footnote{A hop is an edge in a graph. One \textit{hops} between vertices, hence the commonly used term.}’ fashion
between users. Apart from in future 5G urban networks, they are seen in
deployed in disaster zones since they can offer a larger range than simple
point-to-point radio communication. There are whole journals related to this topic, but see e.g. \cite{farago2012} for a good discussion.

\item \textbf{Wireless sensor networks}, which are densely deployed networks
of low-power sensors. They are used for many different sensing tasks.
Wirelessly communicating with each other, sensors with limited battery
resources and limited computational performance collect data from an
environment and route it in a multi-hop fashion toward an elected
\textit{sink} sensor. The sink then sends the collected information toward a
cloud or distant base station. For for further discussion, see e.g. \cite{li2009}.
\end{itemize}

The major result which has lead to such great interest in random networks
around the world was published surrounding a 1957 conference on Monte Carlo
techniques in the UK. Broadbent and Hammersley showed that a non-trivial,
first order phase transition of `sudden connectivity' occurred when an array
of vertices is adorned at random with links between neighbouring pairs with a
critically high probability $p$ \cite{broadbent1957}. This is called
\textit{percolation}\footnote{We introduce this model with greater detail in Appendix \ref{appendix:percolation}.}. Interestingly, with an infinite array of vertices, the
phase transition is `sharp' \cite{grimmettpercolation}. Mathematical study of
this simple model has lead to many important results in modern science, with
Fields medals awarded recently for proving certain results concerning a
conformal invariance property of extremely fine lattices exactly at the
moment of percolation, called \textit{criticality} \cite{grimmettbook}. Similarly, both the unit disk and random connection model mentioned above
display a percolation phenomenon\footnote{More formally, taking
$H(\norm{x-y})= \mathbf{1}\{\norm{x-y} < r_0\}$, in a unit disk graph of
expected $n$ vertices distributed over a square of side $\sqrt{n}$ there
exist supercritical connection ranges $r_0(n) > r_c(n)$ for which a giant
cluster exists asymptotically almost surely (i.e. with probability $1$ in the
limit $n \to \infty$). There exists a similar asymptotic connection rule in
the random connection model \cite{penrosebook,mao2011}.}. 

Now, these graphs are
excellent models of randomly constructed
wireless networks of radio transmitters dispersed over a large
geographic region \cite{gilbert1961}. One can therefore investigate analogous results to percolation,
but in a wireless network setting \cite{cef2012}. Indeed, this is the most common
application of the theory today. Results, generally, concern theoretical
methods for the optimal use of what are now very expensive communication
resources.

The field is broad:

\begin{enumerate}
\item ``(There) is a sense of missing unification and a lack of general
methods that apply to large families of different models.''
\cite{farago2012,flaxman2007}.

\item The task of then applying mathematical results such as performance
bounds to actual networks and the way they run is not straightforward, nor is
it commonly performed in industry. This has lead to a lack of focus.

\item The probabilistic analysis of spatial communication networks, sometimes
involving dynamical processes, is often very difficult. Well developed
techniques from random graph theory can rarely be applied effectively
\cite{bollobas2005}.
\end{enumerate}

Despite these difficulties, there are well cited results over the last two
decades which constitute a research frontier, either concerning mathematical
features of commonly applied random graph models, or the specific use of
random graphs in communication problems.
\section{Current state of the art}\label{sec:stateoftheart}
As a broad introduction to a rich field, we first introduce some of the recent trends in the developing theory of random networks. Then, in Chapter \ref{c:concepts}, we focus on topics which provide the background to the contribution of this thesis, focusing on the areas of connectivity and centrality in communication networks.
\begin{enumerate}
\item \textbf{Mobility models:} Networks with moving vertices are used as
models for \textit{mobile} ad hoc networks (MANETs). For example, consider the
\textit{random walk model}: vertices are distributed uniformly at random on
the surface of a flat torus, but move, at each time step, to a new location a
distance $d$ (a parameter of the model) and angle $\theta$ (selected
uniformly at random for each time step) away from where they currently are
\cite{diaz2009}. One can then ask certain questions related to the communication theory of mobile networks,
such as to what extent a static graph can be used as an approximate model of a mobile systems, amongst related topics.

Mobility studies, for example, is particularly important in the case of \textit{delay tolerant
networking} \cite{magaia2015}. Here, mobile users offload and collect their data 
from a network of intermediate (and also mobile) devices, which store,
cary and forward data as the relay opportunity arises. The specific mobility
model on which the dynamic network is based will significantly alter the
analysis of performance. State of the art research concerns the more
realistic models of human mobility, such as the Levy Walk model (LW)
\cite{rhee2011}, and \textit{Self-Similar Least Action Human Walk} (SLAW)
\cite{lee2012}, both of which capture, for example, the self-similarity of
human walks, and the distribution of inter-meeting times within which message
relays can occur.

Also, the use of betweenness centrality \cite{freeman1977}, is already established within this field.
Numerical calculation of these network centrality indices can assist certain
delay tolerant networking protocols by encouraging mobile vertices to relay
data toward exceptionally central mobile vertices whenever possible
\cite{magaia2015}.

\item \textbf{Symmetric motifs:} A symmetric motif is a potentially
disconnected subset of vertices with the property that any permutation of
indices preserves adjacency \cite{macarthur2009}. They can be identified in
the graph Laplacian by identifying the integer eigenvalues and reading off
the non-zero components of corresponding eigenvectors.  The induced
subgraph\footnote{The graph consisting of those vertices, and edges whose end
points form a subset of those vertices.} of these components is the
corresponding symmetric motif \cite{nyberg2015}.

They occur more frequently in spatial networks due to the geometric factor. Their characteristic nature, mixed with
a simple connection to the graph's Laplacian spectrum, make them a
fascinating and much undervalued research approach to processes running on
dense communication networks.

\item \textbf{The Rado graph:} Infinite random geometric graphs concern
networks which have an infinite number of vertices, rather than specifically taking the
\textit{thermodynamic} or high-density \textit{continuum limit}. As an example, the unit disk graph has a version of this, taking the vertex set $V$ to be a countably infinite subset
of the Euclidean metric space, and $\norm{V} \to
\infty$. Graphs formed this way are not necessarily isomorphic, contrasting a
similar theory concerning non-spatial random graphs (the Erd\H{o}s-R\'{e}nyi
or \textit{random graph} \cite{grimmettbook}) presented in one of the
earliest works on the subject, by Erd\H{o}s and R\'{e}nyi \cite{erdos1959},
where all infinite random graphs are of a unique isomorphism type denoted the
\textit{Rado graph}. For the geometric case, this infinite limit can help
reveal large-scale structure and long term behaviour \cite{bonato2016}.

\item \textbf{Non-Homogeneous point processes:} The simple Poisson point
process model can be extended to address more realistic situations where the
distribution of points is not uniform in space. The most common variant is
the \textit{Strauss process}, which is a point process defined by a Markov
chain Monte Carlo algorithm involving repelling points; densely deployed
stations are randomly dispersed, but never right next to each other, which is
the main lack of realism encountered in the Poisson point process model. See
Sec \ref{sec:bpp}, \ref{sec:mpp} and \ref{sec:gpa}.

\item \textbf{Measures of centrality:} This involves introducing centrality
analysis into communication networks \cite{kog105}. Example
applications include boundary detection, and as part of routing
protocols which assign centrality indices to vertices in order to asses a
sort of `routing potential'. Also, techniques which networks can employ to
determine the centrality of their constituent vertices are of current
interest \cite{lulli2015}.

\item \textbf{Aerial networks:} One may consider employing aerial base
stations, such as drones. This can provide spontaneously adaptive coverage.
Exactly whether or not this solution can outperform ground station-reliant
networks can be assessed using analysis of random geometric networks
\cite{mozaffari2015}.

\item \textbf{Localisation:} Wireless devices often demand essential location
information. Indoor localisation cannot use GPS, given a line-of-sight
requirement with multiple satellites. Underground, aerial and secret military
localisation techniques are currently being developed. In a random network,
one can \textit{multilaterate} by first measuring hop-counts to location-aware
anchor vertices (such as access points), and then converting hop-counts to
Euclidean distances via a conditional distribution parameterised by Euclidean distance. \cite{nguyen2015}.

\item \textbf{Secure communication:} The use of random secrecy graphs, which
consist of a superimposed point process of eavesdroppers and users on which a
geometric graph is formed, can provide useful limits on the ability of
wireless networks to communicate in secure environments \cite{pinto2012}.

\end{enumerate}

\section{What is done in this thesis}
In this thesis, we focus on two key, developing areas:
\begin{enumerate}
\item[1)] \textit{Connectivity in non-convex domains}: Since many real world
domains are not convex subsets of $\mathbb{R}^{2}$, we extend the analysis of
high-density random network connectivity into the more realistic scenario of
a domain containing randomly arranged circular obstacles. As with other work
concerning non-convexity \cite{orestis2013,almiron2013}, the obstructions in
the domain introduce specific features to the enclosed graph. Specifically, obstacles
act like internal perimeters, encouraging isolated vertices. We describe this
effect in detail in chapter \ref{chapter:connectivity}, showing how small
obstacles\footnote{Small and large here mean on a length scale significantly
lesser or greater than the typical distance $r_0$ over which vertices
connect.} encourage isolated vertices according to a factor proportional to
their area $\frac{\pi^{d/2}}{\Gamma\left(\frac{d}{2}+1\right)}
(r_{\text{obstacle}}/r_{0})^{d}$, while large obstacles attract isolated
vertices according to a factor proportional to their perimeter
$\frac{\pi^{d+1/2}}{\Gamma\left(\frac{d+1}{2}\right)}
(r_{\text{obstacle}}/r_{0})^{d-1}$. This implies small obstacles have a
negligible effect on connectivity, since one sees a power of the ratio
$r_{\text{obstacle}}/r_{0}$, and hence the contribution vanishes for small
obstacles $r_{\text{obstacle}} \ll r_{0}$. Also, we show the compound effect
of obstacles is a linear combination of separate contributions, given they
are not too close.

We then observe that vertex isolation near an inner boundary can be more likely than other parts of the domain, and more likely than the outer boundary when the obstacles are numerous.
This is because vertices on the inner perimeter have an exceptionally high
betweenness centrality, acting as popular bridging vertices between distant
parts of the domain.\cite{freeman1977}.
\end{enumerate}
Noticing the importance of this sort of structural concept to communication
processes running on dense urban networks, we begin to analyse it directly.
\begin{enumerate}
\item[2)] In Chapter \ref{chapter:betweenness}, we consider analytically
quantifying betweenness centrality in a random geometric graph. In a limiting
density scenario, which we describe as a \textit{continuum limit}, shortest
paths between vertices are approximated by the convex hulls of their
endpoints i.e. by straight line segments. Based on the assumption that only a
single geodesic path will join any two vertices, by simply counting the
number of convex hulls which 1) can be formed by any pair of points in $x,y
\in \mathcal{V}$ and 2) run through $k$ we can then count the number of
geodesic paths which intersect $k \in \mathcal{V}$.

Using delta functions, we then show how the (expected\footnote{Note in our
continuum limit, betweenness and expected betweenness are equal.},
normalised) betweenness centrality $g^{\star}\left(k\right)$ of some polar
point $k = (\epsilon,\theta)$ in a disk domain metric space is, in fact, a
known integral:
\begin{eqnarray}
g^{\star}\left(k\right) =
\frac{2}{\pi}(1-\epsilon^2)\int_{0}^{\pi/2}\sqrt{1-\epsilon^2\sin^2\left(\theta\right)}\mathrm{d}\theta
\end{eqnarray} 
i.e. the elliptic integral of the second kind scaled by a quadratic function
of displacement from the disk's center. We then discuss numerical convergence
of expected betweenness to this integral as the point process density goes to
infinity, showing how centrality is well approximated, but never exactly
equal to, our convex hull based continuum approximation.
\end{enumerate}
Given the importance of centrality, we then develop the analysis at lower
densities:
\begin{enumerate}
\item[3)] In Chapter \ref{chapter:geodesicpaths}, in light of the above, we
evaluate the expected number of geodesic paths $\sigma_{r_{ij}}$ which run
between two vertices at displacement $r_{xy}=\norm{x-y}$ in a unit disk
graph. These are the \textit{pair dependencies} which are summed in the
evaluation of betweenness, see section \ref{sec:computation}. Taking the
point process density $\lambda < \infty$ in a domain of arbitrary dimension
$d$, this is approximately equal to a polynomial in $\ceil{r_{xy}}-r_{xy}$ of
order $\frac{1}{2}\floor{r_{xy}}\left(d+1\right)$
\begin{multline}\label{e:final}
	\mathbb{E}\left(\sigma_{r_{xy}}\right) \approx
		\frac{\rho^{\floor{r_{xy}}}
		\left(2\pi\right)^{\frac{1}{2}\floor{r_{xy}}\left(d-1\right)}
\ceil{r_{xy}}^{\frac{1}{2}\left(1-d\right)}}{\Gamma\left(\frac{ \ceil{r_{xy}}
+1}{2} + \frac{\floor{r_{xy}}d}{2}\right)}
\left(\ceil{r_{xy}}-r_{xy}\right)^{\frac{1}{2}\floor{r_{xy}}\left(d+1\right)}
\nonumber
\end{multline}We highlight the difficulties in providing a better
approximation. We also numerically corroborate our formulas, and discuss some
interesting features which appear when non-geodesic paths are incorporated.
\end{enumerate}
\begin{enumerate}
\item[4)] Finally, in Chapter \ref{chapter:discussion}, we discuss potential
applications of this research, and make the relationship between theory and
practice clearer.
\end{enumerate}

\section{Contribution to the field}
The statistical analysis of centrality metrics on random communication
networks can provide useful insights into performance. What is already known
is how betweenness can be algorithmically determined, either by a central
computer, or by the constituent vertices in a centralised way. What is
advanced is how this can ultimately be done analytically. This avoids an
infamously expensive computation (see Section \ref{sec:computation}).

We now list our key contributions:

\begin{enumerate}
\item[a)] In ultra-dense networks, we show how non-convex features of bounded
domains can be highlighted by centrality indices.

\item[b)] $(d-1)$-dimensional subgraphs meandering around urban obstructions
can become critical to optimal performance.

\item[c)] Isolated vertices do not determine connection in
quasi-one-dimensional random geometric graphs.

\item[d)] A simple approximation to network-theoretic betweenness
centrality can be used by various high-layer wireless network routing
protocols in order to omit the expensive and often impossible computation
which would normally occur during operation. We also discuss interesting
applications of analytic formulas.

\item[e)] An approximation to the expected number of geodesic paths joining
two distant vertices will lead to a more complete understanding of the value
of these sorts of statistics in straightforward interference management
techniques applied to 5G small cell deployments. Also, the expected number of
geodesics can be used by ad hoc vertices as part of estimating the distance
to an underground point, where GPS is unavailable. This may develop into
overground localisation where appropriate.
\end{enumerate}
We finally highlight the need to develop the field of spatial probabilistic
combinatorics with a greater focus on applications, in, for example, communication theory, particularly concerning issues such
as the spatial dependence of network observables.
\section{Thesis structure}
The rest of this thesis is structured as follows:

\begin{enumerate}
\item[] \textbf{In Chapter Two} we introduce a number of concepts in
communication theory and random geometric graphs which are relevant to the
arguments in this thesis.

\item[] \textbf{In Chapter Three} we study the effect of non-convexity on the
random connection model.

\item[] \textbf{In Chapter Four} we introduce \textit{betweenness centrality}
in asymptotically dense random geometric graphs.

\item[] \textbf{In Chapter Five} the continuum limit is relaxed, and we
approach betweenness centrality at finite density, finding the expected
number of geodesic paths between two vertices in a unit disk graph (of
general dimension $d$).

\item[] \textbf{In Chapter Six} we discuss the applications of this research to
dense networks, and conclude.
\end{enumerate}

\noindent Also, \textbf{Appendix A} reviews the basic results of percolation, and \textbf{Appendix B} details a proof related to the isolated vertices theorem discussed in Section \ref{sec:connect}.

\chapter{Concepts}\label{c:concepts}

In this chapter, we introduce a number of concepts that will prove important for the work that follows in this thesis. In particular, we introduce the rich field of random networks and their current relation to communication theory. This provides a background and historical context for the our contribution. We first discuss connectivity, isolated vertices and non-convexity in random geometric graphs, and then discuss more recent developments in communication networks which have emerged in network science, specifically centrality indices, which for example measure the structural importance of networked communication nodes. We then briefly discuss current trends in the area.

\section{Connectivity}\label{sec:connect}
One of the early graph-theoretic communication problems was a derivation of 
the probability that a random geometric graph was connected in the
\textit{thermodynamic limit}, where both the number of vertices and the
domain volume $n,V \to \infty$ in such a way that the vertex density is
constant \cite{gupta1998}. The result goes as follows: take a simple Poisson
point process $\mathcal{Y}$ of expected $n$ vertices inside a square of side
$\sqrt{n}$, and form a graph by linking pairs of this process whenever they
are within Euclidean distance $r_0$ of each other. Call this graph
$\mathcal{G}(n,\pi r_0^2)$. According to Penrose \cite{penrose1997}, and
later Gupta and Kumar \cite{gupta1998}, the asymptotic connection probability
of the model with \textit{logarithmic} growth of the connection disks, i.e.
$\mathcal{G}\left(n,\log n + c\left(n\right)\right)$, is given by
\begin{eqnarray}
\mathbb{P}\left(\mathcal{G}\left(n,\log n + c\left(n\right)\right)\text{is
connected}\right) \to e^{-e^{-c\left(n\right)}}
\end{eqnarray}
as $n \to \infty$. Therefore, if $c(n)$ goes to infinity with $n$, the graph
will connect asymptotically almost surely.

To see this, consider a disk of area $\pi r_0^2$ centred at some vertex $x
\in \mathcal{Y}$: this contains no other points of $\mathcal{Y}$ with
probability
$\exp\left(-\pi r_{0}^{2}\right)$. For large $n$, these empty disks, which
are isolated vertices, occur independently in the limit. Therefore the number
of isolated vertices $n_0$ is distributed as a binomial random variable
$\text{bin}\left(n,\exp{\left(-\pi r_0^2\right)}\right)$. A simple Poisson
process $\mathcal{Y}_{n_0} \subset \mathcal{V}$ of isolated vertices will
therefore be observed in the limit. This process is empty with probability:
\begin{eqnarray}\label{e:probabilityofisolation}
\left.\mathbb{P}\left(n_0 = k\right)\right|_{k=0}&=&
\left.\frac{1}{k!}e^{-ne^{-\pi r_0^2}}\left(ne^{-\pi
r_0^2}\right)^{k}\right|_{k=0} \nonumber \\
&=& \left.\frac{1}{k!}e^{-ne^{-\left(\log n +
c(n)\right)}}\left(ne^{-\left(\log n + c(n)\right)}\right)^{k}\right|_{k=0}
\nonumber \\
&=& \left.\frac{1}{k!}e^{-e^{-c(n)}}\left(e^{-c(n)}\right)^{k}\right|_{k=0}
\nonumber \\
&=& e^{-e^{-c\left(n\right)}}
\end{eqnarray}
The question is, does this lack of isolated vertices imply the graph is
connected?
\subsection{Isolated vertices}\label{sec:isolatedvertices}
Consider the following two lemmas \cite{penrose1997}:
\begin{lemma}[No two components are large]
Assuming $c > 0$, there exists a $C$ such that asymptotically almost surely
the random graph $G\left(n, c \log n\right)$ does not consist of two or more
connected components each with \textit{Euclidean diameter}\footnote{This is
the largest Euclidean distance which can be found between any two vertices in
a component.} at least $C\sqrt{\log n}$ .
	\end{lemma}
\noindent This means that two large enough components will merge
asymptotically almost surely (a.a.s.). But how large?
\begin{lemma}[All small components consist of a single
vertex.]\label{l:theyareisolatedvertices}
The graph $G\left(n, \log n - \frac{1}{2}\log \log n \right)$ contains no
components $H$ of more than one vertex and Euclidean diameter strictly less
than $C\sqrt{\log n}$.
\end{lemma}
\noindent These two lemmas imply that there exists a specific phase of the
sub-logarithmic growth of $r_{0}$ where a single giant component forms
surrounded by isolated vertices\footnote{Though we provide a summary of the proof in this section, we provide its detail in Appendix \ref{appendix:isolatedvertices}.}. We have just shown that
these isolated vertices form a Poisson point process of their own. This means
that one can effectively approximate the connection probability $P_{fc}$ as
\begin{eqnarray}
P_{fc} \sim 1 - \mathbb{P}\left(n_{0}=0\right)
\end{eqnarray}

To see the first of these lemmas, tile the square domain with tiles of side
$r_0/\sqrt{20}$. This ensures that any two vertices found in any two adjacent
squares are no more than $r_0/2$ apart. Then argue that
\begin{enumerate}
\item A component $U$ of Euclidean diameter at least $C\sqrt{\log n}$ covers
many tiles as $n \to \infty$.
\item Since the tiles have side $r_0/\sqrt{20}$, all tiles adjacent to $U$
must be empty while the component exists.
	\item There are many empty boundary tiles, given the size of $U$.
	\item As $n \to \infty$, they cannot all be empty.
\end{enumerate}
Thus any component of diameter at least $C\sqrt{\log n}$ will merge with
another component, leaving only small components. We do not explicitly prove
parts $1$-$3$, but refer directly to Walter's review \cite{waltersreview}; it
essentially suffices to count the tiles.

Now, the second part is based on the same sort of argument: that an empty
area must be maintained around a small component in order to keep it small,
and that this area is, asymptotically, filled, unless it is a single vertex.
Thus all large components join together, all small components are isolated
vertices, and $P_{fc} \sim \mathbb{P}\left(n_{0}=0\right)$.
\subsection{The connection probability}\label{sec:connectionprobability}
Given this relation between isolated vertices and connectivity, it suffices
to evaluate the probability that a single isolated vertex exists in the
random connection model, then suggest that as the density $\rho \to \infty$,
\begin{eqnarray}
P_{fc} \sim 1 - \mathbb{P}\left(n_{0}=1\right)
\end{eqnarray}
We now evaluate this. Consider a vertex at some fixed $x \in \mathcal{V}$.
Its degree $k(x)$ can be determined by looking at a \textit{marked} Poisson
point process $\mathcal{Y}^{\star}$, where the marks are $U\left[0,1\right]$
random variables
\begin{equation}
\mathcal{Y}^{\star}=\left\{\zeta,u:\zeta \in \mathcal{Y},u \sim
U\left[0,1\right]\right\}
\end{equation} 
which is of intensity $\rho \textrm{d}x$ on $\mathcal{V} \times
U\left[0,1\right]$, and $\textrm{d}x$ is Lesbegue measure on
$\mathbb{R}^{d+1}$. The degree is given by a sum over this marked point
process:
\begin{equation}\label{e:neweq2}
k\left(x\right) = \sum_{\left(y,u\right) \in \mathcal{V} \times
U\left[0,1\right]} \mathbf{1}\{u < \chi \left(x,y\right)
H\left(\norm{x-y}\right)\}
\end{equation}
According to Campbell's theorem, one of the elementary theorems about point
processes \cite{kingmanbook}, we have that $k(x)$ is Poisson with expectation
\begin{equation}\label{e:neweq3}
\mathbb{E}k(x) = \rho \int_{\mathcal{V}} \chi \left(x,y\right)
H\left(r_{xy}\right) \mathrm{d}y
\end{equation}
and therefore
\begin{eqnarray}\label{e:neweq4}
\mathbb{P}\left(k\left(x\right)=0\right)= \exp\left(-\rho \int_{\mathcal{V}}
\chi \left(x,y\right) H\left(r_{xy}\right) \mathrm{d}y\right)
\end{eqnarray}
\noindent which gives the probability that a vertex find itself isolated.

Now, given the conjecture discussed in \ref{sec:connect} (that isolated
vertices occur as a simple Poisson point process), it is natural in light of
Eq. \ref{e:neweq4} to
conjecture that as $\rho \to \infty$, the total number of isolated vertices
is well approximated by a Poisson distribution with mean
\begin{equation}\label{e:neweq5}
\rho \int_{\mathcal{V}} \exp \left(-\rho \int_{\mathcal{V}} \chi
\left(x,y\right) H\left(r_{xy}\right) \mathrm{d}y \right) \mathrm{d}x
\end{equation}

In this limit, as discussed, the obstacle to connectivity is the presence of
isolated vertices, and so, as $\rho \to \infty$,
\begin{align}\label{e:connectprob}					
P_{fc} \sim \exp \left(-\rho \int_{\mathcal{V}} e^{-\rho \int_{\mathcal{V}}
\chi \left(x,y\right) H\left(r_{xy}\right) \mathrm{d}y } \mathrm{d}x \right)
\nonumber
\end{align}
which is approximately $$1 - \rho \int_{\mathcal{V}} e^{-\rho
\int_{\mathcal{V}} \chi \left(x,y\right) H\left(r_{xy}\right) \mathrm{d}y }
\mathrm{d}x$$ for large $\rho$. One can therefore evaluate the connection probability of the
random connection model with an integral.
\subsection{Boundary effects and non-convex domains} 
We now discuss the introduction of non-convexity. This occurs when at least
one straight line segment with endpoints in the domain intersects the
domain's complement. We use this function in order to highlight this:
\begin{equation}\label{e:neweq1}
	\chi \left(x,y\right) =
  \begin{cases}
1 & \quad \text{if } x + \lambda\left(y-x\right) \in \mathcal{V} \text{ for
all } \lambda \in \left[0,1\right] \\
	  0  & \quad \text{otherwise} \\
  \end{cases}
\end{equation}
The perimeter itself, which may be the outside of a building, arena or city,
plays a significant role in connectivity. In \cite{cef2012}, the connectivity
of the random connection model was shown to be influenced by the meandering
boundary's ability to block the viewing angle vertices, decreasing it from
the usual $2 \pi$ radians. In fact, pressed right up on the boundary of a
disk, vertices see less than $\pi$ radians. The extent to which this results
in vertices becoming isolated is quantified by the expected degree of a point
$x$, which is
\begin{eqnarray}\label{e:cmass}
\mathbb{E}k(x) = \int_{\mathcal{V}(x)}H(\norm{x-y})\mathrm{d}y
\end{eqnarray}
where $\mathcal{V}(x)$ is the region of the domain visible to $x$, itself a
domain\footnote{This is often centred at $x$, to make the integral more
tractable.}. This is often called the \textit{connectivity mass} of a point
$x \in \mathcal{V}$.

Thought of as a functional, Eq. \ref{e:cmass} is minimised in the domain's
sharpest corner \cite{cef2012,orestis2013}. In non-convex domains, the
internal perimeter, part of which may be the boundary of an obstacle,
provides another point where the connectivity mass is lower than the bulk. We
investigate this in chapter \ref{chapter:connectivity}.
\section{Betweenness centrality}\label{sec:betweennessintro}
\begin{figure}[t]
\begin{centering}
\includegraphics[scale=0.65]{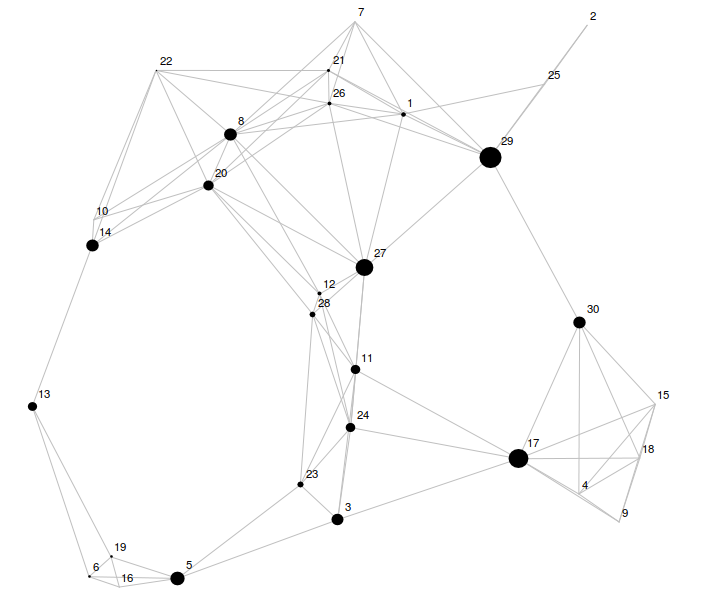}
\caption{Betweenness centrality in a unit disk graph of 30 vertices drawn inside the unit square, with $r_0=0.3$. The size of the vertices is proportional to their betweenness.}
\end{centering}
\label{fig:betweennessexample}
\end{figure}

Shortest path betweenness counts all \textit{shortest} paths between all
pairs of vertices in a graph. Taking $\Sigma_{r_{xy}}$ as the set of all paths between
two vertices $i,j$ (at Euclidean displacement $r_{xy}$), and $\Sigma^{\star}
\subseteq \Sigma_{r_{xy}}$ as that proportion of paths of length, in hops,
shorter than or equal to any other path in $\Sigma_{r_{xy}}$, then
\begin{eqnarray}
\sigma_{r_{xy}} \coloneqq \norm{\Sigma^{\star}}
\end{eqnarray}
where the norm here means the cardinality of the set. In analogy with
differential geometry, we call these paths \textit{geodesic}.

Some of these geodesics will potentially run through some specific vertex $z
\in \mathcal{V}$. This is counted as $\sigma_{r_{xy}}(z)$. For example, if
all geodesic paths from $x$ to $y$ pass through $z$, then the ratio
$\sigma_{r_{xy}}(z)/\sigma_{r_{xy}}$ will equal unity. Otherwise, it will be
a real number between zero and one. To quantify the extent to which a vertex
lies on many geodesic paths, this ratio is summed over all vertex pairs in
$\mathcal{G}$ to produce the \textit{shortest path node betweenness centrality} of $z$
\begin{eqnarray}\label{e:betweennessintro}
\sum_{i \neq j, j \neq k, k \neq i}\frac{\sigma_{r_{xy}}(z)}{\sigma_{r_{xy}}}
\end{eqnarray}
The terms in this sum are called
\textit{pair dependencies} \cite{brandes2001}.

\subsection{Computation}\label{sec:computation}
In order to exploit the sparsity of geodesics, a traversal algorithm is used
to evaluate betweenness, due to Brandes \cite{brandes2001}. Two steps are
required. The first, just discussed, is:
\begin{itemize}
\item[1)] Calculate $\sigma_{xy}$ for all $x,y \in \mathcal{V}$.
\end{itemize}
The second involves the evaluation of the pair dependencies in Eq.
\ref{e:betweennessintro}:
\begin{itemize}
\item[2)] Calculate $\sigma_{xy}(z)$ for all $x,y \in \mathcal{V}$, and sum
to find $\sum_{i \neq j \neq k}\frac{\sigma_{r_{xy}}(z)}{\sigma_{r_{xy}}}$.
\end{itemize}
For a single pair, the first step is at worst
$\mathcal{O}\left(\singlenorm{E}\right)$. A modified version of the
Floyd-Warshall algorithm is used. Brandes incorporates the second step into
the first by recursively updating the pair dependencies during the various
evaluations of $\sigma_{xy}$. A vector of betweenness can then be obtained in
at worst $\mathcal{O}\left(\singlenorm{V}\singlenorm{E}\right)$.

As networks become large, this computation is notoriously intensive. It also
relies on co-operation between all networked devices, which may prove
impossible in practice. Though, interesting random sampling methods are known
to speed up computation, see e.g. \cite{brandes2007}.
\subsection{Shortest paths}
\begin{figure}[t]
\hspace{12mm}\includegraphics[scale=0.18]{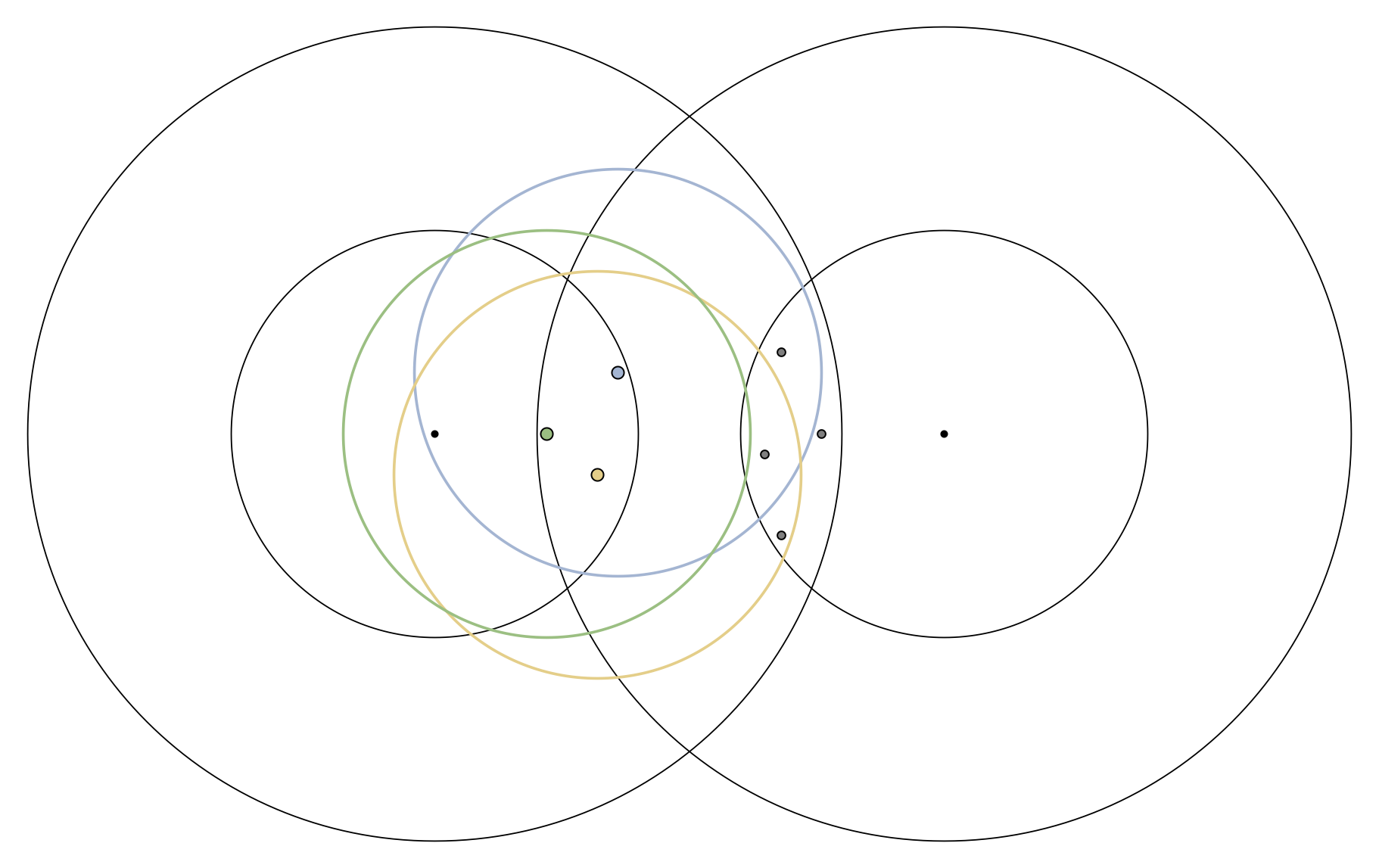}
\caption{Two vertices at center separation 2.5 in a unit disk graph, taking $r_0=1$. The connection ranges, and an extra range (large disks), are drawn. The degrees of the three coloured vertices are dependent random variables. In this example, the green vertex has degree zero, the blue vertex has degree two and the yellow vertex also has degree two. This spatial dependence of vertex degrees in random geometric graphs such as the unit disk graph depicted here is ignored in various independence assumptions, as discussed in the text.}\label{fig:ia}
\end{figure}
Part of the problem of analytically evaluating betweenness centrality
involves the difficulty of enumerating the number of geodesic paths which run
between two vertices in a random geometric graph. In the Erd\H{o}s-R\'{e}nyi
case\footnote{Take $n$ vertices and link them independently with probability
$p$.} \cite{erdos1959}, where the connection probability is just some
constant function, this is relatively simple: one recursively calculates the
probability that some vertex is one hop to its destination, which is $p$, two
hops, which is the probability two vertices form a link to the same vertex
and don't connect directly, which is $(1-p)(n-2)p^2$, and so on.     When the connection function is not constant, however,
nearby vertices have positively correlated vertex degrees. This complex
spatial dependence is characteristic of combinatorial problems with a
geometric element.

Now, in order to avoid this, one can make the so called \textit{independence
assumption} \cite{zhang2012,mao2010}. This has been used to produce
approximations to the distribution of the length of paths between vertices in
the unit disk model \cite{ta2007}, the random connection model
\cite{zhang2012}, and the log-normal fading model \cite{mukherjee2008}. The
issue is depicted in Fig. \ref{fig:ia}.

The \textit{number} of geodesics is an alternative statistic, related, as
discussed, to betweenness centrality. This random \textit{number} of
geodesics is currently used by two important delay tolerant networking
algorithms, BubbleRap and COAR, which algorithmically count geodesics in
order to assess possible routing strategies \cite{magaia2015}.
\section{Extended Concepts}\label{sec:pointprocesses}
We now discuss some extended topics related to point processes, alternative
metrics, and centrality variants.
\subsection{Binomial point processes}\label{sec:bpp}
If we simply fix the number of vertices in a Poisson point process, we have a
binomial point process. This can be a more realistic model of wireless
networks, since the device numbers cannot vary so much over time (in fact,
they don't vary at all).

However, the Poisson point process is used with good reason: with exactly $N$
vertices, the modelling can become intractable, since the distance between
points in a (BPP) formed within a domain $\mathcal{V} \subset \mathbb{R}^d$
is given by the beta distribution, which is only expressible in terms of a
special function. Even worse, sequences of inter-point distances are no
longer independent.

In \cite{srinivasa2010}, potential use of the binomial point process as a
model of device distribution in wireless networks is discussed in detail.

\subsection{Markov point processes}\label{sec:mpp}
Consider a Markov chain $\mathcal{M}$ on the space of unit disk
graphs.\begin{figure}
\begin{centering}
	\includegraphics[scale=0.5]{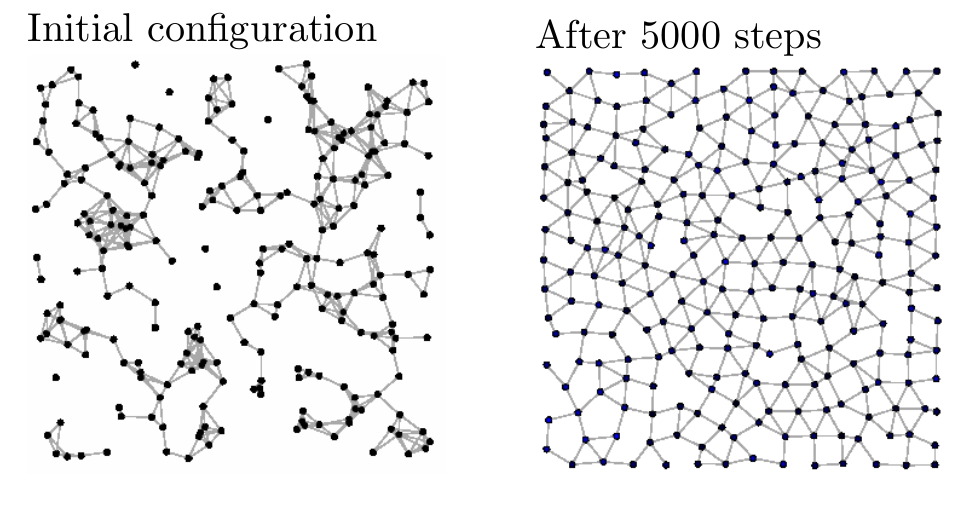}
\par\end{centering}
\caption{The Strauss process. A binomial point process of $n=250$ vertices
(taking $\Omega=0.07$ and $r_0=0.08$) is updated in a step by step fashion by
the Metropolis-Hastings algorithm. It creates a pattern of repelling points.
From Irons et al 2011, referenced in the text.}\label{fig:jordanpic}
\end{figure}
A simple example is the Strauss process. Pairs of points $\{ \zeta,\eta \}
\subset \mathcal{Y}$ `interact' with a geometric potential $\phi$, given by
their Euclidean separation
\begin{eqnarray}
\log \phi \left( \{ \zeta,\eta \} \right) \propto \frac{1}{\norm{\zeta-\eta}}
\end{eqnarray}
The probability of a configuration, $\mathbb{P}(x_1,x_2,\dots,x_N)$, is then
given by a sum over the set $\mathcal{C}$ of all $2$-vertex cliques in the
respective unit disk graph formed on the configuration
\begin{eqnarray}
\log \mathbb{P}(x_1,x_2,\dots,x_N) \propto \sum_{\{ \zeta , \eta \} \in
\mathcal{C}} \log \phi \left( \{ \zeta,\eta \} \right)
\end{eqnarray}
which is associated with a state of $\mathcal{M}$. 

One can generate these configurations with the Metropolis-Hastings algorithm
\cite{irons2011}. The algorithm goes as follows:
\begin{enumerate}
\item Distribute $N$ points over a domain $\mathcal{V} \subseteq
\mathbb{R}^d$.

\item Pick a vertex (call its position $v$), and, picking an angle uniformly,
randomly displace it a distance given by the Beta distribution\footnote{Which
gives the distance between points of a binomial point process
\cite{srinivasa2010}.}, rejecting moves which put the vertex outside the
square. Call this new location $\mu$.

\item Calculate two quantities, first 
\begin{eqnarray}
n_v = \sum_{i \neq v} I_{\{\norm{i - v}<\Omega\}} \frac{\Omega}{\norm{i-v}}
\end{eqnarray}
with $\Omega$ a parameter of the process with units of distance, and then
\begin{eqnarray}
n_\mu = \sum_{i \neq \mu} I_{\{\norm{i - \mu}<\Omega\}}
\frac{\Omega}{\norm{i-\mu}}
\end{eqnarray}
The move is then accepted with probability $\min(1,\omega^{n_{\mu}-n_{v}})$,
where $\omega \in \left[0,1\right]$ is some parameter of the model
quantifying the amount of inter-point repulsion.
\end{enumerate}

The fact that the potential of the whole configuration can be written as a
product of all two vertex cliques in the graph implies that the density
$\mathbb{P}(x_1,x_2,\dots,x_N)$ is a Gibbs ensemble:
\begin{eqnarray}
\mathbb{P}(x_1,x_2,\dots,x_N) \propto \prod_{\{ \zeta , \eta \} \in
\mathcal{C}} \phi \left( \{ \zeta,\eta \} \right)
\end{eqnarray}
which is a factorisation over cliques.

This sort of repelling points model helps combat the most unrealistic
assumption of the Poisson process model, that there will sometimes be nearby
base stations. Determinantal point processes, which are processes whose
spatial distributions are related to the determinant of a matrix, are an
interesting avenue of further research.
\subsection{Geometric preferential attachment}\label{sec:gpa}
Considers a binomial point process $\mathcal{Y}$ of $n$ points
$x_1,x_2,\dots,x_n$ on the surface of a torus $\mathbb{T}$. Allowing multiple
edges, each point is selected in turn such that, at time $t \in
\left\{1,2,\dots,n\right\}$, vertex $x_t \in \Phi$ forms $m$ randomly
selected connections to those vertices within Euclidean distance $r_0$; for
each $i \in \left\{1,2,\dots,m\right\}$, the probability vertex $v$ is
selected is
\begin{eqnarray}
P\left(v \text{ is selected}\right) =
\frac{\text{deg}_{t}v}{\max\left(\sum_{\norm{x_{t}-v} \leq
r_0}\text{deg}_{t}v,\gamma\right)}
\end{eqnarray}
with
\begin{eqnarray}
P\left(x_t \text{ connects to itself}\right) = 1 - \frac{\sum_{\norm{x_{t}-v}
\leq r_0}\text{deg}_{t}v}{\max\left(\sum_{\norm{x_{t}-v} \leq
r_0}\text{deg}_{t}v,\gamma\right)}
\end{eqnarray}
The factor $\gamma = \alpha m \pi r_0^2 t$ tunes the probability of loop
formation (i.e. self connection). This propensity for loop formation
increases with time, the number of edges $m$ we intend to add, and the
connection radius $r_0$.

This model is an extension of the \textit{online nearest neighbour graph}
(see e.g. \cite{penrose2008}, or the earlier\footnote{This was published in
2007, a year before \cite{penrose2008}.}paper of Berger, Bollob\'{a}s, Borgs,
Chayes and Riordan \cite{berger2007}), which is a simple growing spatial
network: vertices are placed one by one in some domain, and each time joined
to their nearest neighbour \cite{penrosebook}.

\subsection{Hyperbolic random geometric graphs}\label{sec:hyperbolic}

Apart from altering the nature of the point process as in subsection
\ref{sec:pointprocesses}, one can transform the underlying metric space in
which the graphs are embedded. An interesting example of this is the
\textit{hyperbolic random geometric graph} \cite{krioukov2010}, where a
simple Poisson point process is formed on the hyperbolic plane with points retained inside a bounded region. Pairs of points
are then joined according to the unit disk rule, but now considering the
distance $g$ between two points $\left(x_1,y_1\right)$ and
$\left(x_2,y_2\right)$ to be
\begin{eqnarray*}
g\left(\left(x_1,y_1\right),\left(x_2,y_2\right)\right) =
\arccosh\left(\cosh(y_1)\cosh(x_2-x_1)\cosh(y_2)-\sinh(y_1)\sinh(y_2)\right)
\end{eqnarray*}
where the Cartesian axis are embedded in the hyperbolic plane, and the two
arguments of $g$ read off accordingly.

As a current application area, these non-standard RGG’s are used to study theories concerning \textit{complex systems} \cite{newmanbook,ladyman2013}.

\subsection{Current flow betweenness}\label{sec:flow}
\begin{figure}[!]
\begin{center}
\hspace{15mm}\includegraphics[scale=0.28]{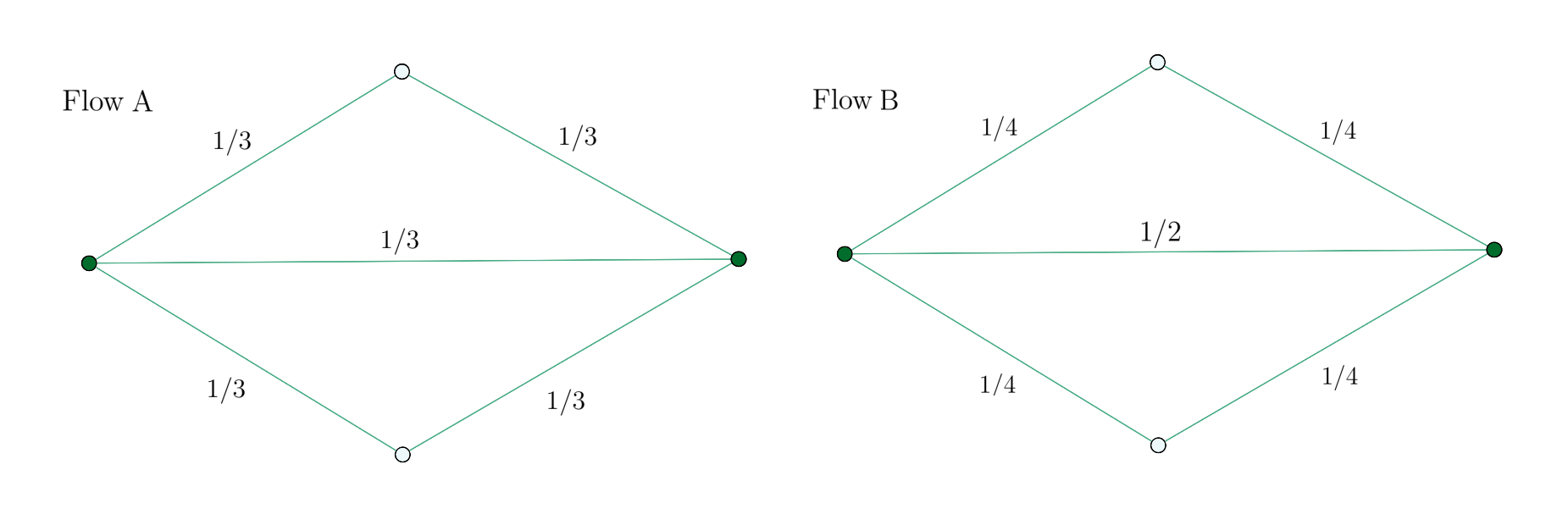}
\end{center}
\caption{A current of 1$A$ flows into the left most vertex of a network of
unit resistors. It is then extracted from the right most vertex. Two
different current conserving flows are shown on the same simple network. The
values on the edges are the currents flowing along those edges. Flow $B$
minimises the total dissipated energy, and so satisfies the Kirchhoff
laws.}\label{fig:flowpic}
\end{figure} Two alternative forms of betweenness, the first introduced by Freeman
\cite{freeman1991} called \textit{flow betweenness}, and the second by Newman
\cite{newman2003} called \textit{current flow betweenness}, aim to quantify
the effect of non-geodesic paths on the centrality of vertices:
\begin{enumerate}
\item Flow betweenness considers a function $f: V \to \mathbb{R}$ quantifying
the current passing through some vertex $v \in V$ while the flow through the
network is maximal. The capacity of all edges must be defined, and one sums
of all possible source-destination pairs.

\item Current flow betweenness considers the unique `unit flow', which is a
current of 1$A$ (one Ampere), passing through the network while satisfying the Kirchhoff
laws. Note there is only one flow which minimises the dissipated energy, and
hence only one flow which satisfies Kirchhoff's laws \cite{grimmettbook}. Again, one sums this over all
source-destination pairs.
\end{enumerate}

Note, dissipated energy is the product of current and voltage. Equivalently
\begin{eqnarray}
E = \sum_{e \in E}I_e^2R_e
\end{eqnarray}
where $I_e$ is the current flowing over edge $e$, and $R_e$ is the resistance
of edge $e$. Minimising this using Lagrange multipliers with the Kirchhoff
current law as a constraint \cite{doyle1984}, see for example Fig.
\ref{fig:flowpic}, provides currents flowing through each vertex with respect to two
other vertices $v_i,v_j \in V$, the sum of which, over all $i,j$, or `source
sets', is the current flow betweenness.

For example, in Fig. \ref{fig:flowpic}, the dissipated energy for flow $A$ is
$5/9$, while for flow $B$, it is $1/2$. Therefore, given we consider a unit flow,
flow $B$ is the unique solution to the Kirchhoff laws.
This would contribute to the vector of current flow betweenness centralities
as one of the many source-sink pair contributions.

Also, note that current flow betweenness is equivalent to random walk
betweenness \cite{newman2003}. Where data is flooded through a network of
small cells, or through a vehicle network, this current flow betweenness can
be useful in estimating the routing load on set of vertices when betweenness
is shown to be uncharacteristic i.e. when data is often diverted from
shortest paths. It is difficult to analyse analytically, and may be more
appropriately obtained algorithmically in realistic settings
\cite{lulli2015}.
\chapter{Connectivity}\label{chapter:connectivity}

\section{Introduction} 
Soft random geometric graphs consist of a set of vertices placed according to
a point process in some domain $\mathcal{{V}}\subseteq\mathbb{R}^{d}$ which
are coupled with a probability dependent on their Euclidean separation. The more common
deterministic connection is generalised to probabilistic connection \cite{penrose2016,cef2012,mao2011} in order, in our case, to model
signal fading. Commonly known as the \textit{random connection model}, we now
have a connection function $H\left(\norm{x-y}\right)$ giving the probability
that links will form between nodes $x,y \in \mathcal{V}$ of a certain
Euclidean displacement $\norm{x-y}$. In a band-limited world of wireless communications continuously pressed
for the theoretical advances that can enable 5G cellular performance, this is
an important new flexibility in the model.

Connectivity has been shown the initial interest \cite{cef2012,orestis2013, mao2011}.
For example, in \cite{cef2012}, using a cluster expansion technique from
statistical physics, at high vertex density $\rho$
the connection probability of a soft random geometric graph formed within a
bounded domain
$\mathcal{{V}}$ is approximated as (the complement of) the
probability that exactly one isolated vertex appears in an otherwise
connected graph.
This is justified by a conjecture of Penrose \cite{penrose2016}, asserting
that the number of isolated vertices follows a Poisson distribution whose
mean
quickly decays as $\rho\rightarrow\infty$, thus highlighting the
impact of the domain's enclosing boundary \cite{cef2012, mao2011,kog105}
where isolation is most common.

Internal boundaries, such as obstacles, cause similar problems. In this
chapter, we focus our efforts on how this particular aspect of the domain
effects the graph behaviour. We therefore extend recent work on connectivity
within non-convex domains, such as those incorporating internal walls
\cite{orestis2013} or a complex, fractal
boundary \cite{dettmann2014}, deriving analytic formulas for the connection
probability
$P_{fc}$ of soft random geometric graphs formed within the annulus and
spherical shell geometries,
quantifying how simple convex obstacles affect connectivity. 
Specifically, we consider the situation where nodes connect with
a probability decaying \textit{exponentially} with their mutual Euclidean
separation. This models the \textit{Rayleigh fading} commonly observed in
mobile communications.\footnote{Most of this chapter has been accepted with
minor revisions for publication at \textit{The Journal of Statistical
Physics} \cite{giles2016}. The work in this chapter remains solely that of the author of
this thesis unless otherwise indicated, given collaboration with supervisors
in the appropriate fashion.}.

\begin{figure}
\noindent \begin{center}
\includegraphics[scale=0.18]{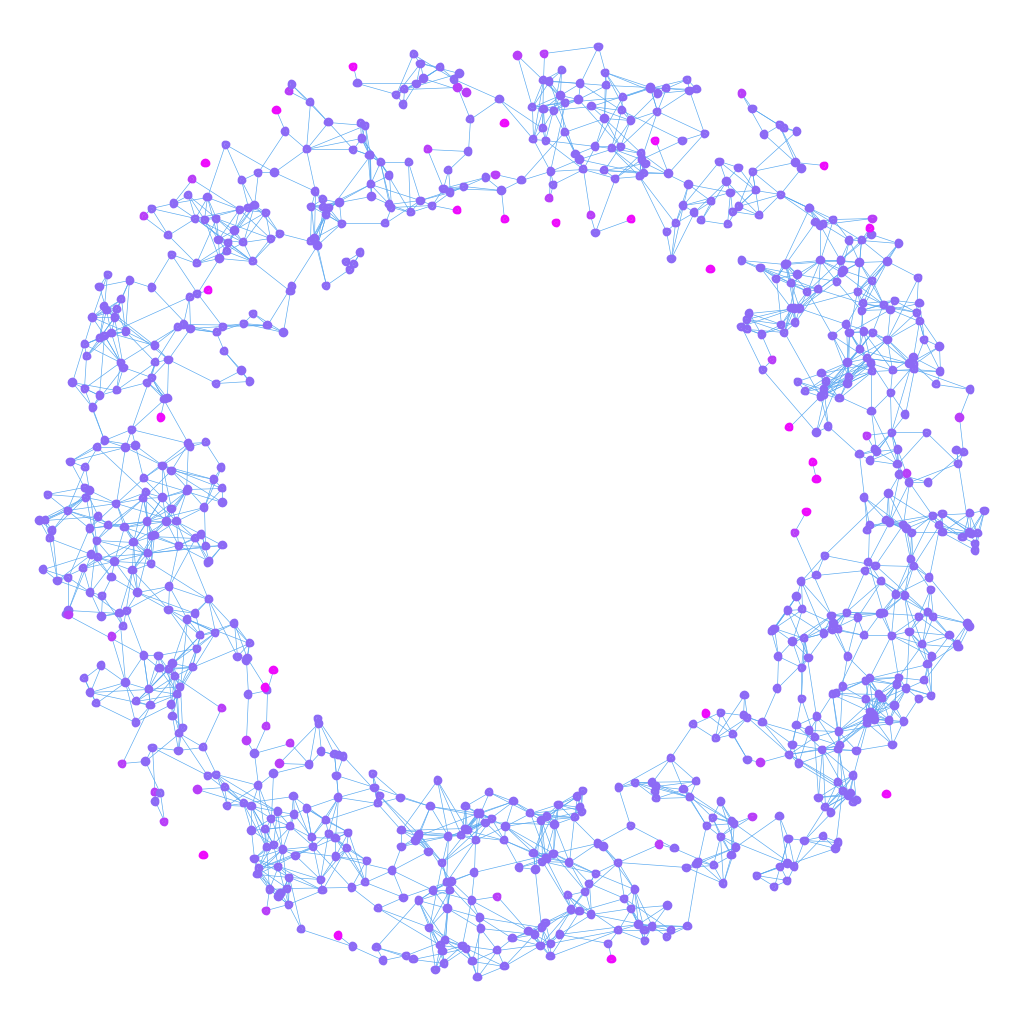}
\includegraphics[scale=0.18]{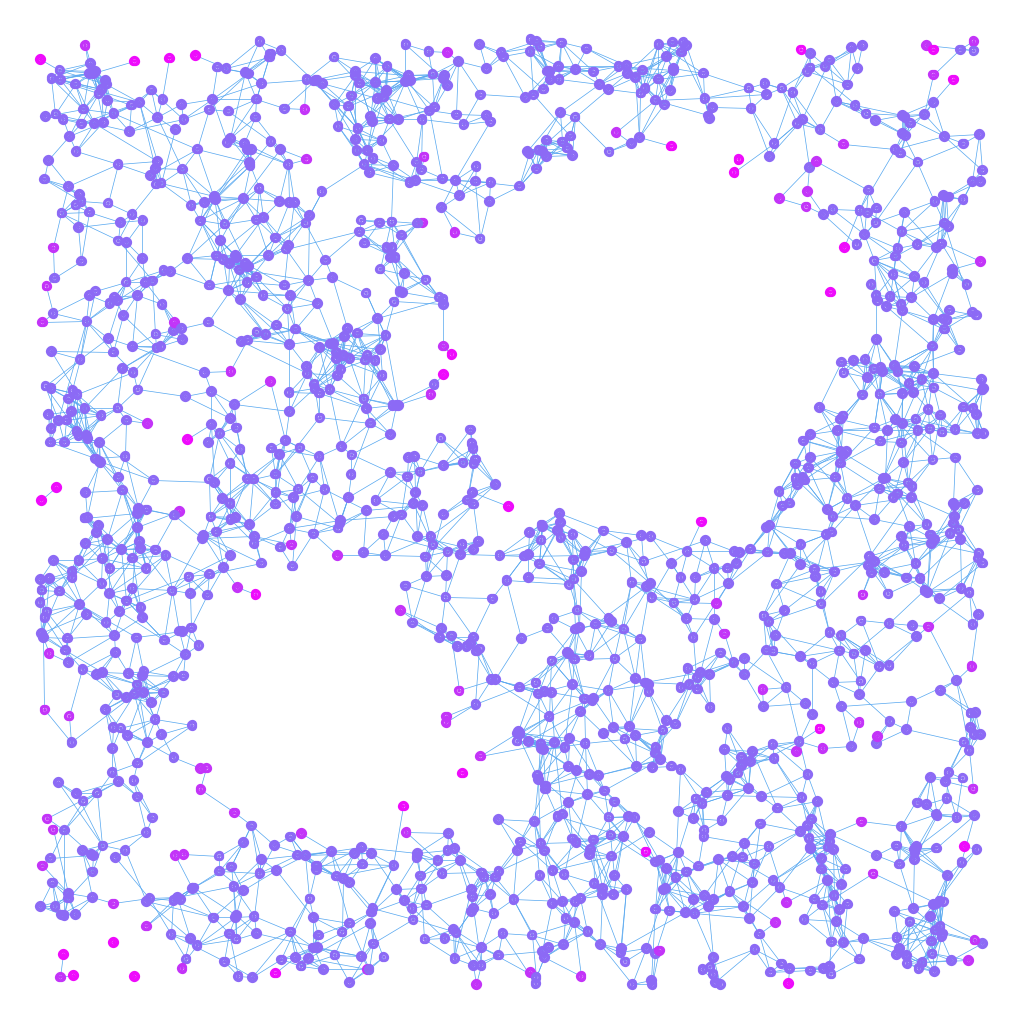}
\par\end{center}
\caption{A soft random geometric graph, first defined at the beginning of Chapter \ref{c:intro}, inside the annulus (large obstacle
case), and inside a square with two circular obstacles. Vertices with low
degree are highlighted in purple. We derive the graph connection probability
in these simple obstructed domains.}
\label{fig:annulus}
\end{figure}

The advantage of soft graphs is this ability to incorporate a fading model.
In built up urban environments where signals scatter repeatedly off walls,
and offering no direct line of sight between transmitter and receiver, the
relevant fading statistics are, as mentioned, those of the Rayleigh
distribution, implying the connection function
\begin{eqnarray}
H(\norm{x-y}) = \exp\left(-\beta \norm{x-y}^2\right)
\end{eqnarray} is to be implemented. We now discuss how this comes about. 

\section{Rayleigh fading}\label{sec:rayleigh}
In Rayliegh fading, the channel \textit{impulse response} \cite{tsebook} is
modelled as a \textit{complex Gaussian process} $\xi$, a sequence of
complex-valued random variables
\begin{equation}
\xi =:\{z_n|n \in \mathbb{N}\}, \quad z_n = U_n(t) + iV_n(t)
\end{equation}
with $U_n\left(t\right)$ and $V_n\left(t\right)$ Gaussian for all $n$,
independent of both each other, and the remaining elements of $\xi$, for all
$n$. One can think of this as a sequence of vectors, each real component of
which is the randomly attenuated amplitude of the quadrature and in-phase
components of a frequency modulated radio signal.

The amplitudes of the impulse response are therefore given by a Pythagorean
relation between these Gaussianly distributed quadrature and in-phase
amplitudes:
\begin{equation}
	r = \sqrt{\norm{U}^2 + \norm{V}^2}
\end{equation}
i.e. the radial component of the sum of two independent Gaussian variables.
The amplitudes $r$ are thus \textit{Rayleigh distributed}. To see this,
integrate the joint density of $U$ and $V$ along the perimeter of a disk of
radius $\lambda$
\begin{eqnarray}\label{e:neweq7}
P\left(r = \lambda \right) &=& \frac{1}{2\pi\sigma^2} \int_{-\infty}^\infty
du \, \int_{-\infty}^\infty dv \, f_U \left(u\right) f_V \left(v\right)
\delta(\lambda-\sqrt{u^2+v^2}) \nonumber \\
&=& \frac{1}{2\pi\sigma^2} \int_{-\infty}^\infty du \, \int_{-\infty}^\infty
dv \, e^{-u^2/2\sigma^2} e^{-v^2/2\sigma^2} \delta(\lambda-\sqrt{u^2+v^2})
\nonumber \\
				 &=&  \frac{\lambda}{\sigma^2} e^{-\lambda^{2}/2\sigma^2}
\end{eqnarray}
which is the density of the Rayleigh distribution.

Now, conisder the so called outage probability $P_{\text{out}}$, which is the
proportion of time the information-theoretic decoding error at the receiver
falls below a critical rate $\Upsilon$:
\begin{eqnarray}\label{e:neweq8}
P_{out} & = &
P\left[\log_{2}\left(1+\frac{P}{N_{0}}\hspace{0.5mm}\norm{h^{2}}\right)<\Upsilon\right]
\end{eqnarray}
We are thus interested in the random channel gain $\norm{h}^2$, and the
signal power $P$. These power gains are proportional to the \textit{square}
of the impulse response's Rayleigh distributed amplitudes, so are
exponentially distributed:
\begin{eqnarray}\label{e:neweq7}
P\left(\norm{h}^2 = \lambda \right) &=&
\frac{1}{2\sigma^2}\exp\left(\frac{-\lambda}{2\sigma^2}\right)
\end{eqnarray}
Also, the signal suffers from a propagation decay. Writing the signal power
as $P$, and the noise power as $N_0$, the signal-to-noise ratio decays as a
power $\eta$ of the transmitter-receiver propagation distance $\norm{x-y}$
\begin{eqnarray}
 \frac{P}{N_0} = c \norm{x-y}^{-\eta}
\end{eqnarray}
This is the path loss exponent, with $\eta =2$ related to free-space
propagation. Rearranging \ref{e:neweq8} gives the source-destination
connection probability as the complement of the outage probability
\begin{eqnarray}\label{e:neweq10}
H\left(\norm{x-y}\right) &=& 1 - P\left[ \norm{h}^2 < \frac{ N_{0}
\left(2^{\Upsilon}-1\right)}{P}\right]
\end{eqnarray}
Finally, extracting the constant $\beta$:
\begin{eqnarray}
\frac{ N_{0} \left(2^{\Upsilon}-1\right)}{P} &=&
\frac{N_{0}\norm{x-y}^{\eta}}{c}\left(2^{\Upsilon}-1\right) \nonumber \\
&=& \left(\frac{N_{0}}{c} \left(2^{\Upsilon}-1\right)\right)
\norm{x-y}^{\eta} \nonumber \\
&=&\beta \norm{x-y}^{\eta} \nonumber
\end{eqnarray}
This implies the Rayleigh fading connection function is given by
\begin{eqnarray}
H(r_{xy}) &=& 1 - P\left[ \norm{h}^2 < \beta \norm{x-y}^{\eta} \right]
\nonumber \\
&=& \exp\left({-\beta \norm{x-y}^{\eta}}\right)
\end{eqnarray}												 
since $\norm{h}^2$ is exponentially distributed, as discussed.\\

To clarify notation, we sometimes refer to the constant
\begin{eqnarray}
	r_{0} = \beta^{-1/\eta}
\end{eqnarray}
to signify the length scale over which nodes typically connect, since the
exponent $\beta \norm{x-y}^\eta>1$ whenever $\norm{x-y}>r_0$, and so the
connection probability is low.

\section{The annulus domain $\mathcal{{A}}$}\label{sec:three}
Take our domain to be the annulus $\mathcal{{A}}$ of inner radius $r$ and
outer radius $R$, two examples of which are depicted in Fig. \ref{fig:trip}.
Consider the outer radius of this annulus to be large compared to the typical
connection range. We are interested in evaluating
\begin{eqnarray}\label{e:connectprob2}					
\mathbb{P}\left(n_{0}=0\right) \approx 1 - \rho \int_{\mathcal{V}} e^{-\rho
\int_{\mathcal{V}} \chi \left(x,y\right) H\left(r_{xy}\right) \mathrm{d}y }
\mathrm{d}x
\end{eqnarray}
as discussed in Section \ref{sec:connectionprobability}. There, we defined
the connectivity mass at a point $x \in \mathcal{A}$, and its analogue over
the region visible to $x$
\begin{eqnarray}
\mathcal{M}\left(x\right) &=& \int_{\mathcal{A}\left(x\right)}
H\left(\norm{x-y}\right)\textrm{d}y
\end{eqnarray}
This mass is approximated within two obstacle-size regimes, the first
where $r\ll r_{0}$, and the second where $r\gg r_{0}$. In each regime we can
make some assumptions about the geometry of the region
$\mathcal{A}\left(x\right)$ visible to $x$, which yields tractable formulas
for the connectivity mass in terms of
powers of the distance $\epsilon$ from the obstacle's perimeter. We then have
$P_{fc}$ in the annulus $\mathcal{A}$. This complements the result of the
disk domain, presented first in \cite{cef2012}.

\subsection{No obstacles}
\begin{figure}
\noindent \begin{centering}
\includegraphics[scale=0.35]{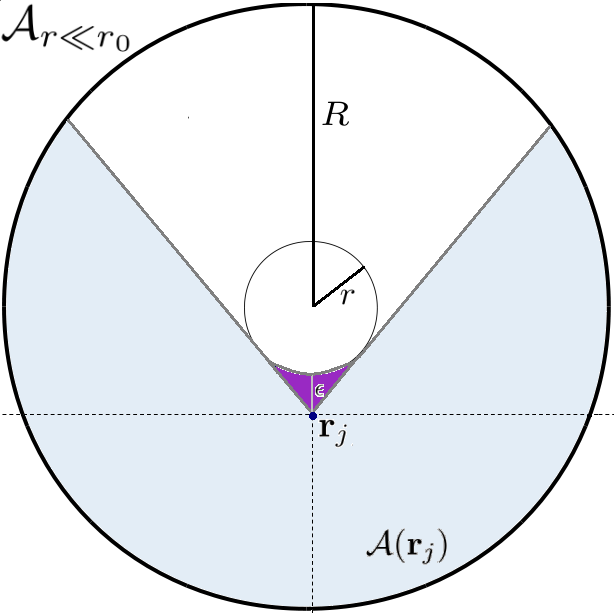}\qquad\includegraphics[scale=0.35]{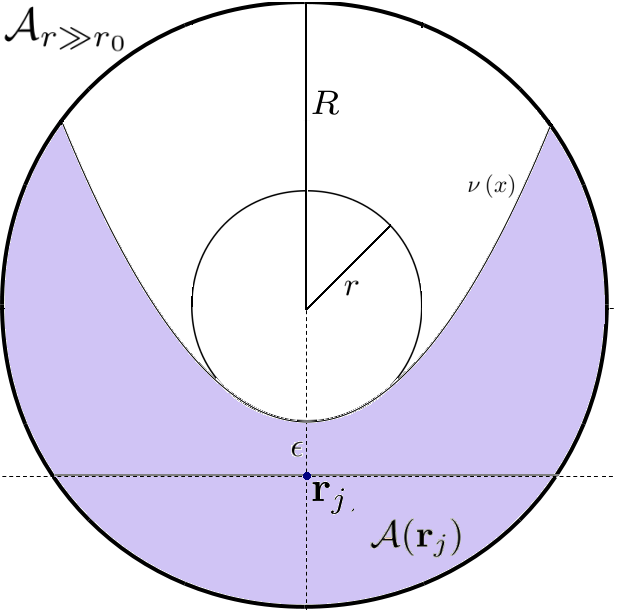}
\par\end{centering}
\caption{A depiction of the integration regions used for the annulus domain
$\mathcal{{A}}$ with small obstruction (middle panel) and large obstruction
(right panel), with the integration regions highlighted. The small, cone-like
region in the middle domain
$\mathcal{A}$ is highlighted in purple.}\label{fig:trip}
\end{figure}
We first take the case where $r=0$ depicted in Fig \ref{fig:trip}. This is
the disk $\mathcal{D}$. We first derive an approximation to $P_{fc}$ in this
limiting domain.

\begin{figure}
\noindent \begin{centering}
\includegraphics[scale=0.3]{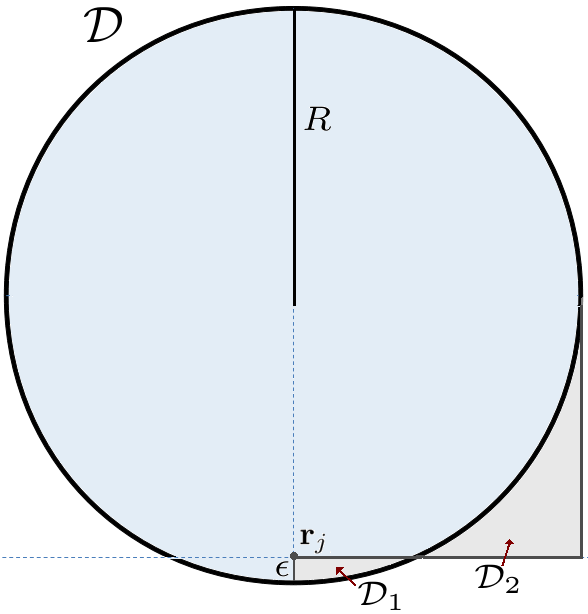}
\par\end{centering}
\caption{A depiction of the integration regions used for the disk domain
$\mathcal{{D}}$. The small, cone-like region in the middle domain
$\mathcal{A}$ is highlighted in purple.}\label{fig:diskregions}
\end{figure}

Firstly, the connectivity mass a distance $\epsilon$ from the disk's centre
is
\begin{eqnarray} \label{eq:neweq16}
\mathcal{{M}}\left(\epsilon\right) & = &
\frac{\pi}{2\beta}+2\left(\int_{\mathcal{{D}}_{1}}
e^{-\beta\left(x^{2}+y^{2}\right)}\textrm{d}y\textrm{d}x-\int_{\mathcal{{D}}_{2}}
e^{-\beta\left(x^{2}+y^{2}\right)}\textrm{d}y\textrm{d}x\right) \nonumber \\
& = & \frac{\pi}{2\beta}-2\int_{0}^{R}\int_{0}^{\epsilon-\sqrt{R^{2}-x^{2}}}	e^{-\beta\left(x^{2}+y^{2}\right)}\textrm{d}y\textrm{d}x
\end{eqnarray}
since the integral over $\mathcal{{D}}_{1}$ cancels. Now, consider two
regimes for the distance $\epsilon$:
in the first, where $\epsilon \approx R$ (close to the boundary), we can make
the approximation
$\exp\left({-\beta y^{2}}\right)\approx1$, since the distances $y$ from
the horizontal to the lower semi-circle in Fig. \ref{fig:diskregions}
will be small, so we can approximate the integral in Eq. \ref{eq:neweq16}
\begin{eqnarray}\label{eq:diskmass}
\int_{0}^{R}\int_{0}^{\epsilon-\sqrt{R^{2}-x^{2}}}e^{-\beta\left(x^{2}+y^{2}\right)}\textrm{d}y\textrm{d}x
& \approx & \int_{0}^{R}\int_{0}^{\epsilon-\sqrt{R^{2}-x^{2}}}e^{-\beta
x^{2}}\textrm{d}y\textrm{d}x \nonumber \\
& = & \frac{\sqrt{\pi}}{2\sqrt{\beta}}\epsilon-\int_{0}^{R}e^{-\beta
x^{2}}\sqrt{R^{2}-x^{2}}\textrm{d}x \nonumber
\end{eqnarray}
such that
\begin{equation}
\mathcal{{M}}(\epsilon \approx R) =
\frac{\pi}{2\beta}-\frac{1}{R\sqrt{\beta}}\left(\frac{\sqrt{\pi}}{4\beta}\right)+\left(R-\epsilon\right)\sqrt{\frac{\pi}{\beta}}+\mathcal{O}\left(\left(R-\epsilon\right)^2\right)
\end{equation}
after Taylor expanding Eq. \ref{eq:diskmass} for $\epsilon \approx R$, since
the mass is smallest on the boundary, dominating Eq. \ref{e:connectprob2}

For the other regime where $\epsilon \ll R$:
\begin{eqnarray}\label{eq:neweq18}
\mathcal{{M}}(\epsilon \ll R) &\approx& \int_{0}^{2\pi}\int_{0}^{\infty}r'
e^{-\beta r'^{2}} \textrm{d}r'\textrm{d}\theta \nonumber \\
	&=& \pi / \beta
\end{eqnarray}
due to the exponential decay of the connectivity function, and so we have the
probability of connection $P_{fc}$ in the disk domain:
\begin{eqnarray}\label{eq:4}
P_{fc} & \approx & 1 - \rho \int_{\mathcal{V}} e^{-\rho \int_{\mathcal{V}}
\chi \left(x,y\right) H\left(r_{xy}\right) \mathrm{d}y } \mathrm{d}x
\nonumber \\
& = & 1 - \rho
\int_{0}^{2\pi}\int_{0}^{L^{+}}\epsilon\exp\left(-\frac{\rho\pi}{\beta}\right)\textrm{d}\epsilon\textrm{d}\theta
\nonumber \\
& & -
\rho\int_{L^{+}}^{R}\exp\left(-\rho\left(\left(\frac{\pi}{2\beta}-\frac{1}{R\sqrt{\beta}}\left(\frac{\sqrt{\pi}}{4\beta}\right)\right)+\left(R-\epsilon\right)\sqrt{\frac{\pi}{\beta}}\right)\right)\epsilon\textrm{d}\epsilon\textrm{d}\theta
\nonumber \\
& \approx & 1-\pi R^{2}\rho e^{-\frac{\rho\pi}{\beta}}-2\pi
R\sqrt{\frac{\beta}{\pi}}e^{-\frac{\rho}{\beta}\left(\frac{\pi}{2}-\frac{1}{R\sqrt{\beta}}\left(\frac{\sqrt{\pi}}{4}\right)\right)}
\end{eqnarray}
where $L^+$ is the point where the two mass approximations equate. 
This approaches equation Eq. 38 of reference \cite{cef2012} as
$R\sqrt{\beta}\rightarrow\infty$. The correction
$\frac{\pi}{2}-\frac{1}{R\sqrt{\beta}}\left(\frac{\sqrt{\pi}}{4}\right)$ in
the exponential is a curvature correction to the previous result, since
before the boundary mass is expanded in a Taylor series truncated to its
first term, equivalent to ignoring the curvature of the boundary.

Monte-Carlo simulations, where graphs are drawn algorithmically and
enumerated if they connect, are presented in Fig. \ref{fig:montecarlo}
alongside our approximation in Eq. \ref{eq:4}, corroborating our
approximation. The simulations show an improvement on previous result in
\cite{cef2012}. The discrepancy at low density is expected since we only
consider the probability of a single isolated vertex. We also highlight the
interesting composition \ref{eq:4}. There is a bulk term
(whose coefficient is proportional to the area of $\mathcal{D}$) and a
boundary term (proportional to the circumference of $\mathcal{D}$).
This is discussed in greater detail in e.g. \cite{cef2012}, though we again
emphasise the dominance
of the boundary term as $\rho \to \infty$.
\subsection{Small obstacles}
Consider the case where $r \ll r_0$. We make the approximation that the small
cone-like domain $\mathcal{{A}}_{c}$ making up a portion
of the region visible to $x$, denoted $\mathcal{A}\left(x\right)$, and shown
in the middle panel of
Fig. \ref{fig:trip}, is only significantly
contributing to the connectivity mass at small displacements $\epsilon$ from
the obstacle, since at larger
displacements it thins and the
wedge-like region $\mathcal{A}\left(x\right) \setminus \mathcal{A}_{c}$ 
dominates.
Practically, it is $\mathcal{{A}}_{c}$ that presents the main integration
difficulties,
so we approximate $H\left(\norm{x-y}\right)$ over this region
where the radial coordinate $r'\ll1$, using $\exp{-\left(\beta r'^{2}\right)}
\approx 1$
\begin{eqnarray}\label{eq:neweq17}
\mathcal{M}\left( \epsilon \ll r_0 \right) & \approx &
\int_{-\pi+\arcsin\left(\frac{r}{r+\epsilon}\right)}^{\pi-\arcsin\left(\frac{r}{r+\epsilon}\right)}\int_{0}^{\infty}
e^{-\beta r'^{2}}r'\textrm{d}r'\textrm{d}\theta \nonumber \\
&+&2\int_{0}^{\arcsin\left(\frac{r}{r+\epsilon}\right)}\int_{0}^{\left(r+\epsilon\right)\cos(\theta)-\sqrt{r^{2}-\left(r+\epsilon\right)^{2}\sin^{2}(\theta)}}r'\textrm{d}r'\textrm{d}\theta
\nonumber \\
& = &
\frac{1}{\beta}\left(\pi-\arcsin\left(\frac{r}{r+\epsilon}\right)\right)
\nonumber \\
	&+&\int_{0}^{\arcsin\left(\frac{r}{r+\epsilon}\right)}
\left(\left(r+\epsilon\right)\cos\left(\theta\right)-\sqrt{r^{2}-\left(r+\epsilon\right)^{2}\sin^{2}\left(\theta\right)}\right)^2\textrm{d}\theta
\nonumber \\
& = &
\frac{\pi}{\beta}+\left(r^{2}-\frac{1}{\beta}\right)\arcsin\left(\frac{r}{r+\epsilon}\right)+r\sqrt{2r\epsilon+\epsilon^{2}}-\frac{\pi}{2}r^{2}
\nonumber
\end{eqnarray}
and then expand, giving
\begin{equation}
\mathcal{M}\left(\epsilon\ll r_0 \right) =
\frac{\pi}{2\beta}+\frac{\sqrt{2}}{\beta\sqrt{r}}\epsilon^{1/2}+\frac{8\beta
r^{2}-5}{6\beta\sqrt{2}r^{3/2}}\epsilon^{3/2}+\mathcal{O}\left(\epsilon^{2}\right)
\end{equation}
leaving us to integrate over the annulus
\begin{eqnarray} \label{eq:7}
P_{fc} & \approx & 1 - \rho \int_{\mathcal{{A}}}
e^{-\rho\mathcal{{M}}\left(x\right)}\textrm{d}x\nonumber \\
& \approx & 1 \nonumber \\ &-&
\rho\int_{0}^{2\pi}\int_{0}^{L^{-}}\left(r+\epsilon\right)\exp\left(-\rho\left(\frac{\pi}{2\beta}+\frac{\sqrt{2}}{\beta\sqrt{r}}\epsilon^{1/2}+\frac{8\beta
r^{2}-5}{6\beta\sqrt{2}r^{3/2}}\epsilon^{3/2}\right)\right)\textrm{d}\epsilon
\textrm{d}\theta \nonumber \\
& & - \pi R^{2}\rho e^{-\frac{\rho\pi}{\beta}}-2\pi
R\sqrt{\frac{\beta}{\pi}}e^{-\frac{\rho}{\beta}\left(\frac{\pi}{2}-\frac{1}{R\sqrt{\beta}}\left(\frac{\sqrt{\pi}}{4}\right)\right)}
\nonumber \\
& \approx & 1 \nonumber \\ &-&
2\pi\rho\int_{0}^{L^{-}}\left(r+\epsilon\right)
e^{-\frac{\rho\pi}{2\beta}}e^{-\rho\frac{\sqrt{2}}{\beta\sqrt{r}}\epsilon^{1/2}}\left(1-\rho\frac{8\beta
r^{2}-5}{6\beta\sqrt{2}r^{3/2}}\epsilon^{3/2}\right)\textrm{d}\epsilon - \pi
R^{2}\rho e^{-\frac{\rho\pi}{\beta}} \nonumber \\
&-&2\pi
R\sqrt{\frac{\beta}{\pi}}e^{-\frac{\rho}{\beta}\left(\frac{\pi}{2}-\frac{1}{R\sqrt{\beta}}\left(\frac{\sqrt{\pi}}{4}\right)\right)}
\nonumber \\
& \approx & 1 - \pi
r^{2}\frac{2\beta^{2}}{\rho}e^{-\frac{\rho\pi}{2\beta}}-\pi R^{2}\rho
e^{-\frac{\rho\pi}{\beta}}-2\pi
R\sqrt{\frac{\beta}{\pi}}e^{-\frac{\rho}{\beta}\left(\frac{\pi}{2}-\frac{1}{R\sqrt{\beta}}\left(\frac{\sqrt{\pi}}{4}\right)\right)}
 \end{eqnarray}
where $L^{-}$ is the point where the connectivity mass in the bulk meets our
approximation $\mathcal{M}\left(\epsilon\ll r_0 \right)$ near the obstacle.
We numerically corroborate Eq. \ref{eq:7} in Fig. \ref{fig:montecarlo} using
Monte Carlo simulations. 

Note that this obstacle term is extremely small
compared to the other contributions in Eq. \ref{eq:7}, given its
coefficient decays linearly with $\rho$ and the factor of
$\left(r\sqrt{\beta}\right)^2 \ll 1$. We conclude that
a small internal perimeter of radius $r$ in any convex domain $\mathcal{V}$
results in
a negligible effect on connectivity. 
\subsection{Large obstacles}
For the large obstacle case $r \gg r_0$, the relevant mass is
\begin{eqnarray}\label{eq:8}
\mathcal{M}\left(\epsilon \ll r_0 \right) & \approx & 
2\int_{0}^{\infty}\int_{0}^{\infty}e^{-\beta(x^{2}+y^{2})}\textrm{d}x\textrm{d}y+\int_{-\infty}^{\infty}\int_{0}^{
\epsilon +
\frac{1}{2r}x^2}e^{-\beta(x^{2}+y^{2})}\textrm{d}y\textrm{d}x\nonumber \\
 & = & \frac{\pi}{2\beta}+\frac{\sqrt{\pi}}
{2\sqrt{\beta}}\int_{-\infty}^{\infty}e^{-\beta
x^{2}}\textrm{erf}\left[\left(\epsilon+\frac{1}{2r}x^{2}\right)\sqrt{\beta}\right]\textrm{d}x\nonumber
\\
\end{eqnarray}
yielding a power series in $\epsilon$
\begin{eqnarray}
\mathcal{M}\left(\epsilon \ll r_0 \right) & \approx &
\frac{\pi}{2\beta}+\frac{\pi}{2\beta}\textrm{erf}\left[\sqrt{\beta}\epsilon\right]+\frac{1}{r\sqrt{\beta}}\left(\frac{\sqrt{\pi}}{4\beta}e^{-\beta\epsilon^{2}}\right)\nonumber
\\
& = &
\frac{\pi}{2\beta}+\frac{1}{r\sqrt{\beta}}\left(\frac{\sqrt{\pi}}{4\beta}\right)+\frac{\sqrt{\pi}}{\sqrt{\beta}}\epsilon+\mathcal{O}\left(\epsilon^{3/2}\right)
\end{eqnarray}
This implies the connectivity mass is scaling in the same way as
for the outer boundary, but where the curvature correction is of opposite
sign. We therefore have
\begin{equation}\label{eq:11}
P_{fc} \approx 1-2\pi
r\sqrt{\frac{\beta}{\pi}}e^{-\frac{\rho}{\beta}\left(\frac{\pi}{2}+\frac{1}{r\sqrt{\beta}}\left(\frac{\sqrt{\pi}}{4}\right)\right)}-
\pi R^{2}\rho e^{-\frac{\rho\pi}{\beta}}-2\pi
R\sqrt{\frac{\beta}{\pi}}e^{-\frac{\rho}{\beta}\left(\frac{\pi}{2}-\frac{1}{R\sqrt{\beta}}\left(\frac{\sqrt{\pi}}{4}\right)\right)}
\end{equation}
which is corroborated numerically in Fig. \ref{fig:montecarlo}.

This implies that large obstacles behave like separate, internal perimeters.
In the large-domain limit (where the node numbers go
to infinity and the connection range is tiny compared to the large domain
geometry), we can thus use
\begin{equation}\label{e:neweq15}
P_{fc} \approx 1-2\pi
\left(R+r\right)\sqrt{\frac{\beta}{\pi}}e^{-\frac{\rho\pi}{2\beta}}-\pi(R^{2}-r^{2})\rho
e^{-\frac{\rho\pi}{\beta}}
\end{equation}
\section{The spherical shell $\mathcal{{S}}$} \label{sec:four}
Consider now the spherical shell domain $\mathcal{{S}}$
of inner radius $r$ and outer radius $R$, which is the
three-dimensional analogue of the annulus.
\subsection{Small spherical obstacles}
The region visible to the node at $x$ is again decomposed into two parts,
the three-dimensional version of $\mathcal{{A}}_{c}$,
called $\mathcal{{S}}_{c}$, and the rest of the region visible to $x$,
denoted
$\mathcal{{S}}(x)\setminus\mathcal{{S}}_{c}$. As in the annulus with the
small obstacle, we
approximate $H\left(\norm{x-y}\right)$ over this region
where the radial coordinate $r'\ll1$, remembering the axis of integration is
centered at $x$:
\begin{eqnarray}
\mathcal{M}_{\mathcal{S}_{c}}\left(\epsilon\right) & = &
\int_{\mathcal{{S}}_{c}}r'^{2}e^{-\beta r'^{2}}\sin\theta
\textrm{d}r'\textrm{d}\theta \textrm{d}\varphi \nonumber \\
& \approx & \int_{\mathcal{{S}}_{c}}r'^{2}\sin\theta
\textrm{d}r'\textrm{d}\theta \textrm{d}\varphi \label{e:spheremass}
\end{eqnarray}
Eq. \ref{e:spheremass} is then evaluated by breaking up $\mathcal{{S}}_{c}$
into the area of a cone of
radius $\lambda$, height $h$ and apex angle $2\theta_{c}$
\begin{eqnarray}
\lambda&=&\frac{r}{r+\epsilon}\sqrt{2r\epsilon+\epsilon^{2}} \nonumber \\
h&=&\frac{2r\epsilon+\epsilon^{2}}{r+\epsilon} \nonumber  \\
\theta_{c}&=&\arcsin\left(\frac{r}{r+\epsilon}\right) \nonumber 
\end{eqnarray}
Note that the apex is at a distance $\epsilon$ from the obstacle. The volume
of the cone like region is given by
\begin{eqnarray}
\mathcal{{M}}_{\mathcal{{S_{C}}}}\left(\epsilon\right) & = & \frac{1}{3}\pi
\lambda^{2}h-\frac{1}{6}\pi\left(r+\epsilon-h\right)\left(3\lambda^{2}+\left(r+\epsilon-h\right)^{2}\right)\nonumber
\\
& = & \frac{\epsilon^{2}\pi
r^{2}\left(\epsilon+2r\right)^{2}}{3\left(\epsilon+r\right)^{3}}-\frac{\epsilon^{2}\pi
r^{3}\left(2\epsilon+3r\right)}{3\left(\epsilon+r\right)^{3}}\nonumber \\
 & = & \frac{\epsilon^{2}\pi r^{2}}{3\left(\epsilon+r\right)}
\end{eqnarray}
Adding the mass over $\mathcal{{S}}\left( x
\right)\setminus\mathcal{{S}}_{c}$, we use the fact that
the full solid angle available to a bulk node is $4\pi$, and that the angle
$\omega\leq\Omega$
available to the node at $x$ is
\begin{eqnarray}
\omega & = &
\frac{1}{4\pi}\int_{0}^{2\pi}\int_{0}^{\theta_{c}}\sin\left(\theta\right)\textrm{d}\theta
\textrm{d}\varphi \nonumber \\
& = &
\frac{1}{2}\left(1-\cos\left(\arcsin\left(\frac{r}{r+\epsilon}\right)\right)\right)\nonumber
\\
& = &
\frac{1}{2}\left(1-\sqrt{1-\left(\frac{r}{r+\epsilon}\right)^{2}}\right)
\end{eqnarray}
such that
\begin{equation}
\int_{\mathcal{{S}}\left( x \right)\setminus\mathcal{{S}}_{c}}r'^{2}e^{-\beta
r'^{2}}\sin\theta \textrm{d}r'\textrm{d}\theta \textrm{d}\varphi
=\frac{\pi\sqrt{\pi}}{\beta\sqrt{\beta}}\left(1-\frac{1-\sqrt{1-\left(\frac{r}{r+\epsilon}\right)^{2}}}{2}\right)
\nonumber
\end{equation}
We then have $\mathcal{{M}}\left( \epsilon \ll r_0 \right)$
\begin{eqnarray}
\mathcal{{M}}\left( \epsilon \ll r_0 \right)&\approx&\frac{\epsilon^{2}\pi
r^{2}}{3\left(\epsilon+r\right)}+\frac{\pi\sqrt{\pi}}{\beta\sqrt{\beta}}\left(1-\frac{1-\sqrt{1-\left(\frac{r}{r+\epsilon}\right)^{2}}}{2}\right)
\nonumber \\
&=&\frac{\pi\sqrt{\pi}}{2\beta\sqrt{\beta}}+\frac{\pi\sqrt{\pi}}{\beta\sqrt{\beta}}\frac{1}{\sqrt{2r}}\epsilon^{1/2}+\frac{3\pi^{3/2}}{4\sqrt{2}\left(r\beta\right)^{3/2}}\epsilon^{3/2}+\mathcal{O}\left(\epsilon^{2}\right)
\nonumber \end{eqnarray}
which implies that small spherical obstacles reduce the connection
probability
within the unobstructed sphere domain $\mathcal{S}_{r=0}$ to give a
connection probability of
\begin{eqnarray}
P_{fc}^{\mathcal{S}_{r \ll r_0}} & \approx & P_{fc}^{\mathcal{S}_{r=0}}
\nonumber \\
&-& \rho
e^{-\rho\left(\frac{\pi\sqrt{\pi}}{2\beta\sqrt{\beta}}\right)}\int_{0}^{2\pi}\int_{0}^{\pi}\int_{0}^{L_{\mathcal{{S}}}^{-}}\left(r+\epsilon\right)^{2}\sin\left(\theta\right)
\nonumber \\ &&
e^{-\rho\left(\frac{\pi\sqrt{\pi}}{2\beta\sqrt{\beta}}+\frac{\pi\sqrt{\pi}}{\beta\sqrt{\beta}}\frac{1}{\sqrt{2r}}\epsilon^{1/2}+\frac{3\pi^{3/2}}{4\sqrt{2}\left(r\beta\right)^{3/2}}\epsilon^{3/2}\right)}\textrm{d}\epsilon
\textrm{d}\theta \textrm{d}\varphi \nonumber \\
& \approx & P_{fc}^{\mathcal{S}_{r=0}} - \frac{4}{3}\pi
r^{3}\left(\frac{12\beta^{3}}{\rho\pi^{3}}\right)e^{-\rho\left(\frac{\pi\sqrt{\pi}}{2\beta\sqrt{\beta}}\right)}
\end{eqnarray}\newline

\begin{figure}
\begin{centering}\hspace{15mm}\includegraphics[scale=0.21]{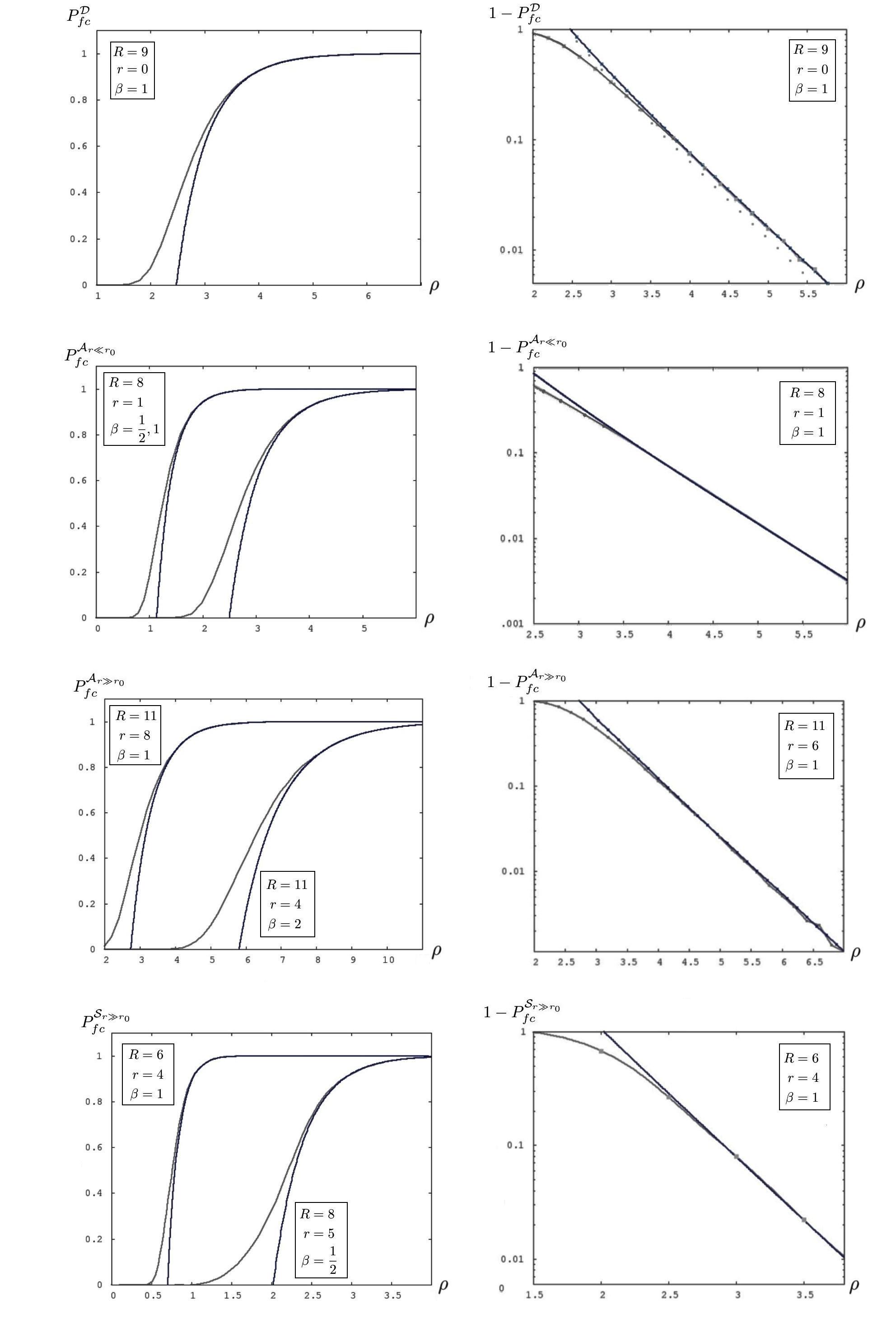}\end{centering}
\caption{We numerically estimate the connection probability of soft random
geometric graphs drawn inside the annuli $\mathcal{A}$ and spherical shell
$\mathcal{S}$. Every curve is
compared with our analytic predictions given by Eqs. \ref{eq:4}, \ref{eq:7},
\ref{eq:11} and \ref{eq:10}, where indicated.
The discrepancy at low density is expected due to the fact we calculate only
the probability of a single isolated vertex, and use its complement as an
approximation to the connection probability in a typical dense network
scenario.}\label{fig:montecarlo}
\end{figure}

\subsection{Large spherical obstacles}
For large obstacles ($r \gg r_0$), we extend Eq. \ref{eq:8} into the third
dimension.
$\mathcal{M}\left( \epsilon \ll r_0 \right)$ thus becomes
\begin{eqnarray}
\mathcal{M}\left( \epsilon \ll r_0 \right) \approx 
4\int_{0}^{\infty}\int_{0}^{\infty}\int_{0}^{\infty}\textrm{d}x\textrm{d}y\textrm{d}ze^{-\beta(x^{2}+y^{2}+z^{2})}
\nonumber \\
+\int_{-\infty}^{\infty}\int_{-\infty}^{\infty}\int_{0}^{\nu\left(x,z\right)}\textrm{d}y\textrm{d}x\textrm{d}ze^{-\beta(x^{2}+y^{2}+z^{2})}
\end{eqnarray}
where $\nu\left(x,z\right)=\epsilon+\frac{1}{2r}\left(x^{2}+z^{2}\right)$,
yielding
\begin{eqnarray}\label{eq:spheremass}
\mathcal{M}\left( \epsilon \ll r_0 \right) & \approx &
\frac{\pi\sqrt{\pi}}{2\beta\sqrt{\beta}}+\frac{\pi\left(r\sqrt{\beta\pi}\textrm{erf}\left[\epsilon\sqrt{\beta}\right]+e^{-\beta\epsilon^{2}}\right)}{2r\beta^{2}}\nonumber
\\
& = &
\frac{\pi\sqrt{\pi}}{2\beta\sqrt{\beta}}+\frac{\pi}{2\beta^{2}r}+\frac{\pi}{\beta}\epsilon+\mathcal{O}\left(\epsilon^{2}\right)
\end{eqnarray}
implying the connection probability is
\begin{eqnarray}
P_{fc}^{\mathcal{S}_{r \gg r_0}} & \approx & P_{fc}^{\mathcal{S}_{r=0}} -
\rho
e^{-\rho\left(\frac{\pi\sqrt{\pi}}{2\beta\sqrt{\beta}}+\frac{\pi}{2\beta^{2}r}\right)}\int_{0}^{2\pi}\int_{0}^{\pi}\int_{0}^{L_{\mathcal{{S}}}^{-}}\left(r+\epsilon\right)^{2}\sin\left(\theta\right)\textrm{d}\epsilon
\textrm{d}\theta \textrm{d}\varphi
e^{-\rho\left(\frac{\pi}{\beta}\epsilon\right)}\nonumber \\
&\approx& P_{fc}^{\mathcal{S}_{r=0}} - 4\pi
r^{2}\left(\frac{\beta}{\pi}\right)e^{-\rho\left(\frac{\pi\sqrt{\pi}}{2\beta\sqrt{\beta}}+\frac{1}{R\sqrt{\beta}}\left(\frac{\pi}{2\beta\sqrt{\beta}}\right)\right)}
\end{eqnarray}\newline
where $L^-$ is the point where our mass approximation in Eq.
\ref{eq:spheremass} is equal to the mass
in the bulk of the sphere $\left(\pi/\beta \right)^{3/2}$.

We now have the connection probability in the spherical shell $\mathcal{S}$
\begin{multline}\label{eq:10}
	P_{fc}^{\mathcal{S}} \approx 
1-\frac{4\pi}{3}\left(R^{3}-r^{3}\right)\rho
e^{-\rho\left(\frac{\pi\sqrt{\pi}}{\beta\sqrt{\beta}}\right)}
-4\pi
R^{2}\left(\frac{\beta}{\pi}\right)e^{-\rho\left(\frac{\pi\sqrt{\pi}}{2\beta\sqrt{\beta}}-\frac{1}{R\sqrt{\beta}}
\left(\frac{\pi}{2\beta\sqrt{\beta}}\right)\right)}\\
-	\begin{cases}
\frac{4}{3}\pi
r^{3}\left(\frac{12\beta^{3}}{\rho\pi^{3}}\right)e^{-\rho\left(\frac{\pi\sqrt{\pi}}{2\beta\sqrt{\beta}}\right)}
& \quad \text{if } r \ll r_0 \\
4\pi
r^{2}\left(\frac{\beta}{\pi}\right)e^{-\rho\left(\frac{\pi\sqrt{\pi}}{2\beta\sqrt{\beta}}
-\frac{1}{R\sqrt{\beta}}\left(\frac{\pi}{2\beta\sqrt{\beta}}\right)\right)} &
\quad \text{if } r \gg r_0
	\end{cases}
\end{multline}
which is corroborated in Fig. \ref{fig:montecarlo} for the large obstacle
case. Just as with the annulus,
small spherical obstacles have little impact on connectivity, and large
spherical obstacles behave like separate perimeters. This
behaviour is likely the same for all dimensions $d>3$, where the geometry is
a hypersphere containing a convex
$d$-dimensional obstacle. 
\begin{figure}[!]
\noindent
\begin{centering}
\includegraphics[scale=0.17]{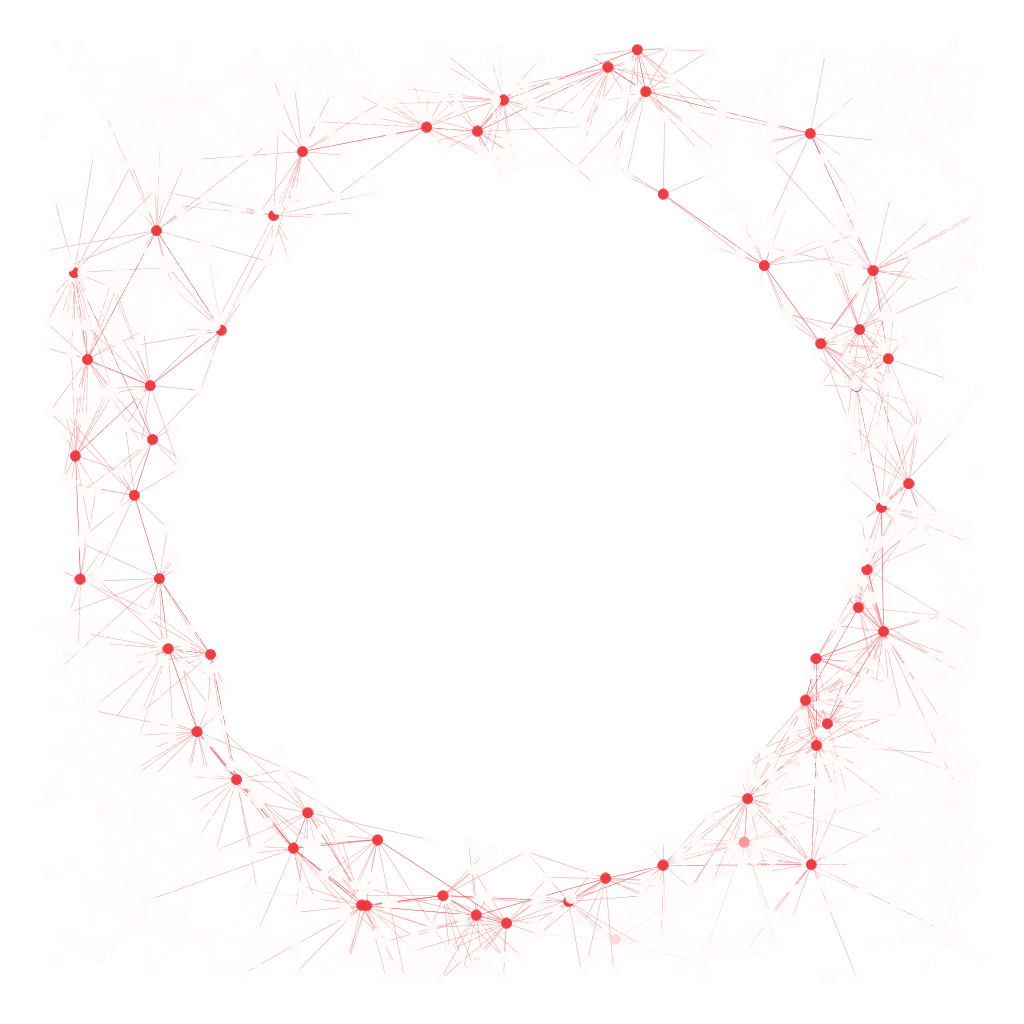}
\includegraphics[scale=0.23]{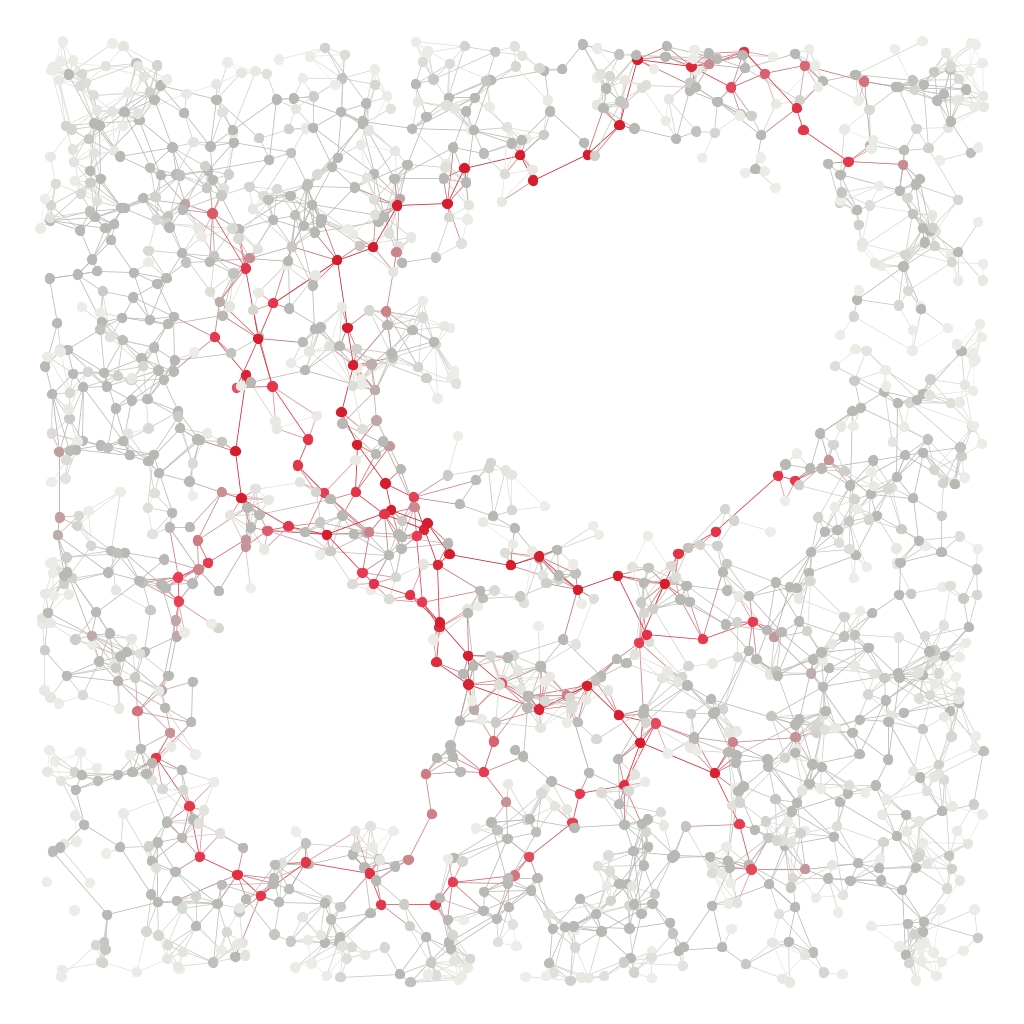}
\par\end{centering}
\caption{A soft random geometric graph inside a Sinai domain, and in a domain
with multiple obstructions. The betweenness centrality is plotted in light
tones (low) to darker tones eventually becoming red (high), showing the skeleton form around the
obstacles.}
\label{fig:skeleton}
\end{figure}

\section{Discussion} \label{sec:five}
We have derived analytic formulas for the connection
probability of soft random geometric graphs drawn inside various annuli and
shells given the link formation probability between two vertices is an
exponentially decaying function of their Euclidean separation. This models
the Rayleigh fading of radio signal propagation within a wireless \textit{ad
hoc} network.
\subsection{Numerical simulation}
Monte Carlo corroboration of Eqs. \ref{eq:4}, \ref{eq:7},
\ref{eq:11} and \ref{eq:10} are presented in Fig. \ref{fig:montecarlo}. We
numerically construct random graphs and count how many connect. This provides
a numerical estimate of the actual connection probability.
\subsection{Multiple obstacles}
We have thus extended the soft connection model into simple non-convex spaces
based on circular or spherical obstacles (rather than fractal boundaries
\cite{dettmann2014}, internal walls \cite{orestis2013} or fixed obstacles on
a grid \cite{almiron2013}). We highlight situations where obstacles are (and
are not) important influences on connectivity:
\begin{enumerate}
\item Small obstacles have little impact on connectivity.
\item Large obstacles have a similar impact on connectivity as the enclosing
perimeter, but their effects are dominated by the boundary as $\rho \to
\infty$.
\end{enumerate}
One may therefore be think that obstacles have little impact on dense network
connectivity. This is not true when they occur in great
numbers: given that obstacles are not too close, their effects add up in a linear
fashion, potentially outweighing
the effect of the boundary. To demonstrate this, take the Sinai-like domain 
in the right hand panel of Fig. \ref{fig:annulus}. Without obstacles, we have\begin{eqnarray}\label{eq:manyobs1}
P_{fc} = 1 - L^{2}\rho e^{-\frac{\pi}{\beta} \rho } - 4
L\sqrt{\frac{\beta}{\pi}} e^{-\frac{\pi}{2\beta} \rho} -
\frac{16\beta}{\rho\pi} e^{-\frac{\pi}{4\beta} \rho}
\end{eqnarray}
taken from \cite{cef2012}. This is composed of a bulk term, a boundary term
and a corner term. As we have seen,
introducing $n$ circular obstacles of various radii $r_i$ will reduce this
connection probability to:
\begin{eqnarray}\label{eq:manyobs2}
1 - \sum_{i=1}^{n}\pi r_{i}^{2}
\left(\frac{2\beta^{2}}{\rho}\right)e^{-\frac{\rho\pi}{2\beta}}
- \left(L^{2}-\sum_{i=1}^{n}\pi r_{i}^{2}\right)\rho e^{-\frac{\pi}{\beta}
\rho } - 4L \sqrt{\frac{\beta}{\pi}} e^{-\frac{\pi}{2\beta} \rho} -
\frac{16\beta}{\rho\pi} e^{-\frac{\pi}{4\beta} \rho} \nonumber
\end{eqnarray}
which holds whenever the obstacles are separated from each other and the
boundary
by at least $2r_0$.
\subsection{Surfaces without boundary}
Boundary effects can be removed by working on surfaces without an enclosing
perimeter. Examples include the flat torus, popular in rigorous studies but
difficult to realise in wireless networks, and the sphere. Thus as $\rho \to
\infty$ the obstacle effects are the dominant contribution to $P_{fc}$. This may 
be of interest to pure mathematicians studying random graphs for purposes
outside communication theory \cite{penrosebook}. Fractal obstacles may be of
particular interest \cite{dettmann2014}.
\subsection{Quasi-one-dimensional regime}\label{sec:1d}
Note that as the width of the annulus goes to zero, the approximation that
connectivity is the same as no isolated vertices breaks down. The graph now
disconnects by forming two clusters separated from each other by two
unpopulated strips of width usually greater than $ r_0$. Studying the
asymptotic connectivity of these quasi-1D random geometric graphs will be
topic of further study.
\subsection{Betweenness centrality near obstacles}
Vertex isolation near obstacles will have a significant effect on network
functionality. This is because vertices near obsatcles have exceptionally
high betweenness centrality, depicted in Fig. \ref{fig:skeleton}. Routing in
obstructed domains must take this effect into account, or vertices will
become overloaded near obstacles as they take on an excessive number of
routing tasks. Studying the connectivity properties of this ring of vertices
meandering around obstacles is another topic of further study.
\chapter{Betweenness centrality}\label{chapter:betweenness}

\section{Introduction}
Betweenness centrality is a graph theoretic measure of how often a vertex $z$
is on a shortest path of links between any other pair of vertices in a graph
\cite{freeman1977}. It is defined according to this sum:
\begin{eqnarray}\label{e:betweennessformula}
\gamma(z)=\sum_{i \neq j, j \neq k, k \neq i}\frac{\sigma_{r_{ij}}(z)}{\sigma_{r_{ij}}}
\end{eqnarray}
$\sigma_{r_{ij}}$ is the total number of shortest paths that join $i$ and
$j$, and $\sigma_{r_{ij}}(z)$ gives
the number of those geodesics that pass through $k$. Intuitively, nodes with
high betweenness can be thought of as
decisive for the functionality of decentralized communication networks, since
they typically route more data packets, based on the assumption that traffic
tries to follow only the shortest available multi-hop paths. Though the relation at this point between structural measures like betweenness, and related measures of actual network performance which follow from the underlying structure, is not clear, the existence of an important relationship is not in doubt. We seek to clarify this as part of our contribution.

This notion
of importance is in sharp contrast to methods which simply
enumerate node degrees, since a bridging
node which connects two large clusters is, for example, of crucial importance
to the whole network, even though
it may only have two neighbours. This sort of information
is brought out by betweenness centrality, but usually goes undetected. 

In router-based communication networks, the router itself has a normalised
betweenness of unity, since all nodes
connect to it directly, while all other nodes have a centrality of zero. In
ad hoc networks, and in sensor networks, betweenness
is distributed randomly at each vertex according to a distribution which
depends on the geometric location of vertices, implying a diverse betweenness
profile. Now, in wireless networks, this diversity can be harnessed in at
least three separate ways: in
2005 Gupta et al. \cite{gupta2005} used betweenness as a criteria for
electing \textit{cluster heads} which communicate to base stations on
behalf of all the cooperating machines. Later, in 2010, Ercsey-Ravasz et al.
\cite{ercsey2010} demonstrated how betweenness
can be used to delineate the network's \textit{skeleton} or
\textit{vulnerability backbone}, see also the more recent paper \cite{liu2015}, which is a percolating
cluster of the most structurally important vertices. Finally, in
2006, Wang et al. \cite{wang2006} researched the use of betweenness for
boundary detection, since at high vertex density the
betweenness of devices exhibits a bi-modal behaviour near the domain
boundary, and can therefore elucidate its location.

Given the insights from the last chapter concerning the structural importance
of vertices near non-convex features of domains such as connectivity
obstacles, in this chapter we develop an understanding of how the expected
betweenness of a vertex at
some domain location changes with the parameters of the model to which it
takes part, evaluating analytic formulas for betweenness
as a function of domain position.

\begin{figure}
\includegraphics[scale=0.22]{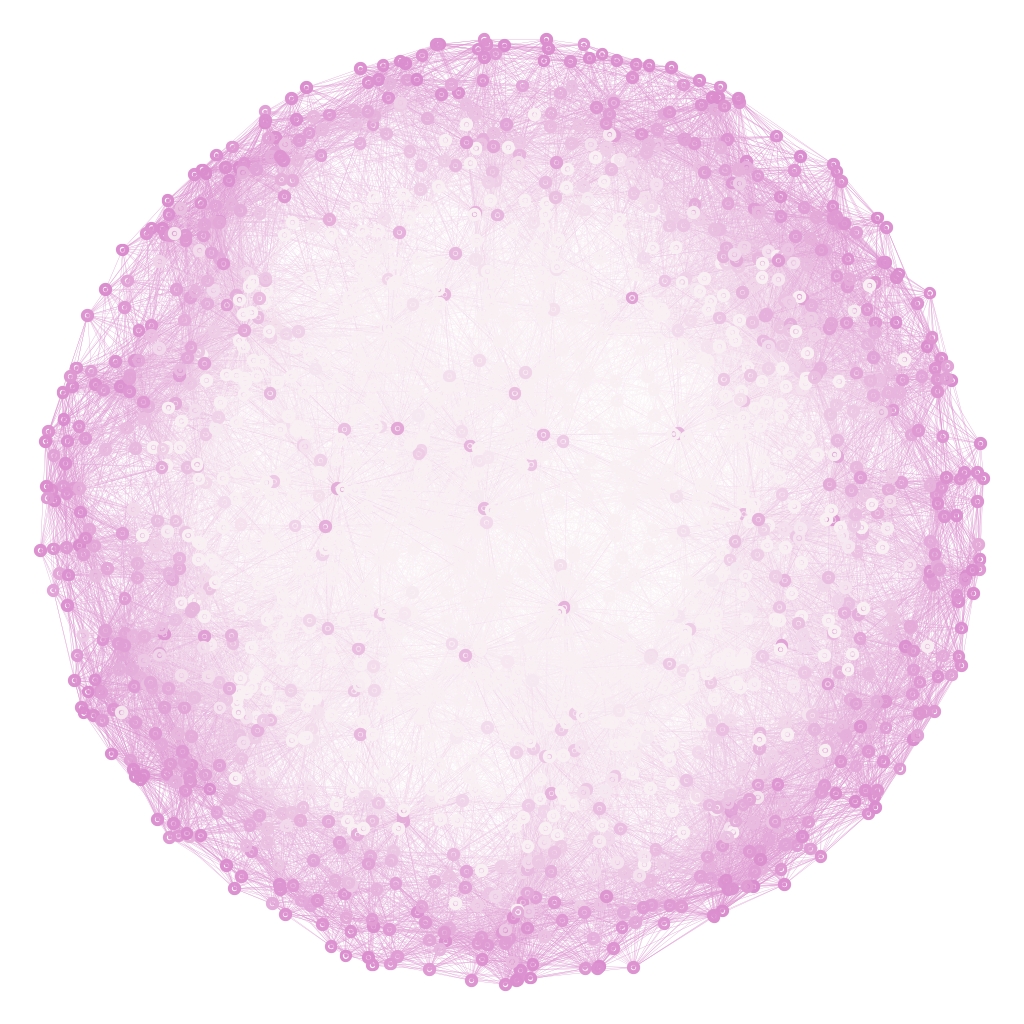}\qquad\qquad\quad
\includegraphics[scale=0.22]{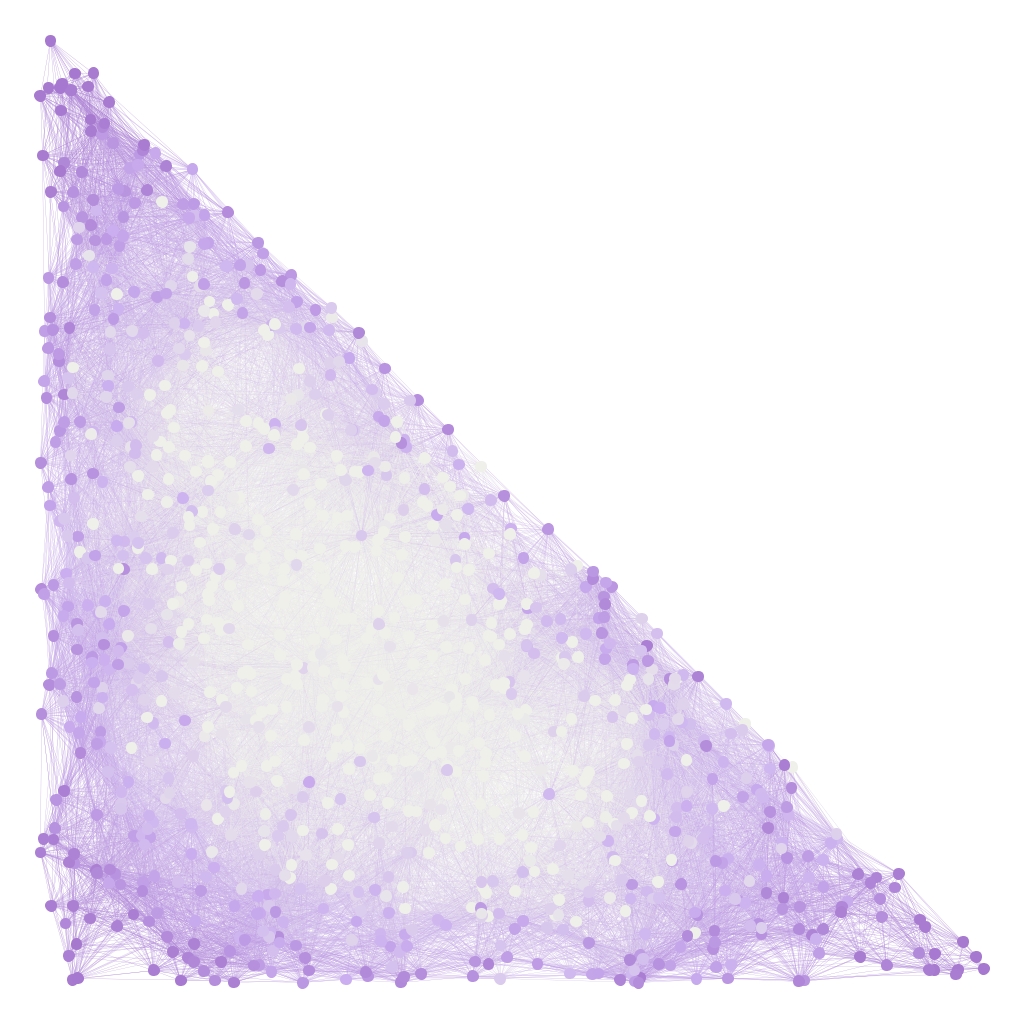}
\vspace{5mm}
\includegraphics[scale=0.175]{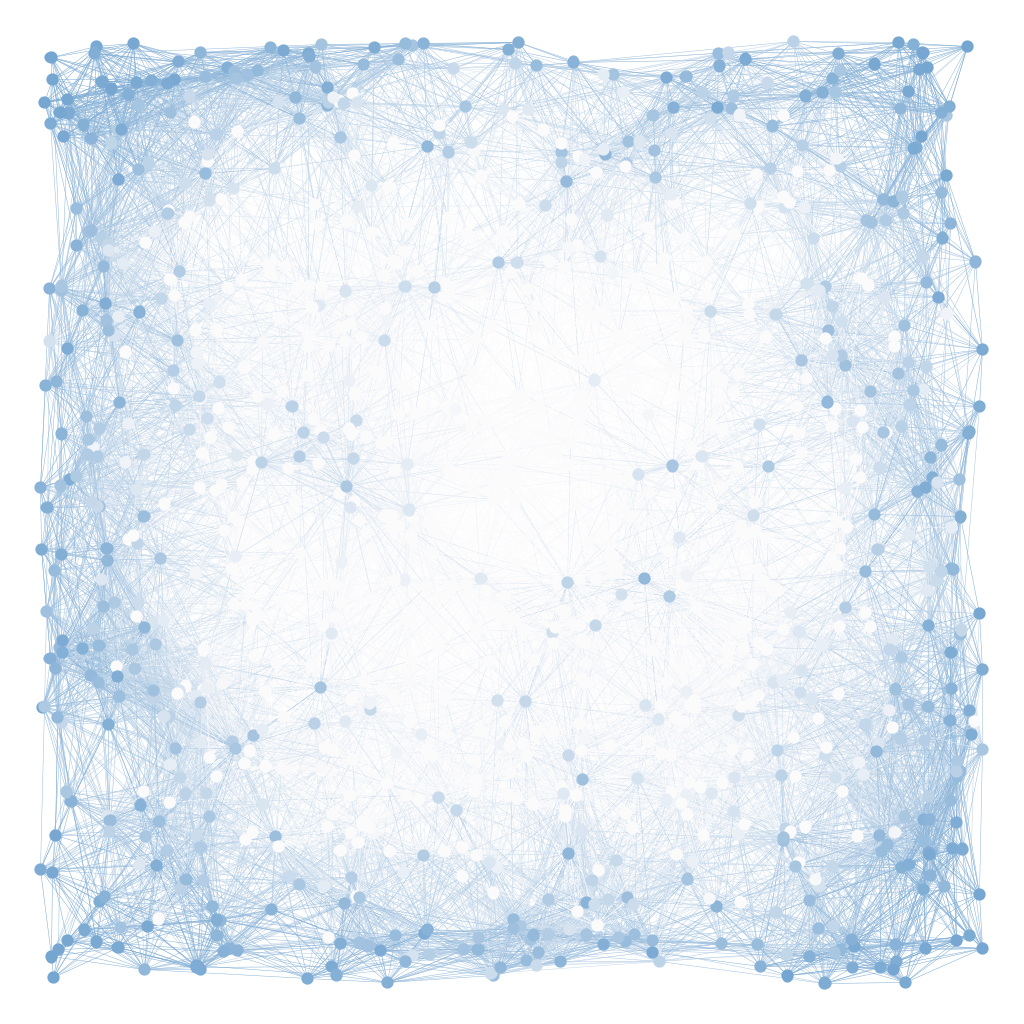}\qquad\qquad
\includegraphics[scale=0.22]{betweenness/corridor5}\newline
\caption{Four realisations of soft random geometric graphs and their
betweenness centrality bounded within various domains, including the disk
$\mathcal{{D}}$, square, right-angled triangle and square domain containing
two circular obstacles: in both the left and upper right figures the darker
colour represents low centrality, whereas the lighter colour high centrality,
whereas in the obstructed square domain (lower right) the least central nodes
are faded to grey and the most central are highlighted in red. Note that the
boundaries of the domains are locations where betweenness is at a minimum.
The link colours are based on the average betweenness of the two connected
nodes.}
\label{fig:one}
\end{figure}

We start our derivation with the disk domain $\mathcal{D}$ of
radius $R$ (Fig. \ref{fig:one}). We will consider a `continuum' density of
vertices, with vanishing connection range. This is for two reasons:
\begin{enumerate}
\item For the sake of tractability. 
\item To model a dense network.
\end{enumerate}

We then argue that betweenness, a computationally intensive operation with
possibly high communication overheads, can be well approximated by our
analytical closed form predictions, and can therefore prove useful in
practice.

This chapter is structured as follows: in Section \ref{sec:model} we present
our basic network model and state our main assumptions. In Section
\ref{sec:delta} we introduce an analytic formula for $\mathbb{E}(\gamma(z))$
in the continuum limit (where the node density $\rho\rightarrow\infty$),
which is our main result. In Section \ref{sec:montecarlo} we present Monte
Carlo simulations which corroborate our predictions, in Section
\ref{sec:discussion} we discuss amongst other issues the applicability of the derived betweenness
centrality formula within multi-hop wireless networks.
\section{The model}\label{sec:model}
Consider a soft random geometric graph formed by distributing vertices in a
bounded region $\mathcal{V} \subseteq \mathbb{R}^{d}$ according to a Poisson
point process $\mathcal{Y}$ of density $\rho$, and then adding an edge
between points $\{x,y\} \in \mathcal{Y}$ with probability $H(\norm{x-y})$,
where $H: \mathbb{R}^{2} \to \left[0,1\right]$ is the connection function
taken to be $\exp\left({-\beta \norm{x-y}^{2}}\right)$, as in the previous
chapter, and $\norm{x-y}$ is the Euclidean distance between vertices.
We consider only the \textit{continuum limit}, which is a sort of double
limit where $\rho\rightarrow\infty$ and the typical connection range goes to
zero \textit{in such a way that the graph remains connected}. In this limit
we make the assumption that all vertices on any straight line between any two
other vertices lie on the shortest path, measured in hops, that links those
two endpoint vertices. We seek the continuum analogue of Eq.
\ref{e:betweennessformula}.

Consider the domain $\mathcal{V}$ to have volume $V$ according to the
Lebesgue measure. The probability that some vertex is placed at position
$\mathbf{r_{\mathrm{i}}}$ in $\mathcal{V}$ is
$V^{-1}\mathrm{d}\mathbf{r_{\mathrm{i}}}$. Thus, the probability that any
vertex pair will simultaneously be placed at
$\{\mathbf{r}_{i},\mathbf{r}_{j}\}$ and construct between itself a shortest
path which passes through $z$ is
$V^{-2}\mathrm{d}\mathbf{r_{\mathrm{i}}}\mathrm{d}\mathbf{r_{\mathrm{j}}}
\chi_{ij}(z)$, where the characteristic function $\chi_{ij}(z)$ equates to
unity whenever $k$ lies on the path $i\rightarrow j$ given by the straight
line segment $\mathbf{r}_{ij}$ that joins $\mathbf{r_{\mathrm{i}}}$ and
$\mathbf{r_{\mathrm{j}}}$, and is zero otherwise. Summing this over all
possible $\{\mathbf{r}_{i},\mathbf{r}_{j}\}$ pair locations gives a continuum
approximation to the expected betweenness centrality of $z$:
\begin{equation}\label{e:2}
g(z)=\frac{1}{2V^{2}}\int_{\mathcal{{V}}}\mathrm{d}\mathbf{r_{\mathrm{i}}}\int_{\mathcal{{V}}}d\mathbf{r_{\mathrm{j}}}\,\chi_{ij}(z)
\end{equation}
We take $\mathcal{V}$ to be the disk doman $\mathcal{D}$
 of radius $R$, so $V=\pi R^{2}$. Note also
that due to the symmetry of the disk we can describe the position of the node
$k$ by its Euclidean distance $\epsilon$ from the disk's centre.
\begin{figure}
\begin{center} 
\includegraphics[scale=0.48]{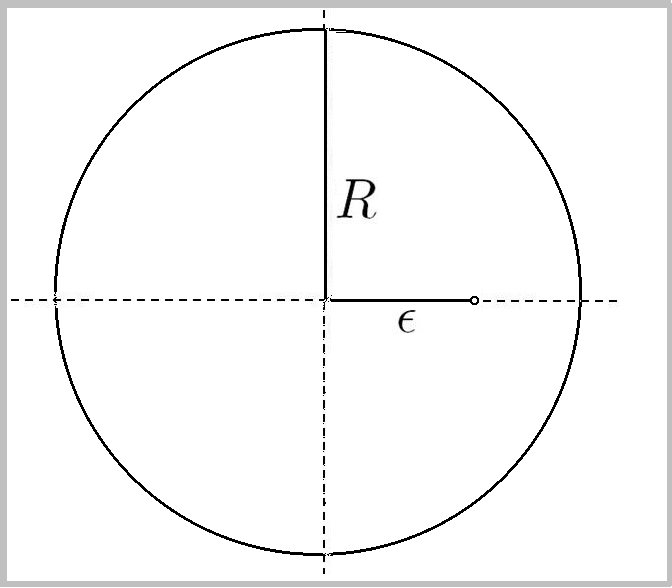} \\ \vspace{3mm}
\includegraphics[scale=0.5]{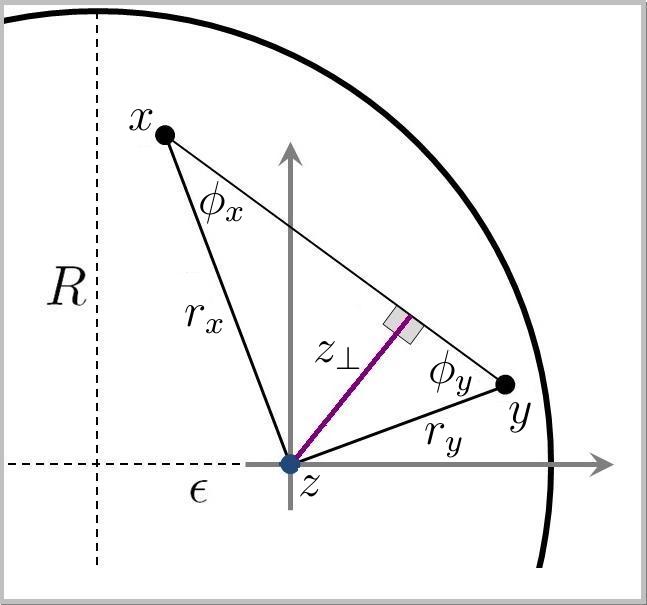}
\end{center}
\caption{The disk domain $\mathcal{D}$, and the three vertices $i$,$j$ and
$k$. We are trying to obtain a dense network approximation to the betweenness
centrality of a vertex placed at $k$. The scalar $k_{\bot}$ represents the
smallest Euclidean distance from $k$ to any point on the straight line
joining $i$ and $j$. The axis are centred on $k$, while the circle is centred
at $\left(-\epsilon,0\right)$. The angles $\phi_{i}$ and $\phi_{j}$ and the
distances $r_{i}$ and $r_{j}$ are also shown.}
\label{fig:kperp}
\end{figure}

Consider Fig. \ref{fig:kperp}. We define the scalar $z_{\bot}$ as the
distance of $z$ from the straight line $\mathbf{r}_{ij}$. Defining the delta
function
$\delta\left(z_{\perp}\left(\mathbf{r}_{\mathrm{i}},\mathbf{r}_{j}\right)\right)$,
we then write the following:
\begin{equation}\label{e:4}
\int_{\mathcal{{D}}}\mathrm{d}\mathbf{r_{\mathrm{i}}}\int_{\mathcal{{D}}}d\mathbf{r_{\mathrm{j}}}\,\chi_{ij}=\int_{\mathcal{{D}}}d\mathbf{r_{\mathrm{i}}}\int_{\mathcal{{D}}}d\mathbf{r_{\mathrm{j}}}\,\delta\left(z_{\perp}\right)
\end{equation}
The delta function will only contribute to the integral of Eq. \ref{e:4} when
its argument $k_{\perp}$ is a zero of $\delta\left(k_{\perp}\right)$. As
such, if we then describe $k_{\perp}$ such that it has a unique zero whenever
$k$ lies on the path $i\rightarrow j$, integrating
$\delta\left(k_{\bot}\right)$ over the
space of all node pairs $\{\mathbf{r}_{i},\mathbf{r}_{j}\}$ should
return $g(k)$ as required.

\section{The $\delta$ function}\label{sec:delta}
Fig. \ref{fig:kperp} shows $z$ located a distance $\epsilon$ from the centre
of $\mathcal{{D}}$, with the coordinate system centred on $z$ and orientated
such that the disk centre is at $(-\epsilon,0)$. Considering nodes $x$ and
$y$ at distances $r_{x}$ and $r_{y}$ from $x$ respectively, we have that the
internal angles $\phi_x$, $\phi_y$ and $(\theta_x - \theta_y)$ sum to $\pi$.
The perpendicular distance $z_\perp$ from $z$ to the line $\mathbf{r}_{xy}$
then satisfies both
\begin{equation}\label{e:101}
\frac{z_\perp}{r_x} = \sin(\phi_x) \nonumber
\end{equation}
and
\begin{equation}\label{e:104}
\frac{z_\perp}{r_y} = \sin(\phi_y) \nonumber
\end{equation}
Adding the above and taking small angle approximations (since we are
interested in the case where $ z_\perp\ll 1$) we have that
\begin{equation}\label{e:102}
\phi_x + \phi_y= \pi - \theta_y + \theta_x =z_\perp \left( \frac{1 }{r_x} +
\frac{1 }{r_y} \right)
\end{equation}
whenever $z_{\bot}\ll1$. This approximation presents a unique
zero of $z_{\bot}$ whenever $\theta_{x}-\theta_{y}+\pi=0$,
allowing\begin{eqnarray}\label{e:8}
\delta\left(z_{\bot}\right) & = &
\delta\left(\frac{\theta_{x}-\theta_{y}+\pi}{\frac{1}{r_{x}}+\frac{1}{r_{y}}}\right)\nonumber\\
& = &
\delta\left(\theta_{i}-\theta_{j}+\pi\right)\left(\frac{1}{r_{i}}+\frac{1}{r_{j}}\right)
\end{eqnarray}
due to the trivial scaling laws of the delta function. Eq. \ref{e:4}, a
double volume integral, becomes a quadruple integral
\begin{eqnarray}\label{e:9}
g(\epsilon) & = &
\frac{1}{2V^{2}}\int_{\mathcal{{D}}}\mathrm{d}\mathbf{r_{\mathrm{x}}}\int_{\mathcal{{D}}}\mathrm{d}\mathbf{r_{\mathrm{y}}}\,
\chi_{xy}\left( z \right) \nonumber\\
& = &
\frac{1}{2V^{2}}\int_{0}^{2\pi}\mathrm{d}\theta_{x}\int_{0}^{2\pi}\mathrm{d}\theta_{y}
\int_{0}^{r(\theta_{x})}r_{x}\mathrm{d}r_{x}\int_{0}^{r(\theta_{y})}r_{y}\mathrm{d}r_{y}\delta\left(z_{\bot}\right)
 \end{eqnarray}
Taking
$r(\theta)=\sqrt{R^{2}-\epsilon^{2}\sin^{2}(\theta)}-\epsilon\cos\left(\theta\right)$,
the polar equation of the circle bounding $\mathcal{{D}}$, Eq. \ref{e:9}
becomes
\begin{eqnarray}\label{e:901}
\frac{1}{2V^{2}}\int_{0}^{2\pi}\mathrm{d}\theta_{x}\int_{0}^{2\pi}\mathrm{d}\theta_{y}\delta\left(\theta_{x}-\theta_{y}+\pi\right)
\int_{0}^{r(\theta_{y})}r_{y}\mathrm{d}r_{y}\int_{0}^{r(\theta_{y})}\left(\frac{1}{r_{x}}+\frac{1}{r_{y}}\right)r_{x}\mathrm{d}r_{x}
\nonumber
\end{eqnarray}
which is
\begin{eqnarray}
\frac{1}{2V^{2}}\int_{0}^{2\pi}\mathrm{d}\theta_{x}\int_{0}^{2\pi}\mathrm{d}\theta_{y}\delta\left(\theta_{x}-\theta_{y}+\pi\right)
\left(r(\theta_{x})\frac{r^{2}(\theta_{y})}{2}+r(\theta_{y})\frac{r^{2}(\theta_{x})}{2}\right)
\nonumber
 \end{eqnarray}
then we integrate the delta function:
 \begin{eqnarray*}\label{e:902}
g(\epsilon) & = &
\frac{1}{4V^{2}}\int_{0}^{2\pi}\mathrm{d}\theta_{x}r(\theta_{x})r(\theta_{x}+\pi)\left(r(\theta_{x})+r(\theta_{x}+\pi)\right)\nonumber\\
& = &
\frac{1}{2V^{2}}\int_{0}^{2\pi}\mathrm{d}\theta_{x}\left(R^{2}-\epsilon^{2}\right)\sqrt{R^{2}-\epsilon^{2}\sin^{2}\left(\theta_{x}\right)}\end{eqnarray*}
leaving
\begin{equation}\label{e:unnormalised}
g(\epsilon)=\frac{2\left(R^{2}-\epsilon^{2}\right)}{\pi^{2}R^{3}}E\left(\frac{\epsilon}{R}\right)
\end{equation}
where
\begin{equation}\label{e:11}
E\left(k\right)=\int_{0}^{\pi/2}d\theta\sqrt{1-k^{2}\sin^{2}\left(\theta\right)}
\end{equation}
is the complete elliptic integral of the second kind (which is related to the
perimeter of an ellipse \cite{adlaj2012}). We normalise this to
$g^{\star}(\epsilon)$ by dividing Eq. \ref{e:unnormalised} by its maximum
value (which is at $\epsilon = 0$), to obtain our main result
\begin{eqnarray}
g^{\star}(\epsilon)
=\frac{2}{\pi}\left(1-\epsilon^{2}\right)E\left(\epsilon\right)
\label{e:12}
\end{eqnarray}
with $\epsilon$ in units of $R$ (and with betweenness now an element of the
unit interval).

Elliptic integrals cannot be easily visualised, so for clarification we can
expand Eq. \ref{e:12} near the origin (i.e. when $\epsilon \ll 1$) to obtain
\begin{equation}\label{e:quadratic}
g^\star (\epsilon\ll1) = 1- \frac{5 \epsilon^2}{R^2} + \frac{13
\epsilon^4}{64 R^4} + \mathcal{O}(\epsilon^6)
\end{equation}
while near the boundary (i.e. when $\epsilon \approx R$)
\begin{equation}
g^\star (\epsilon\approx R) = \frac{4(R-\epsilon)}{\pi R} +
\mathcal{O}((R-\epsilon)^2)
\end{equation}
which implies a quadratic scaling of betweenness near the centre, and a
linear scaling near the periphery.
\section{Discussion}\label{sec:discussion}
In this section we discuss numerical corroboration, and potential
applications of this centrality analysis.
\subsection{Numerical simulation}\label{sec:montecarlo}
\begin{figure}
\begin{center}
\includegraphics[scale=0.55]{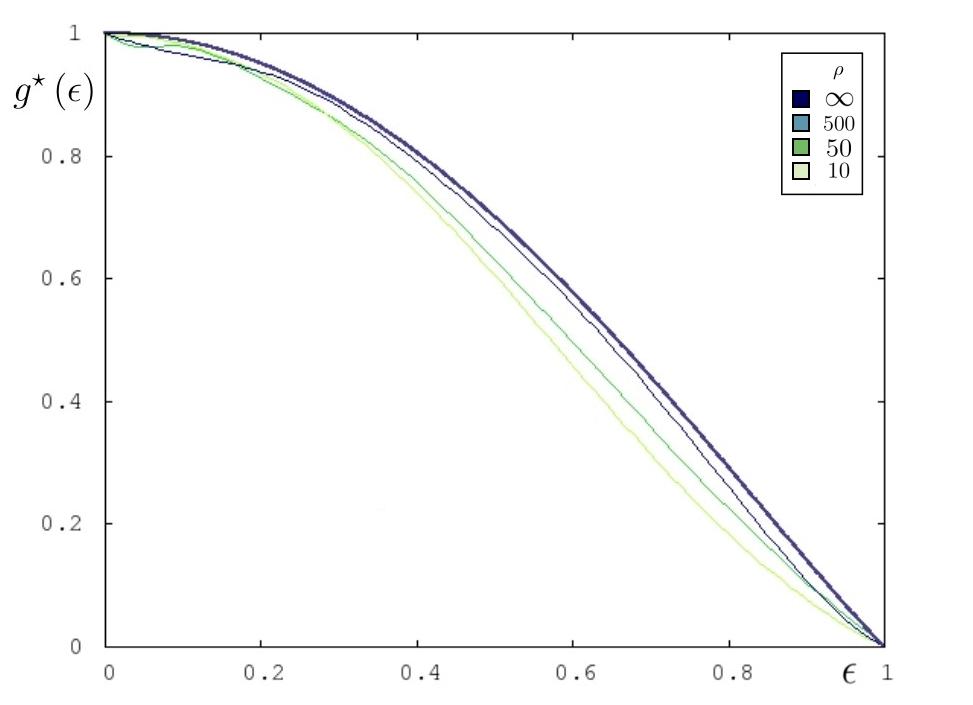}
\end{center}
\caption{Numerical evaluation of our continuum approximation Eq. \ref{e:12}.
Increasing the density of the point process, we ensure it contains a point
$z$ at distance $\epsilon$ from the centre of $\mathcal{D}$. A soft random
geometric graph with a Rayleigh fading connection function of vanishing range
$\beta \gg 1$, as in Chapter \ref{chapter:connectivity}, is then formed, and
the betweenness centrality of $z$ is calculated numerically. The thicker line
at the top is Eq. \ref{e:12}. The finite density simulations never converge
to this approximation, but the error is small.}
\label{fig:iccmain}
\end{figure}

Fig. \ref{fig:iccmain} shows that the betweenness $\gamma(\kappa)$ of nodes
situated in the bulk of $\mathcal{D}$ is typically high. Binning the
centrality in small increments of displacement from the domain centre and
averaging over many network realizations, we can plot this computationally acquired sample mean of
betweenness against location, and thus demonstrate how at finite densities
this betweenness approaches the continuum prediction of Eq. \ref{e:12}, which
we do in this figure. We take $\beta$ to be the largest value required for
full network connectivity, and increase $\rho$ from $10$ to $500$, each time
evaluating betweenness numerically using the algorithm defined in the
Mathematica 10 language.

The limit is never reached, only approached. At high density the discrepancy
is small. We propose this is simply due to the geodesic path not approaching
a straight line with density. There are also many geodesic paths in the
limit. So it appears this evident discrepancy will not go away by simply
increasing the density. Nevertheless, the betweenness does indeed appear to
be proportional to this mass of convex hulls which intersect $k$ and lie
within $\mathcal{D}$, and so this simple approximation may prove useful in
practice.
\subsection{Applications} \label{sec:applications}
We now discuss some possible applications of Eq. \ref{e:12}. Betweenness centrality has been used for electing cluster heads (CH)
\cite{gupta2005}). These are special vertices which take on the task of
routing collected data, often sensor data, to some distant point on behalf of
the network. In sensor networks, they are often called \textit{sinks}. A
number of routing protocols are usually implemented. For example, the basic
LEACH (Low Energy Adaptive Clustering Hierarchy \cite{heinzelman2000})
protocol uses a random selection of cluster heads at each round (i.e.
time-step), the vertices each taking turns in bearing the burden of cluster
head status, or, alternatively, EECS (Energy Efficient Clustering Scheme
\cite{liu2012}), which requires vertices to broadcast their remaining power
to first-degree neighbours, asking machines that find themselves with the
most battery power to then elect themselves to cluster head status.

However, in large networks using a vanishing transmitter range these
protocols don't work: far too many cluster heads get elected due to the huge
vertex numbers and the efficiency problem that this technique is trying to
mitigate re-arises. Potentially increasing transmitter range could resolve
the problem (since the usual techniques are based on one-hop access to the
head vertex), though this introduces interference problems, so one searches
for another solution.

Betweenness is a possible alternative election criteria, since it is
proportional to power consumption, and to routing load, unlike most other
centrality measures. This allows idle boundary vertices to act as cluster
heads whenever power minimisation is preferred, or busy domain-center
vertices whenever optimisation of vertex-to-vertex communication overheads is
tasked. Knowledge of betweenness as a function of position helps in the
selection of positions which, when occupied by vertices, results in CH
election.

Note also that, based on the intuition that central vertices are easier to
reach, communication-based resource consumption is minimised whenever
high-betweenness vertices are, in general, used as cluster heads.

Boundary detection is another application. This is an important field, with
various applications \cite{wang2006,dong2009,chen2012}. One potential use of
betweenness as a boundary detector is for mitigating the so called boundary
effect phenomenon \cite{cef2012}, where high-density network connectivity is
hampered through vertices becoming isolated near the domain boundaries. One
potential mitigation technique is to increase the device transmit power at
the domain boundary e.g. we can harness some spare power in the relatively
idle boundary vertices, increasing machine transmit power appropriately with
betweenness. This does not require the sharing of routing tables or other
connectivity information, since betweenness is directly proportional to the
devices current routing load. Finding the optimal function of the betweenness
(or perhaps other centrality measures) is beyond the scope of this thesis,
but we highlight that this is an interesting and important open problem.

Skeleton extraction is a third application of betweenness analysis in wireless networks. The skeleton \cite{liu2015} consists of the
most central vertices, given by a rather arbitrary percentage (the top $5
\%$, for example). Note the bottom right panel of Fig. \ref{fig:one}, where
betweenness is plotted over a square domain containing two circular obstacles
(which restrict line-of-sight connections between vertices, e.g.
\cite{orestis2013}). The skeleton \cite{liu2015} forms around the circular
obstacles. We discuss this further in Section \ref{sec:skeleton}.
\subsection{The advantages of betweenness}\label{sec:advantages}
The betweenness centrality of vertices stands alongside a variety of popular
centrality variants, many of which originated in sociology. See e.g. \cite{magaia2015} for a detailed discussion of those related to betweenness itself in the scope of delay tolerant communication, defined in section \ref{sec:stateoftheart}. We now discuss
the possible alternatives, though note that there are hundreds of existing variants that could be studied, e.g. analytically, with those briefly discussed in this section characteristic of the main themes. These themes include  considering edges or nodes, flows, or analysis of matrices such as the Laplacian see e.g. \cite{newman2003}.
\begin{enumerate}
\item \textbf{Edge betweenness centrality:} the number of shortest paths
which involve a particular \textit{edge} are considered \cite{wang2013}.
\item \textbf{Flow betweenness centrality:} a flow of material through the
network elucidates the importance of vertices \cite{freeman1991}.
\item \textbf{Node closeness centrality:} the sum of all geodesic distances
to the remaining $N-1$ vertices determines centrality, see e.g. \cite{wang2015}.
\item \textbf{Node eigenvector centrality:} a sophisticated extension of
degree centrality, a form of which is implemented in PageRank. Vertices are given centrality scores based on their proximity to other central vertices, which admits analysis using matrices \cite{pagerank}.
\item \textbf{Node degree centrality:} vertices are ranked simply by the
number of neighbours they have \cite{saxena2016}.
\end{enumerate}
Edge betweenness highlights communication channels which may become over
subscribed. This would be a good alternative investigation: it is just as
difficult to compute on a sparse graph \cite{wang2013}, and just as desirable in a closed,
analytic form, for example near an obstacle. It can then be used as a
complement to an algorithmic determination of edge betweenness centrality
(for various reasons related to the prediction of e.g. channel
characteristics). But as part of the current scope of centrality in wireless
networks, it is perhaps less relevant. It cannot be used for boundary
detection, for example, and it is sensitive to an exact specification of link
efficiency, which is a difficult to characterise element of wireless
performance.

Current flow betweenness is perhaps the better candidate for wireless
network-theoretic development than edge betweenness, see e.g. \cite{lulli2015}. It characterises the
intuitive reason vertices, or edges, are central: they can't help but be put
under pressure during the operation of a flow. It is, however, too
mathematically involved for an analysis of this sort. It may be more
appropriate to algorithmically determine, as which is discussed in the aforementioned reference.

Closeness is also available as a topic of study in communication networks \cite{santos2016}. This is
for the simple reason that short geodesic hop counts to all vertices
characterise multi-hop routing in a straightforward way. But it does not
capture the `bridging vertex' issues picked up by betweenness.

Degree centrality, and sophisticated extensions such as PageRank \cite{pagerank}, potentially suffer from
the same sort of tractability issues as current flow betweenness if we intended to go beyond calculating expectations, and instead look for a distribution of a vector of degrees, since spatial dependence between nearby vertex degrees is difficult to deal with. Though the
expected degree of a vertex is easier to calculate. The trade off
between the characterisation of flows which betweenness can provide, the
overall tractability of its first moment, and its
established use in ad hoc networks, provides the motivation for an analytic
study.
\subsection{Convergence}\label{sec:converge}
Despite the disparity between betweenness centrality and our continuum
approximation Eq. \ref{e:12}, Fig. \ref{fig:iccmain} shows a remarkably good
approximation, given we must take into account:
\begin{enumerate}
\item That there are many geodesic paths which join two vertices as $\rho \to
\infty$.
\item That the paths are not straight, but appear to form geometric functions (see e.g. the discussion of Schramm–Loewner Evolutions in \cite{grimmettbook}.
\end{enumerate}
Thus, the small discrepancy is a fair price to pay given the difficulty that
a more accurate, finite density approximation appears to present.

Also, we have pointed out the quadratic scaling of betweenness as vertices
moves away from the central region of the disk domain in Eq.
\ref{e:quadratic}. It is interesting to ask how sensitive this first order
term is to details of the boundary.
\chapter{Geodesic Paths}\label{chapter:geodesicpaths}

\section{Introduction}
Minimising the total number of sequential transmissions required for two
vertices to communicate is a common concern in multi-hop communication. This
is due to:
\begin{enumerate}
\item The distortion caused by excessive amplify-and-forward activity,	
\item The loss of SINR (signal-to-interference-plus noise ratio) caused by an overly intense spatio-temporal point
process of relay transmissions.
\end{enumerate}
Probabilistic modelling of the network’s route statistics can prove important
for the management of these issues, as well as surrounding concerns in multi-hop communication \cite{ta2007,mao2010,zhang2012}. This motivates a number of problems concerning the statistics of
paths, most prominently that of finding a function
$g : \mathbb{Z} \to \{\Omega,\mathcal{F},\mathbb{P}\}$ relating the minimum
number of hops between two vertices $x,y$ in a random geometric graph to a
probability distribution on the set of Euclidean distances which could
separate them in the metric space.
\begin{figure}[t]
\noindent \begin{centering}
\includegraphics[scale=0.57]{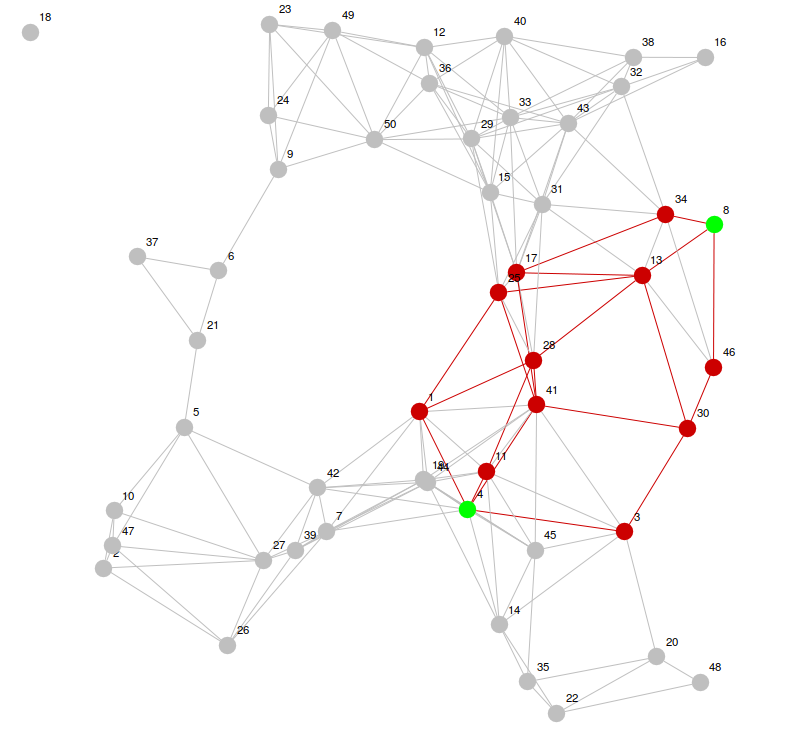}
\par\end{centering}
\caption{Vertices 4 and 8 are separated by Euclidean distance $2.4$ in a unit
disk graph, with the connection range $r_{0}=1$. They are joined by ten
geodesic paths, each of four hops. We evaluate the expectation of this
quantity in terms of the mutual Euclidean separation of vertices.}
\label{fig:paths}
\end{figure}
The origins of this problem go back to Chandler’s letter in the 1980’s
\cite{chandler1989}. As an example application, this
relationship can be used to provide location estimates to vertices possessing
hop-counts to ‘anchor’
devices distributed around the domain. This is a
particularly useful skill in wireless sensor networks where GPS is often
unavailable, see e.g. \cite{nguyen2015}.

The \textit{number} of geodesics is also potentially of interest. An example set of
geodesics is depicted in Fig. \ref{fig:paths}. This is the more general
statistic, since no $k$ hops paths suggest paths are at least $k+1$ hops. A
distribution on the number of $k$ hop paths can provide these probabilities.
Therefore, let $\mathcal{V} \subseteq \mathbb{R}^{d}$ be a bounded region of
volume $V$ associated with both the Lebesgue measure
$\mathrm{d}x$ and the Euclidean metric $r_{xy} = \norm{x-y}$ for any $x,y \in
\mathcal{V}$. Construct a \textit{unit disk graph} $\mathcal{G}\left(n,\pi
r_{0}^2\right)$ in $\mathcal{V}$ by deterministically linking pairs of a
Poisson point process $\mathcal{Y} \subseteq \mathcal{V}$ of density $\rho$
whenever $\norm{x-y} < r_{0}$.

Set $r_{0}=1$. Consider two nodes $x$ and $y$ in
$\mathcal{G}\left(n,\pi\right)$, and call these \textit{terminal} nodes. In
any graph realisation, these terminal nodes are either disconnected such that
no path exists between them, or they are linked by a variety of different
paths of various lengths. Consider now only the geodesic paths. Write
$\sigma_{r_{xy}}$ for the number of paths of this geodesic length joining $x$
and $y$, and separated by Euclidean distance $r_{xy}$. We name
$\sigma_{r_{xy}}$ the \textit{geodesic cardinality} of the distance $r_{xy}$,
since it gives the cardinality or `number of elements' $\sigma_{r_{xy}}$ in
the set of geodesic paths joining $x$ and $y$. Vertices lying mutually
displaced by $r_{xy}$ are expected to be linked by
$\mathbb{E}\left(\sigma_{r_{xy}}\right)$ geodesic paths. This is
\textit{expected geodesic cardinality}.
\begin{figure}[t]
\noindent \begin{centering}
\includegraphics[scale=0.5]{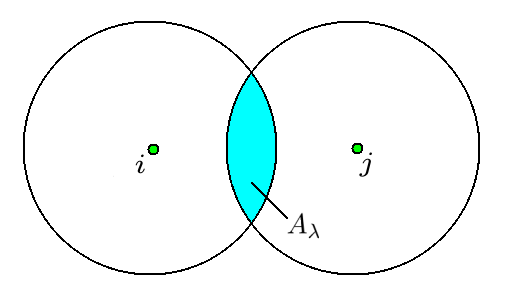}
\par\end{centering}
\caption{Points $x$ and $y$ separated by $\lambda \in [1,2)$ create an area
highlighted in blue, which is $A_{\lambda}$. Vertices falling within this
blue region lie on one of potentially many (what will necessarily be) shortest paths which run between these two vertices, each of two
hops.}
\label{fig:lens1}
\end{figure}

How long are these paths? They can be of potentially any number of hops
\textit{larger} than $\ceil{r_{xy}}$ (which is the smallest integer larger
than or equal to $r_{xy}$ known as its `ceiling'). To clarify this, we relate in Table
\ref{table:hops} a range of Euclidean distances and the minimum
number of hops required to join two vertices of that displacement $r_{xy}$.     
\begin{table}[!]
\begin{center}
	\begin{tabular}{| l | l | l | l |} \hline
    Displacement $r_{xy}$ & Minimum number of hops $i \to j$\\ \hline
    $r_{xy} \in [0,1)$ & $1$ hop \\ \hline
    $r_{xy} \in [1,2)$ & $2$ hops \\ \hline
	$r_{xy} \in [2,3)$ & $3$ hops \\ \hline
	$r_{xy} \in [3,4)$ & $4$ hops \\ \hline
	$r_{xy} \in [4,5)$ & $5$ hops \\ \hline
	\vdots & \vdots \\ \hline
    \end{tabular}
\end{center}
\caption{The minimum number of hops required for two
vertices to communicate in a unit disk graph, which is always the ceiling of the Euclidean distance separating the vertices.}
\label{table:hops}
\end{table}
We now look for an analytic expression $\sigma_{r_{xy}}$ in terms of
displacement $r_{xy}$, and density $\rho$. We only consider a density regime
such that all vertices display at least one path to every other vertex in the
graph which is of the \textit{shortest possible length} $k=\ceil{r_{xy}}$
hops. This ensures the work is relevant to the ultra-dense 5G small cell
scenario, and the associated communications problems, such as interference
minimisation, discussed in the introduction.
\section{A recursive formula for $\sigma_{r_{xy}}$}
We now work through the calculation for the expected number of
$\ceil{r_{xy}}$-hop paths in dimension $d=2$, which is asymptotic to the
number of geodesic paths as $\rho \to \infty$. Firstly,
\begin{equation}\label{e:2mk1}
	\mathbb{E}\left(\sigma_{r_{xy} \in \left[0,1\right)}\right)=1
\end{equation}
since nodes connect directly. For $r_{xy} \in \left[1,2\right]$, notice that
all nodes lying within the intersection of the two connection loci $x$ and
$y$, called `lenses', will necessarily lie on a geodesic two-hop path. This
is depicted in Fig \ref{fig:lens1}. Therefore,
\begin{eqnarray}\label{e:int1to2d2}
\mathbb{E}\left(\sigma_{r_{xy} \in \left[1,2\right)}\right)&=&\rho A_{r_{xy}}
\nonumber \\
                                          &=& \rho \left (2 \arccos \left(
\frac{\lambda}{2} \right) - \lambda \sqrt{ 1 - \frac{\lambda ^ 2}{4}}\right)
\end{eqnarray}
where $ A_{r_{xy}}x$ is the area of intersection of two unit circles at
center separation $\lambda$.
\begin{figure}
\noindent \begin{centering}
\includegraphics[scale=0.4]{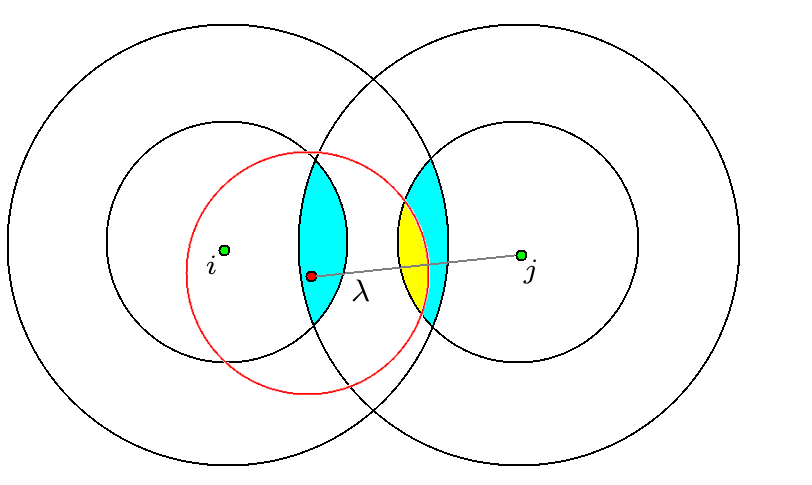}
\par\end{centering}
\caption{Some node (red) lying a distance $\lambda$ from $j$ (green) will
connect directly to all nodes within $r_0$ of its position (the red circle).
All of those vertices which simultaneously lie within a distance $r_0$ of
$j$, denoted by the yellow area, will form two hop paths, and thus each three
hop paths $i \leftrightarrow j$.}
\label{fig:lens2}
\end{figure}

For $r_{xy} \in \left[2,3\right)$, things are significantly more involved.
Fig. \ref{fig:lens2} shows the setup, which involves two lens, and we want
the expected number of bridging links between them. In order to evaluate this
expectation, construct a new area $A_{\lambda}$ in the right-hand lens.
Within this sub-region, vertices are Poissonly distributed with mean $\rho
A_{\lambda}$, and so:
\begin{eqnarray}
\mathbb{E}\left(\sigma_{r_{xy} \in
\left[2,3\right)}\right)&=&\int_{r_{xy}-1}^{2}\rho A_{\lambda} \rho
l_{\lambda} \mathrm{d}\lambda \label{e:int2to3d2}
\end{eqnarray}
where
\begin{equation}\label{e:arclength}
l_{\lambda}=2 \lambda \arccos\left(\frac{\lambda^2+r_{xy}^2-1}{2r_{xy}
\lambda}\right)
\end{equation}
is a contour in the first lens at a distance $\lambda$ from $j$, and
\begin{equation}\label{e:lensarea}
A_{\lambda}=2 \arccos \left( \frac{\lambda}{2} \right) - \lambda \sqrt{ 1 -
\frac{\lambda ^ 2}{4}}
\end{equation}
The expected number of three-hop paths is therefore given by an integral over
all positions in the left lens:
\begin{multline}
\mathbb{E}\left(\sigma_{r_{xy} \in \left[2,3\right)}\right) = \\
4\rho^2\int_{r_{xy}-1}^{2} \lambda
\arccos\left(\frac{\lambda^2+r_{xy}^2-1}{2r_{xy} \lambda}\right)\left(
\arccos \left( \frac{\lambda}{2} \right) - \frac{\lambda}{2} \sqrt{ 1 -
\frac{\lambda ^ 2}{4}}\right) \mathrm{d}\lambda
\end{multline}
Now, this integral has no closed form. Performing asymptotic analysis, Taylor
expand $A_{\lambda}$ at $\lambda=\ceil{r_{xy}}$ and $l_{\lambda}$ at
$\lambda=r_{xy}-1$. After dropping all but the first term, we should be able
to extract the leading order behaviour
\begin{multline}
\mathbb{E}\left(\sigma_{r_{xy} \in
\left[2,3\right]}\right)\approx4\rho^2\int_{r_{xy}-1}^{2}
\left(\sqrt{\frac{r_{xy}-1}{2r_{xy}}\left(\lambda+1-r_{xy}\right)}+\mathcal{O}\left(\left(\lambda+1-r_{xy}\right)^{3/2}\right)\right)
\\
\left(\frac{1}{3}\left(2-\lambda\right)^{3/2}+\mathcal{O}\left(\left(2-\lambda\right)^{5/2}\right)\right)\mathrm{d}\lambda
\nonumber
\end{multline}
which evaluates to
\begin{eqnarray}
	\rho^2 \frac{\pi}{3\sqrt{3}}\left(3-r_{xy}\right)^3 +
    \mathcal{O}\left(\left(3-r_{xy}\right)^4\right) \label{e:8}
\end{eqnarray}
after finally expanding about $r_{xy}=\ceil{r_{xy}}$ and taking the first
term. We can iterate this procedure to obtain an expression for
$\mathbb{E}\left(\sigma_{r_{xy} \in \left[3,4\right)}\right)$
\begin{eqnarray}\label{e:int3to4d2}
\mathbb{E}\left(\sigma_{r_{xy} \in
\left[3,4\right)}\right)&=&\rho^3\int_{r_{xy}-1}^{3}l_{\lambda_{1}}\mathrm{d}\lambda_{1}\int_{\lambda_{1}-1}^{2}
A_{\lambda_{2}} l_{\lambda_{2}} \mathrm{d}\lambda_{2} \nonumber
\end{eqnarray}
where 
\begin{eqnarray}
l_{\lambda_{1}}=2 \lambda_{1}
\arccos\left(\frac{\lambda_{1}^2+r_{xy}^2-1}{2r_{xy} \lambda_{1}}\right)
\label{e:10} \nonumber\\
l_{\lambda_{2}}=2 \lambda_{2}
\arccos\left(\frac{\lambda_{2}^2+\lambda_{1}^2-1}{2\lambda_{1}
\lambda_{2}}\right) \nonumber \label{e:11}
\end{eqnarray}
After expanding $l_{\lambda_{1}}$ about $r_{xy}-1$ and $l_{\lambda_{2}}$
about $=\lambda_{1}-1$ (as before), we have
\begin{eqnarray}
\mathbb{E}\left(\sigma_{r_{xy} \in \left[3,4\right)}\right)&=& \rho^3
    \frac{32\pi\sqrt{2}}{945}\left(4-r_{xy}\right)^{9/2} +
    \mathcal{O}\left(\left(4-r_{xy}\right)^{5}\right) 
\end{eqnarray}
after Taylor expanding Eq. \ref{e:int3to4d2} about $r_{xy}=\ceil{r_{xy}}$ and
taking the leading term. Clearly we can go on and produce a recursive formula
\begin{equation}\label{e:recursionrelation2d}
   \mathbb{E}\left(\sigma_{r_{xy}}\right) = 
\rho
\int_{r_{xy}-1}^{\floor{r_{xy}}}l_{\lambda}\mathbb{E}\left(\sigma_{\lambda}\right)\textrm{d}\lambda
\end{equation}
We have already performed a few of these integrals. Table \ref{table:2d}
lists them up to $r_{xy} \in [5,6)$.
\begin{table}
	\begin{center}
    \begin{tabular}{| l | l | l | l |}
    \hline
Displacement $r_{xy}$ & Expected geodesic cardinality $\sigma_{r_{xy}}$ &    Order of error \\ \hline
    $r_{xy} \in [0,1)$ & 1 & 
(no error) \\ \hline
    $r_{xy} \in [1,2)$ & $\frac{4\rho}{3}\left(2-r_{xy}\right)^{3/2}$ & 
    $\mathcal{O}\left(\left(2-r_{xy}\right)^{5/2}\right)$ \\ \hline
$r_{xy} \in [2,3]$ & $\frac{\pi\rho^2}{3\sqrt{3}}\left(3-r_{xy}\right)^{3}$    & $\mathcal{O}\left(\left(3-r_{xy}\right)^{8/2}\right)$ \\ \hline
    $r_{xy} \in [3,4)$ &
    $\frac{32\pi\rho^3\sqrt{2}}{945}\left(4-r_{xy}\right)^{9/2}$ 
    & $\mathcal{O}\left(\left(4-r_{xy}\right)^{11/2}\right)$ \\ \hline
$r_{xy} \in [4,5)$ &
$\frac{\pi^2\rho^4}{180\sqrt{5}}\left(5-r_{xy}\right)^{6}$
    & $\mathcal{O}\left(\left(5-r_{xy}\right)^{14/2}\right)$\\ \hline
$r_{xy} \in [5,6)$ &
$\frac{1,024\pi^2\rho^5}{2,027,025\sqrt{3}}\left(6-r_{xy}\right)^{15/2}$
    & $\mathcal{O}\left(\left(6-r_{xy}\right)^{17/2}\right)$ \\ \hline
    \end{tabular}
\caption{Solutions to the recursion relation
\ref{e:recursionrelation2d}.}\label{table:2d}
	\end{center}
\end{table}
A pattern is apparent in the coefficients. With $d=2$, the pattern leads us
directly to a general term:
\begin{multline} \label{e:gcd2}
	\mathbb{E}\left(\sigma_{r_{xy}}\right) = \\
 \frac{\rho^{\floor{r_{xy}}}\left(2\pi
    \right)^   {\frac{1}{2}\floor{r_{xy}}}}{\Gamma\left(
\frac{3}{2}\floor{r_{xy}}+1\right)\sqrt{\ceil{r_{xy}}}} \left( \ceil{r_{xy}}    - r_{xy} \right)^{\frac{3}{2}\floor{r_{xy}}}
    + \mathcal{O}\left(\left( \ceil{r_{xy}}
    - r_{xy} \right)^{\frac{1}{2}\left(3\floor{r_{xy}}+2\right)}\right)
\end{multline}
The coefficient gives detailed information about how the number of paths
scales over a unit interval.
\section{Geodesics in $d$-dimensions}
When the vertices are not constrained to lie in $\mathbb{R}^{2}$, which is
realistic since base stations are often also vertically dispersed around
urban areas, we need a high dimensional analysis. The procedure above is
therefore repeated for the
three-dimensional case.

Starting again with $r_{xy} \in \left[0,1\right)$, we immediately have
\begin{equation}\label{e:2mk1}
	\mathbb{E}\left(\sigma_{r_{xy} \in \left[0,1\right)}\right)=1
\end{equation}
For the next interval, in place of $A_{r_{xy}}$ in Eq. \ref{e:int1to2d2}, put
the volume of intersection $V_{r_{xy}}$ of two unit spheres separated by a distance
$r_{xy}$
\begin{equation}\label{e:3lensvolume}
    V_{r_{xy}} = \frac{\pi}{12}\left(4 + r_{xy}\right)\left(2 -
    r_{xy}\right)^{2}
\end{equation}
leaving
\begin{eqnarray} \label{e:term1to2d3}
\mathbb{E}\left(\sigma_{r_{xy} \in \left[1,2\right)}\right)&=& \frac{\pi \rho}{12}\left(4 +
r_{xy}\right)\left(2 - r_{xy}\right)^{2}
\end{eqnarray}
in a similar manner to the case $d=2$. Increasing the distance again, in place of $A_{\lambda}$ in Eq.
\ref{e:int2to3d2} put $V_{\lambda} = \frac{\pi}{12}\left(4 + \lambda\right)\left(2 -
    \lambda\right)^{2}$, and in place of $l_{\lambda}$ we put the surface area $S_{\lambda}$ of the
spherical cap of a sphere-segment, which, omitting details, is
\begin{eqnarray}\label{e:sphericalcap}
S_{\lambda} &=& 2 \pi \lambda^{2} \left(1 - \frac{r_{xy}^2 + \lambda^2 -1}{2
    r_{xy} \lambda}\right)
\end{eqnarray}
The necessary integral is therefore
\begin{multline}
\rho^2\int_{r_{xy}-1}^{2} \left(\frac{\pi}{12}\left(4 + \lambda\right)\left(2
-
\lambda\right)^{2}\right) \left(2 \pi \lambda^{2} \left(1 - \frac{r_{xy}^2 +
\lambda^2 -1}{2
r_{xy} \lambda}\right)\right)
\mathrm{d}\lambda \\
 =\frac{\rho^{2}\pi ^2}{1260}\left(\left(r_{xy}+3\right)
\left(r_{xy}+9\right)-\frac{6}{r_{xy}}\right)\left(3 - r_{xy}\right)^4
\end{multline}
and this is exact, since the volume of intersection of two unit spheres is a
polynomial
in an integer power of the separation of their centres, and so can be
integrated without the need for any expansions. The next term goes according to
\begin{eqnarray}
\mathbb{E}\left(\sigma_{r_{xy} \in
\left[3,4\right)}\right)&=&\rho^3\int_{r_{xy}-1}^{3}S_{\lambda_{1}}\mathrm{d}\lambda_{1}\int_{\lambda_{1}-1}^{2}
    V_{\lambda_{2}} S_{\lambda_{2}} \mathrm{d}\lambda_{2}
\end{eqnarray}
which is
\begin{eqnarray} \label{e:term3to4d3}
\frac{\rho^3 \pi^3}{453,600} \left( \left(2 + r_{xy}\right)\left(8 +
r_{xy}\right)\left(14 + r_{xy}\right) - \frac{144}{r_{xy}}\right)\left(4 -
r_{xy}\right)^{6}
\end{eqnarray}
We then perform the same expansions as before. These expansions are listed in
Table \ref{table:3dintegrals}.
\begin{table}
	\begin{center}
    \begin{tabular}{| l | l | l | l |}
    \hline
Displacement $r_{xy}$ & Expected geodesic cardinality $\sigma_{r_{xy}}$ &    Order of error \\ \hline
    $r_{xy} \in [0,1)$ & 1 & 
	(no error) \\ \hline
    $r_{xy} \in [1,2)$ & $\frac{\pi\rho}{2}\left(2-r_{xy}\right)^{2}$ & 
    $\mathcal{O}\left(\left(2-r_{xy}\right)^{3}\right)$ \\ \hline
    $r_{xy} \in [2,3)$ & $\frac{\pi^2\rho^2}{18}\left(3-r_{xy}\right)^{4}$
    & $\mathcal{O}\left(\left(3-r_{xy}\right)^{5}\right)$ \\ \hline
    $r_{xy} \in [3,4)$ &
    $\frac{\pi^3\rho^3}{360}\left(4-r_{xy}\right)^{6}$ 
    & $\mathcal{O}\left(\left(4-r_{xy}\right)^{7}\right)$ \\ \hline
$r_{xy} \in [4,5)$ & $\frac{\pi^4\rho^4}{12,600}\left(5-r_{xy}\right)^{8}$    & $\mathcal{O}\left(\left(5-r_{xy}\right)^{9}\right)$\\ \hline
$r_{xy} \in [5,6)$ & $\frac{\pi^5\rho^5}{680,400}\left(6-r_{xy}\right)^{10}$    & $\mathcal{O}\left(\left(6-r_{xy}\right)^{11}\right)$ \\ \hline
    \end{tabular}
\end{center}
\caption{Solutions to the recursion relation \ref{e:recursionrelation2d} with
$d=3$.}
\label{table:3dintegrals}
\end{table}

In fact, at this point, we can evaluate the $d$-dimensional term directly. The necessary replacements of
the lens area $A_{\lambda}$ and arc length
$l_{\lambda}$ in $d$-dimensional hyperspace cannot be expressed in terms of
elementary functions. We use the
regularised incomplete beta function $I_{x}\left(a,b\right)$
\begin{equation}
I_{x}\left(a,b\right)\int_{0}^{1}t^{a-1}\left(1-t\right)^{b-1}\mathrm{d}t =    \int_{0}^{x}t^{a-1}\left(1-t\right)^{b-1}\mathrm{d}t
\end{equation}
which is the ratio of the incomplete beta function and the complete beta
function \cite{li2011}.

So, the $d-1$ dimensional surface area of the hyperspherical cap of vertex
angle $\phi \in \left[0,\pi\right]$ and radius $\lambda$ is
\begin{equation} \label{e:hyperarea}
	A_{\lambda} =
\frac{\pi^{d/2}}{\Gamma\left(d/2\right)}\lambda^{d-1}I_{\sin^2{\left(\phi\right)}}\left(\frac{d-1}{2},\frac{1}{2}\right)
\end{equation}
where
\begin{equation}
\phi = \arccos{\left(\frac{\lambda^2+r_{xy}^2-1}{2r_{xy} \lambda}\right)}\end{equation}
is the colatitude at which the $d$-sphere is
cut to produce the hyperspherical segment; cutting the hypersphere at $\phi =
\pi/2$ would produce a hemisphere, for example, when $d=3$. $A_{\lambda}$
reduces to an arc length in Eq. \ref{e:arclength} when $d=2$, and the
two-dimensional
surface area of a spherical cap (Eq. \ref{e:sphericalcap}) when $d=3$
\cite{li2011}.

The volume of intersection of two unit hypersphere separated by
$r_{xy}$ is twice the \textit{volume} of this hyperspherical cap (each are
glued to the plane
which cuts through the sphere-sphere intersection). Thus
\begin{equation} \label{e:hypervolume}
    V_{\lambda} =
    \frac{\pi^{d/2}}{\Gamma\left(d/2 +
1\right)}\lambda^{d}I_{\sin^2{\left(\phi\right)}}\left(\frac{d+1}{2},\frac{1}{2}\right)
\end{equation}
which reduces to the area of the lens when $d=2$ (Eq. \ref{e:lensarea}), and
the volume of the sphere-sphere
intersection when $d=3$ (Eq. \ref{e:3lensvolume}). As before, we have
\begin{eqnarray}
\mathbb{E}\left(\sigma_{r_{xy} \in
\left[2,3\right)}\right)&=&\int_{r_{xy}-1}^{2}\rho V_{\lambda}
    \rho A_{\lambda} \mathrm{d}\lambda \label{e:5}
\end{eqnarray}
which is
\begin{eqnarray} \label{e:int2to3dd}
    \frac{\pi^{d}}{\Gamma\left(d/2\right)\Gamma\left(d/2 +
    1\right)} \int_{r_{xy}-1}^{2}
\lambda^{2d-1}I_{\sin^2{\left(\phi\right)}}\left(\frac{d+1}{2},\frac{1}{2}\right)I_{\sin^2{\left(\phi\right)}}\left(\frac{d-1}{2},\frac{1}{2}\right)\mathrm{d}\lambda
\nonumber
\end{eqnarray}
which is intractable, but we can extract the leading order term.

So, in place of $A_{r_{xy}}$ in Eq. \ref{e:int1to2d2}, put the volume of intersection
$V_{r_{xy}}$ of two unit hyperspheres separated by a distance $r_{xy}$
\begin{equation}\label{e:3lensvolume}
    V_{r_{xy}} = \frac{\pi^{d/2}}{\Gamma\left(d/2 +
1\right)}r_{xy}^{d}I_{\sin^2{\left(\phi\right)}}\left(\frac{d+1}{2},\frac{1}{2}\right)
\end{equation}
leaving
\begin{eqnarray} \label{e:term1to2d3}
\mathbb{E}\left(\sigma_{r_{xy} \in \left[1,2\right)}\right)&=&\rho V_{r_{xy}}
\nonumber \\
&=& \frac{\pi^{d/2}\rho}{\Gamma\left(d/2 +1\right)}r_{xy}^{d}I_{\sin^2{\left(\phi\right)}}\left(\frac{d+1}{2},\frac{1}{2}\right)
\end{eqnarray}
in a similar manner to before. Increase the distance again, and in place of $A_{\lambda}$ in Eq.
\ref{e:int2to3d2}, put $V_{\lambda}$
\begin{equation}
    V_{\lambda} = \frac{\pi^{d/2}}{\Gamma\left(d/2 +
1\right)}\lambda^{d}I_{\sin^2{\left(\phi\right)}}\left(\frac{d+1}{2},\frac{1}{2}\right)
\end{equation}
and in place of $l_{\lambda}$ put the surface area of the hyperspherical cap
$A_{\lambda}$ with
\begin{eqnarray}
	 \phi = \arccos{\left(\frac{\lambda^2+r_{xy}^2-1}{2r_{xy} \lambda}\right)}
\end{eqnarray}
as before. For reference this is
\begin{equation} \label{e:hyperarea2}
	A_{\lambda} =
\frac{\pi^{d/2}}{\Gamma\left(d/2\right)}\lambda^{d-1}I_{\sin^2{\left(\phi\right)}}\left(\frac{d-1}{2},\frac{1}{2}\right)
\end{equation}
so since (as before) we have
\begin{eqnarray}
\mathbb{E}\left(\sigma_{r_{xy} \in
\left[2,3\right)}\right)&=&\int_{r_{xy}-1}^{2}\rho V_{\lambda}
    \rho A_{\lambda} \mathrm{d}\lambda \label{e:int2to3dd}
\end{eqnarray}
then we need
\begin{eqnarray}
    \frac{\pi^{d}}{\Gamma\left(d/2\right)\Gamma\left(d/2 +
    1\right)} \int_{r_{xy}-1}^{2}
\lambda^{2d-1}I_{\sin^2{\left(\phi\right)}}\left(\frac{d+1}{2},\frac{1}{2}\right)I_{\sin^2{\left(\phi\right)}}\left(\frac{d-1}{2},\frac{1}{2}\right)\mathrm{d}\lambda
\nonumber
\end{eqnarray}
which was given explicitly in Eq. \ref{e:int2to3dd}. Before (when this was
intractable) we expanded $A_{\lambda}$ at $\lambda=r_{xy}-1$ and $V_{r_{xy}}$
(which is the previous term in the sequence) at $r_{xy}=\ceil{r_{xy}}=2$ (and
then replace $r_{xy}$ with $\lambda$). Given
\begin{equation}
A_{\lambda} = \frac{\left(2
\pi\right)^{\frac{d-1}{2}}r_{xy}\left(r_{xy}-1\right)^{d}}{\Gamma{\left(\frac{1}{2}\left(d+1\right)\right)}}\left(1-r_{xy}+\lambda\right)^{\frac{d-1}{2}}
+ \mathcal{O}\left(\left(\lambda -
\left(r_{xy}-1\right)\right)^{\frac{d+1}{2}}\right) \nonumber
\end{equation}
and
\begin{equation}
V_{\lambda} =
\frac{\pi^{\frac{d-1}{2}}\rho}{\Gamma{\left(\frac{1}{2}\left(d+3\right)\right)}}\left(2-\lambda\right)^{\frac{d+1}{2}}
+ \mathcal{O}\left(\left(\lambda-2\right)^{\frac{d+3}{2}}\right) \nonumber
\end{equation}
the integral in Eq. \ref{e:int2to3dd} evaluates to	
\begin{equation}
\mathbb{E}\left(\sigma_{r_{xy} \in \left[2,3\right)}\right) =
\frac{\left(2\pi\right)^{d-1}3^{\frac{1}{2}\left(1-d\right)}}{\Gamma{\left(2+d\right)}}\left(3-\lambda\right)^{d+1}
+ \mathcal{O}\left(\left(\lambda-3\right)^{d+2}\right)
\end{equation}
after expanding as usual. The rest of these integrals are displayed in Table
\ref{table:ddimensionintegrals}.
\begin{table}
		\begin{center}
    \begin{tabular}{| l | l | l | l |}
    \hline
Displacement $r_{xy}$ & Expected geodesic cardinality $\sigma_{r_{xy}}$ &    Order of error \\ \hline
    $r_{xy} \in [0,1)$ & 1 & 
	(no error) \\ \hline
    $r_{xy} \in [1,2)$ & 

$\frac{\rho\left(2\pi\right)^{\frac{1}{2}\left(d-1\right)}2^{\frac{1}{2}\left(1-d\right)}}{\Gamma{\left(\frac{3+d}{2}\right)}}\left(2-r_{xy}\right)^{\frac{1}{2}\left(d+1\right)}$
&
	
$\mathcal{O}\left(\left(r_{xy}-2\right)^{\frac{1}{2}\left(d+3\right)}\right)$
\\ \hline

$r_{xy} \in [2,3)$ &
$\frac{\rho^{2}\left(2\pi\right)^{d-1}3^{\frac{1}{2}\left(1-d\right)}}{\Gamma{\left(\frac{1}{2}\left(2+d\right)\right)}}\left(3-r_{xy}\right)^{d+1}$

	& $\mathcal{O}\left(\left(r_{xy}-3\right)^{d+2}\right)$ \\ \hline

$r_{xy} \in [3,4)$ &
$\frac{\rho^{3}\left(\pi\right)^{\frac{3}{2}\left(d-1\right)}2^{\frac{1}{2}\left(d-1\right)}}{\Gamma{\left(\frac{1}{2}\left(5+3d\right)\right)}}\left(4-r_{xy}\right)^{\frac{3}{2}\left(d+1\right)}$

&
$\mathcal{O}\left(\left(r_{xy}-4\right)^{\frac{3}{2}\left(d+3\right)}\right)$
\\ \hline

    \end{tabular}
\end{center}
\caption{Solutions to the recursion relation \ref{e:recursionrelation2d} for
general dimension $d$.}
\label{table:ddimensionintegrals}
\end{table}
and, finally, by inspecting a pattern in the coefficients, we have the
expected number of $\ceil{r_{xy}}$ paths in $d$-dimensions.
\begin{multline}\label{e:final}
	\mathbb{E}\left(\sigma_{r_{xy}}\right)=
		\frac{\rho^{\floor{r_{xy}}}
		\left(2\pi\right)^{\frac{1}{2}\floor{r_{xy}}\left(d-1\right)}
\ceil{r_{xy}}^{\frac{1}{2}\left(1-d\right)}}{\Gamma\left(\frac{ \ceil{r_{xy}}
+1}{2} + \frac{\floor{r_{xy}}d}{2}\right)}
\left(\ceil{r_{xy}}-r_{xy}\right)^{\frac{1}{2}\floor{r_{xy}}\left(d+1\right)}\end{multline}
with error
$\mathcal{O}\left(\left(\ceil{r_{xy}}-r_{xy}\right)^{\frac{1}{2}\floor{r_{xy}}\left(d+3\right)}\right)$.

\section{Numerical simulation}

\begin{figure}[t]
\hspace{1.5mm}\includegraphics[scale=0.5]{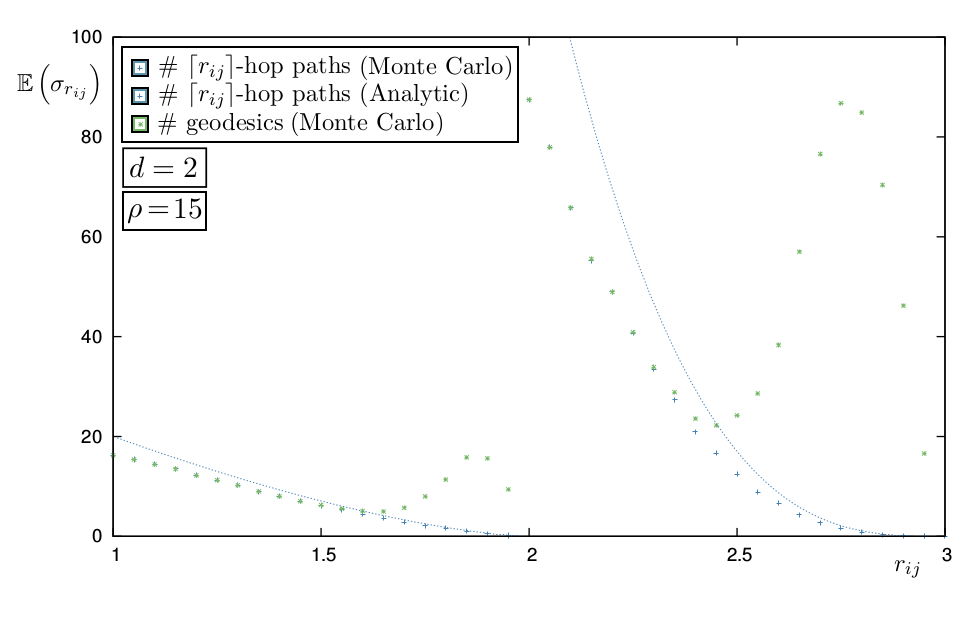}
\caption{The smooth decaying blue lines are Monte Carlo data for
\textit{optimal geodesic cardinality}, with our approximation Eq.
\ref{e:final} with $d=2$. The green stars are the actual shortest path
counts. The blue pluses and green stars are numerical results, while the
smooth curves are our analytic results.}
\label{fig:all}
\end{figure}

\begin{figure}
\includegraphics[scale=0.7]{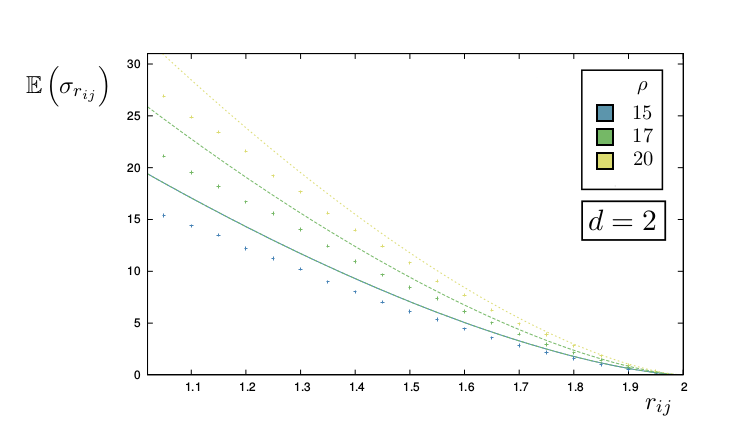}
\caption{The expected number of two-hop paths and our approximation Eq. \ref{e:final}, for $\rho=15, 17,
20$ and $d=2$.}
\label{fig:optimalone}
\includegraphics[scale=0.7]{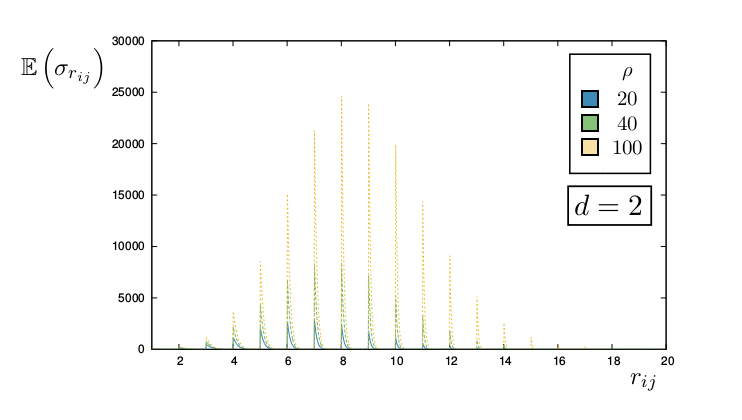}
\caption{Eq. \ref{e:final} plotted for three increasing densities
$\rho=20,40,100$ over a long range of displacements. The expected number of
paths rises and falls as $r_{xy}$ increases.}
\label{fig:updown}
\end{figure}

Fig. \ref{fig:optimalone} gives an idea of the error introduced by the Taylor expansions we perform. Performing these integrals exactly is topic of current research. Notice also Fig. \ref{fig:updown}, which shows how Eq. \ref{e:final} demonstrates a
rising-and-falling behaviour as $r_{xy}$ increases. This occurs because there
are two competing forces which `create and destroy' geodesic cardinality. The
creation force involves increasing $r_{xy}$, which allows more space for geodesics to form. 
This force dominates for small $r_{xy}\leq8$ with respect to the
parameters in Fig. \ref{fig:updown}. The destructive force is caused by the
increasingly difficult but necessary alignment of vertices in the various
lens' which must be achieved for optimal geodesics to form.
This force is dominating for $r_{xy}\geq8$, since more than eight lenses must
accommodate aligning vertices for optimal paths to form.

Fig. \ref{fig:all} shows how the number of geodesic paths is asymptotic to
the number of $\ceil{r_{xy}}$ hop paths as $\rho \to \infty$. This
approximation breaks down as $r_{xy} \to \ceil{r_{xy}}$, because the lens
areas are expected to contain fewer and fewer vertices, and so the geodesic
can become longer than $\ceil{r_{xy}}$ hops.

\chapter{Discussion}\label{chapter:discussion}
In this chapter, we discuss the relationship between our mathematical
results, and communications practice. We also discuss possible future
developments of this research, and how it may have a lasting impact on future
multi-hop communication.
\section{Applications in communication networks}
We focus in this section on detailing the nature of our contributions to the field of wireless networks.

\subsection{Understanding skeletons}\label{sec:skeleton}
In Fig \ref{fig:skeleton}, two example \textit{skeletons} have been extracted from random geometric graphs bound inside
obstructed domains \cite{liu2015}. In Chapter \ref{chapter:betweenness}, we defined a skeleton from its use in e.g. Ercsey-Ravasz et al. \cite{ercsey2010}, who demonstrated how betweenness centrality 
can be used to describe and analyse a network's \textit{skeleton} or
\textit{vulnerability backbone}. This is a percolating
cluster of the most structurally important vertices. See also \cite{liu2015} for more recent work. These appear to be quasi-one-dimensional sub-systems, which are subgraphs embedded in a $d$-dimensional space which extend in only a single dimension such as around the perimeter of a disk or along a corridor region. They may have a restricted extension into the remaining $d-1$ dimensions, hence the use of the phrase \textit{quasi}. Thin systems, such as thin films, have been important in physics more generally \cite{ohring2001}.
In our case, these systems appear to have control over the high performance functionality of the network when contained within bounding geometries with non-convex features. Questions therefore emerge, such as:
\begin{enumerate}
\item How does the connectivity of this sub-graph effect the overall network
performance?

\item How can this extra connectivity be balanced with minimisation of
inter-skeleton interference?

\item Do multiple small obstructions create multiple one-dimensional
sub-systems which all need to be managed?

\item Overall, what is the best way skeleton vertices can be used to optimise
communication?
\end{enumerate}

\subsection{Topology-based geolocalisation}
There is greater sensitivity in the \textit{number} of paths than in the
\textit{length} of paths to small changes in location in random spatial networks.
\noindent As an enhancement of topological-based localisation, where hop counts to location-aware
anchor vertices are converted, preferably analytically, into likely Euclidean
distances, the number of paths can be used as an accuracy enabler \cite{zhang2012,mao2010,nguyen2015}.

This is particularly useful in dense networks, since here the number of hops
from $x$ to $y$ is determined by the ceiling of Euclidean separation, and so does not change over potentially an entire connection
interval. The number of paths, however, does scale, and moreover,
with an increasing resolution at high density. A combination of path length, and path counts, could provide a key solution to localisation in dense networks beyond GPS. See e.g. the landmark paper on DV-hop \cite{niculsecu2003}, which approximates distances via the average distance of a hop multiplied by the number of hops between a pair of nodes.

\subsection{Future vehicular ad hoc networks}
Betweenness centrality can be more easily evaluated if the domain is a one
dimensional strip, such as the interval $\left[0,L\right)$. The combinatorial
enumeration of paths, and evaluation of the lengths of paths, is
significantly less involved. Yet the applications are perhaps more direct, as ad hoc vehicle
networks are among the most promising applications of multi-hop
routing, and key development area of future wireless networks \cite{gerla2014}.

Possible future research topics include:
\begin{enumerate}
\item Can current flow betweenness centrality be expressed as a one-dimensional field of random variables, which is a random field \cite{clifford1990}, in a
one-dimensional random geometric graph on a finite interval of the real line?

\item What is the distribution of geodesic paths between two arbitrarily
separated vertices in a one-dimensional random geometric graph? Can this be used to evalaute expected location-specific betweenness in closed form in one dimension?
\end{enumerate}

\begin{figure}[t]
\begin{center}
\includegraphics[scale=0.5]{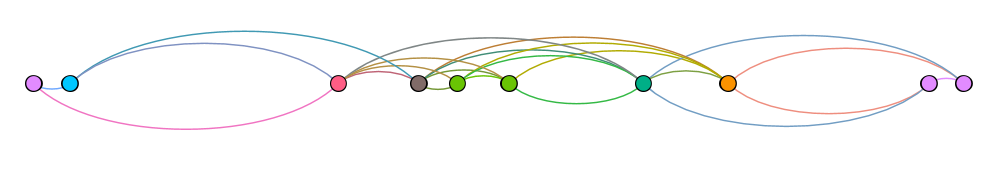} 
\caption{A one-dimensional unit disk graph. The colours are used to help
distinguish paths. The geodesic length between the end points is three
hops. This can be used to model an ad hoc communication network of vehicles
queing in traffic.}\label{fig:discussionskeleton}
\end{center}
\end{figure}

\textit{Distributed algorithms} run on a network in such a way that the resources required for their operation are shared equally between all cooperating devices. This distribution of resource consumption can prove important for battery limited sensor networks, even if they have access to the resources of a central server such as a sink or cluster head \cite{hua2016}.

Analysis of expected centrality as a scalar field (defined on the $d$-dimensional domain positions), as discussed in this thesis, works towards understanding the best way to distribute the necessary computational load of such algorithms. Continued work on this topic, such as working to see how centrality scales near obstacles in the domain, will also highlight the effect that non-convex features can have on the long term performance of these distributed network processes, such as in sensing scenarios.

\section{Mathematical directions}
Below we highlight some research directions which are more directly related to mathematical directions in random networks.
\subsection{The distribution of the number of geodesics}
\begin{figure}
\noindent \begin{centering}
\includegraphics[scale=0.55]{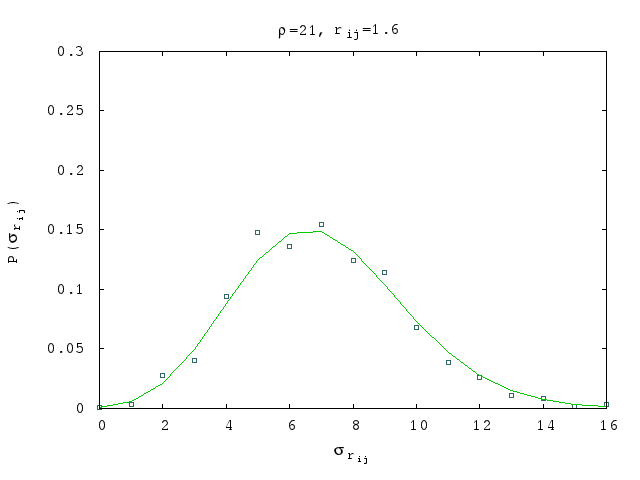} \\
\hspace{24mm}\includegraphics[scale=0.55]{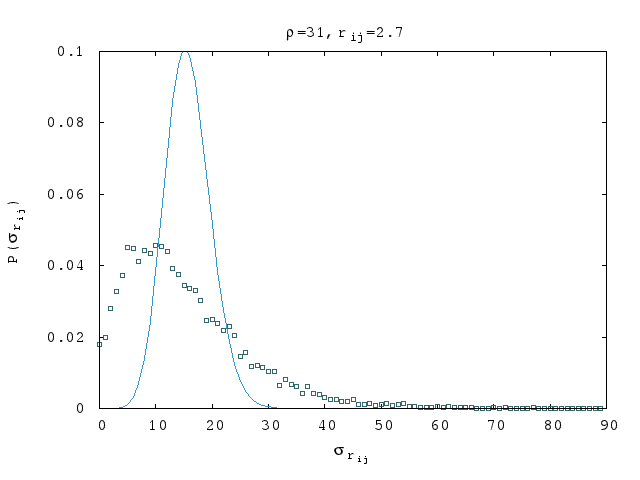}\end{centering}
\caption{Top: The small boxes are Monte Carlo approximations to the
probability $P\left(\sigma_{r_{xy}}=k\right)$ with
$r_{xy}=1.6$. In this distance interval the distribution of geodesic paths is Poisson, an described by the green line Eq. \ref{e:poisson}. Bottom: The blue line is the curve
$\text{Po}\left(\mathbb{E}\left(\sigma_{r_{xy}=2.7}\right)\right)$, and the
boxes Monte Carlo data for this larger displacement. The number of
geodesic paths is therefore not Poisson.}
\label{fig:poissondistribution}
\end{figure}
What is the distribution of shortest paths between two vertices, now considering $d=2$? The first
moment of $\sigma_{r_{xy}}$ is expected number of $\ceil{\norm{x-y}}$-hop
paths. It is interesting to ask if we can get the other moments. Note that
\textit{if} the distribution of paths is Poisson, the second moment is
related to the first in a trivial way
\begin{eqnarray}\label{e:secondmomentdefinition}
\mathbb{E}\left(\sigma^2_{r_{xy}}\right) =
\mathbb{E}\left(\sigma_{r_{xy}}\right)\left(\mathbb{E}\left(\sigma_{r_{xy}}\right)+1\right)
\end{eqnarray}
and $\sigma_{r_{xy}}$ is distributed as	
\begin{eqnarray}\label{e:poisson}
P\left(\sigma_{r_{xy}}=k\right) =
\frac{1}{k!}\exp{\left(-\mathbb{E}\left(\sigma_{r_{xy}}\right)\right)}
\left(\mathbb{E}\left(\sigma_{r_{xy}}\right)\right)^{k}
\end{eqnarray}
This actually \textit{does} hold when $r_{xy} \in \left[1,2\right)$, as
demonstrated in the top panel of Fig. \ref{fig:poissondistribution}, but not
for larger displacements, as demonstrated in the bottom panel of Fig
\ref{fig:poissondistribution}. For these larger displacements, the variance 
\begin{eqnarray}
\mathrm{Var}\left(\sigma_{r_{xy}}\right) =
\mathbb{E}\left(\sigma^2_{r_{xy}}\right) -
\mathbb{E}^2\left(\sigma_{r_{xy}}\right)
\end{eqnarray}
exceeds the mean by an amount proportional to both density $\rho$ and
displacement $r_{xy}$, and so $\sigma_{r_{xy}}$ is no longer Poisson, but
conjectured \textit{negative binomial}.
\begin{conjecture}[Distribution of the number of geodesic
paths]\label{c:negativebinomial}
Given two vertices displaced by $r_{xy}$ in the unit disk graph $G\left(n,\pi
\right)$, the random variable $\sigma_{r_{xy}}$ is distributed as
	\begin{eqnarray}
P\left(\sigma_{r_{xy}}=k\right) = {k + r - 1 \choose
k}p^{r}\left(1-p\right)^{k}
	\end{eqnarray}
with $p,r$ the solution to the simultaneous equations
\begin{eqnarray}
\mathbb{E}\left(\sigma_{r_{xy}}\right)&=&\frac{(1-p)r}{p} \\
\mathbb{E}\left(\sigma^2_{r_{xy}}\right) &=&
\mathbb{E}\left(\sigma_{r_{xy}}\right)\left(\frac{1}{p} +
\mathbb{E}\left(\sigma_{r_{xy}}\right)\right)
\end{eqnarray}.
\end{conjecture}
\noindent An analytic expression for
$\mathbb{E}\left(\sigma^2_{r_{xy}}\right)$ would allow one to numerically
verify this. The basis of the conjecture is mainly the similarity in form of
the probability mass functions, and also the apparent over-dispersion (i.e.
the variance of $\sigma_{r_{xy}}$ is larger than its mean). One would need to
investigate this further as it is not understood.

\subsection{Non-optimal geodesics}
Some paths are geodesics of length $k > \ceil{r_{xy}}$ hops. These
occur when there are no paths of the optimal length. The shortest path is
then found to be one hop longer. We call geodesics of length
$\ceil{r_{xy}}+1$ hops \textit{$\beta$-optimal} paths, and those of length
$\ceil{r_{xy}}+2$ hops \textit{$\gamma$-optimal}, and so on. 
\begin{figure}[!]
\noindent \begin{centering}
\includegraphics[scale=0.8]{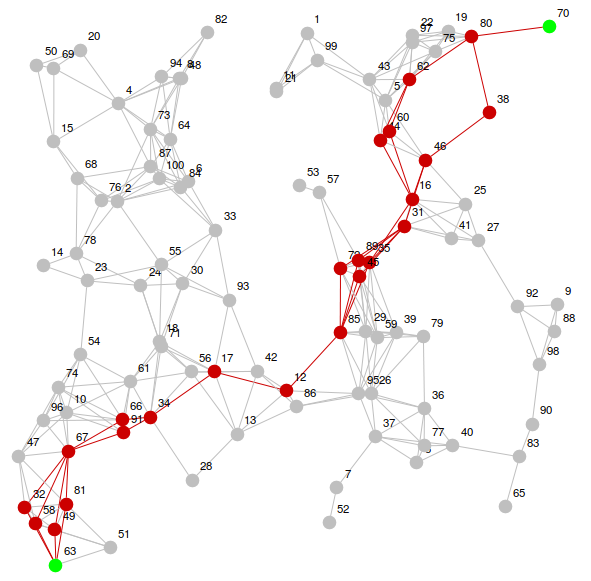}
\par\end{centering}
\caption{The graph shows $96$ geodesics, each of $13$-hops, from vertex 63 to
vertex 70. Since the Euclidean distance between the endpoints is 9.8, these
geodesics are not optimal. As $r_{xy}$
increases, the probability there existing at least one 'optimal' alignment of vertices, which is an alignment of points of $\mathcal{Y}$ that induces a path of $\ceil{r_{xy}}$ hops from $x$ to $y$, vanishes. The chain of links is more likely to break as it grows longer.
This explains the rising and falling effect in Fig. \ref{fig:updown}.}
\label{fig:eighthops}
\end{figure}
The expected number of geodesics with $r_{xy}\in\left[1,2\right)$ is
therefore an approximation to
\begin{eqnarray}\label{e:betaoptimal}
\mathbb{E}\left(\sigma_{r_{xy}\in\left[1,2\right)}\right) = \rho A_{r_{xy}}
(1-e^{-\rho A_{r_{xy}}}) + e^{-\rho A_{r_{xy}}}\left(\rho^2 \int_{1}^{2}
l_\lambda A_\lambda \mathrm{d}\lambda\right)
\end{eqnarray}
with the number of $\gamma$-optimal paths a second order correction.
Investigating this decomposition of geodesic paths into contributions from
paths of all lengths is a topic of further research.

In particular, 

\begin{enumerate}
\item If each edge in a path is likened to a sequence of unit resistors in
series (in an electronic circuit), and all paths of the same length in hops
are considered to be like strings of resistors running in parallel, what is
the resistance of the single resistor which can replace the set of paths
which join $x$ and $y$? There will be a distribution on the possible
resistances, but interestingly, this will be determined in the limit of high
density by only the distribution of geodesics.

\item How successful is the modelling of data flows as a flow of current?
What are the key differences?
\end{enumerate}

\noindent We hope to develop the first of these goals as part of a project on
\textit{random geometric electrical networks}.

\section{Final words}
The potential ultra-dense base station deployments of 5G have been
investigated from the point of view of random geometric graphs, and
specifically betweenness centrality, as well as part of the theory of ad hoc wireless networks in general. It has been shown how spatial
probabilistic combinatorics can play a key role in effective, common sense
wireless interference management techniques, as well as various optimisation methods, such as dynamically tuning the betweenness of vertices, and bounding the length of geodesic paths. This enables delay tolerant networking,
involving sociological ideas such as betweenness centrality, and current
flow betweenness centrality. The future tractability of intelligent vehicles is also
supported by the greater ease of analysis in one-dimensional and quasi-one dimensional random geometric graphs, compared with their
higher-dimensional counterparts.

In the future, the issue of spatial dependence needs to be addressed in order to provide a greater
depth to the applied mathematics of these fascinating graphs. The potential limiting performance
of ad hoc communication networks is not yet understood, but with this sort of graph
theoretic analysis, it may yet be harnessed.
\begin{appendices}

\chapter{Percolation} \label{appendix:percolation}
\begin{figure}
\noindent
\begin{centering}
	\includegraphics[scale=0.17]{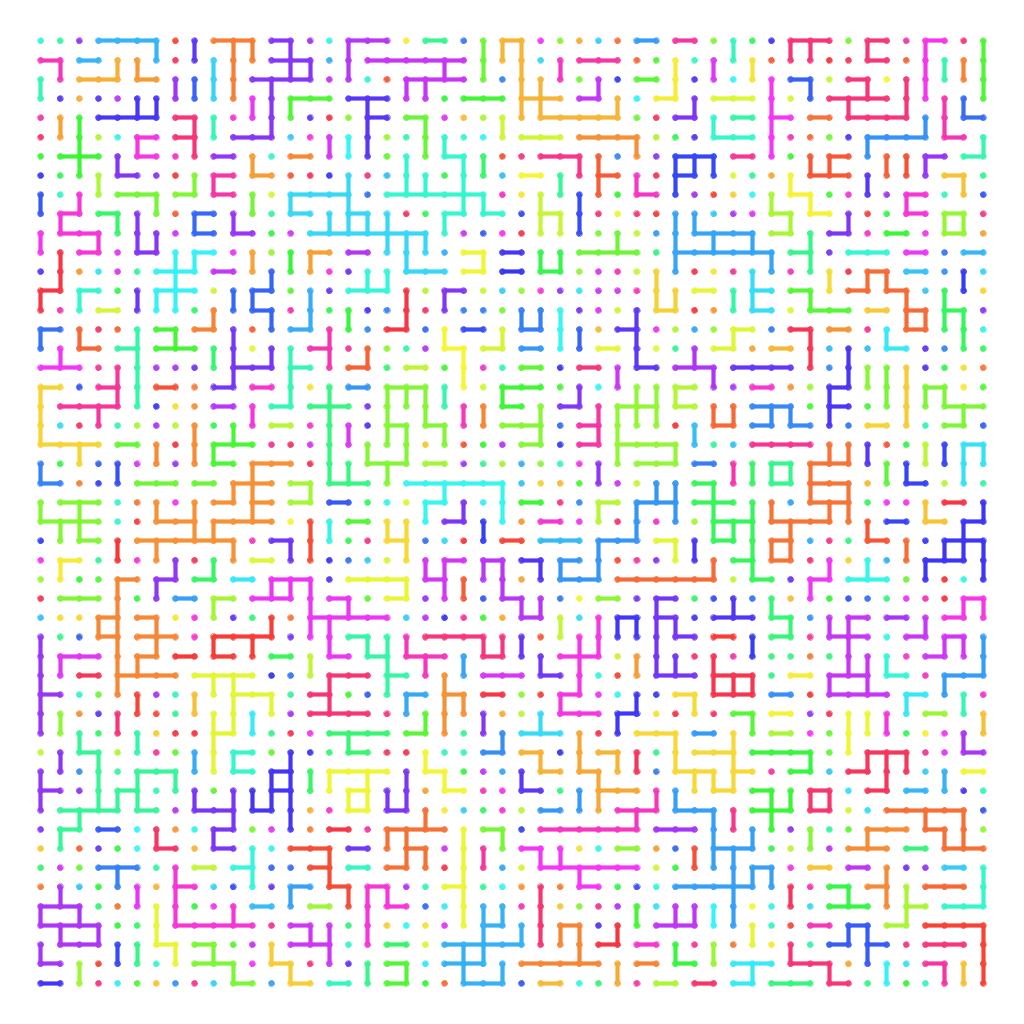}\includegraphics[scale=0.17]{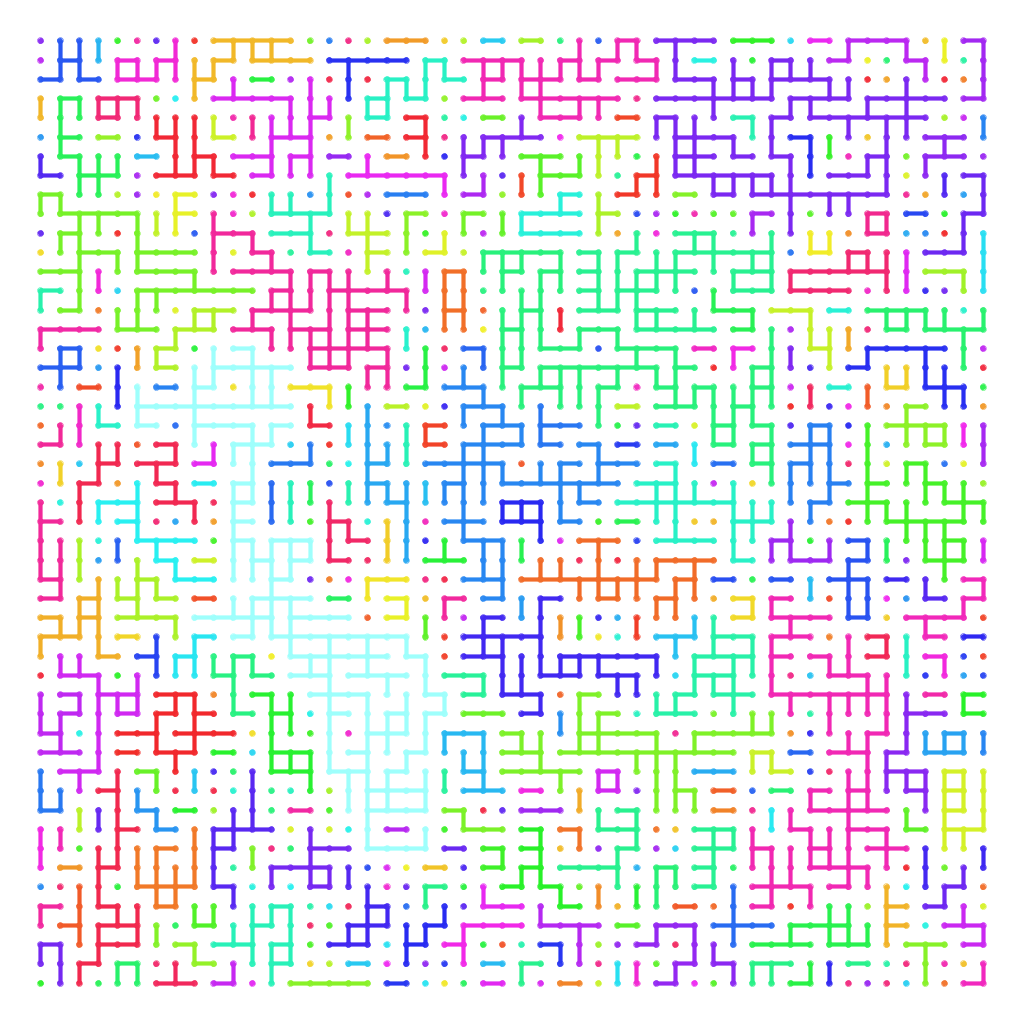} \\
	\includegraphics[scale=0.17]{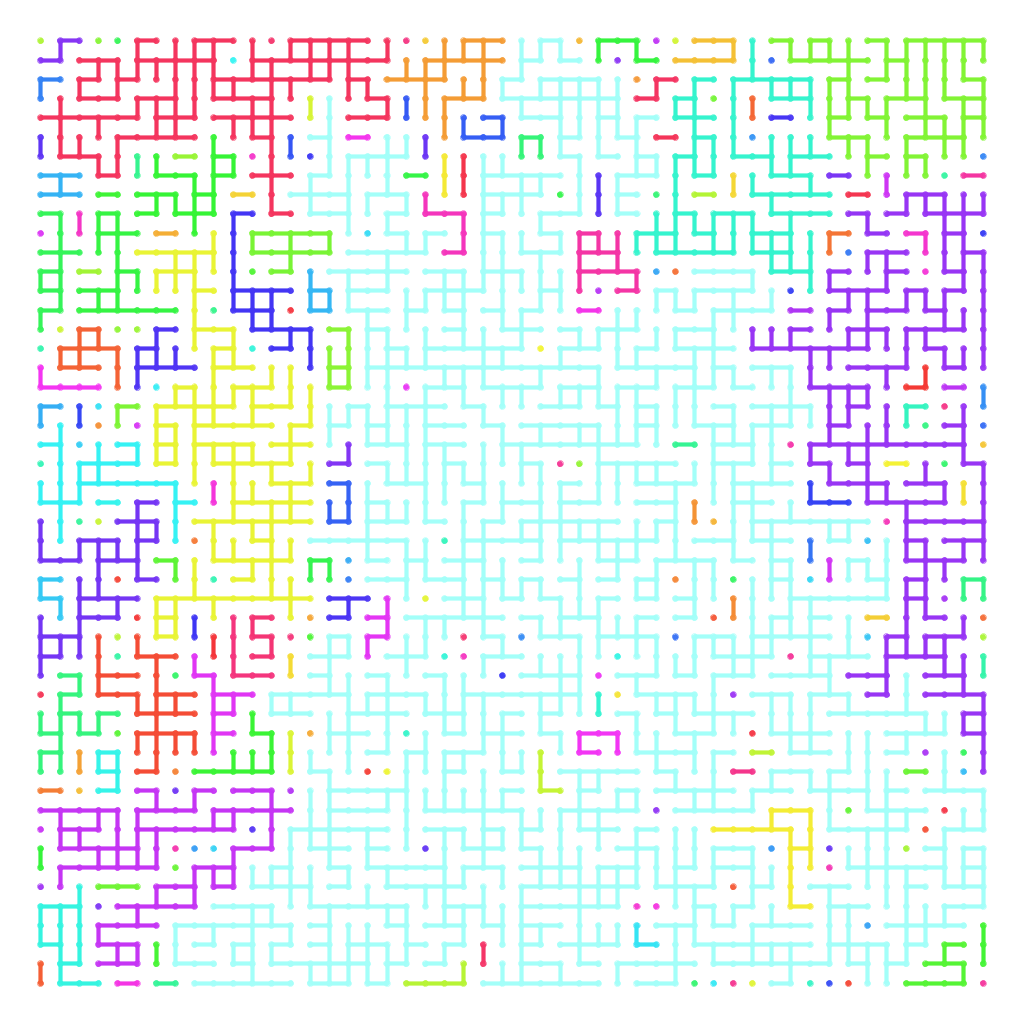}\includegraphics[scale=0.17]{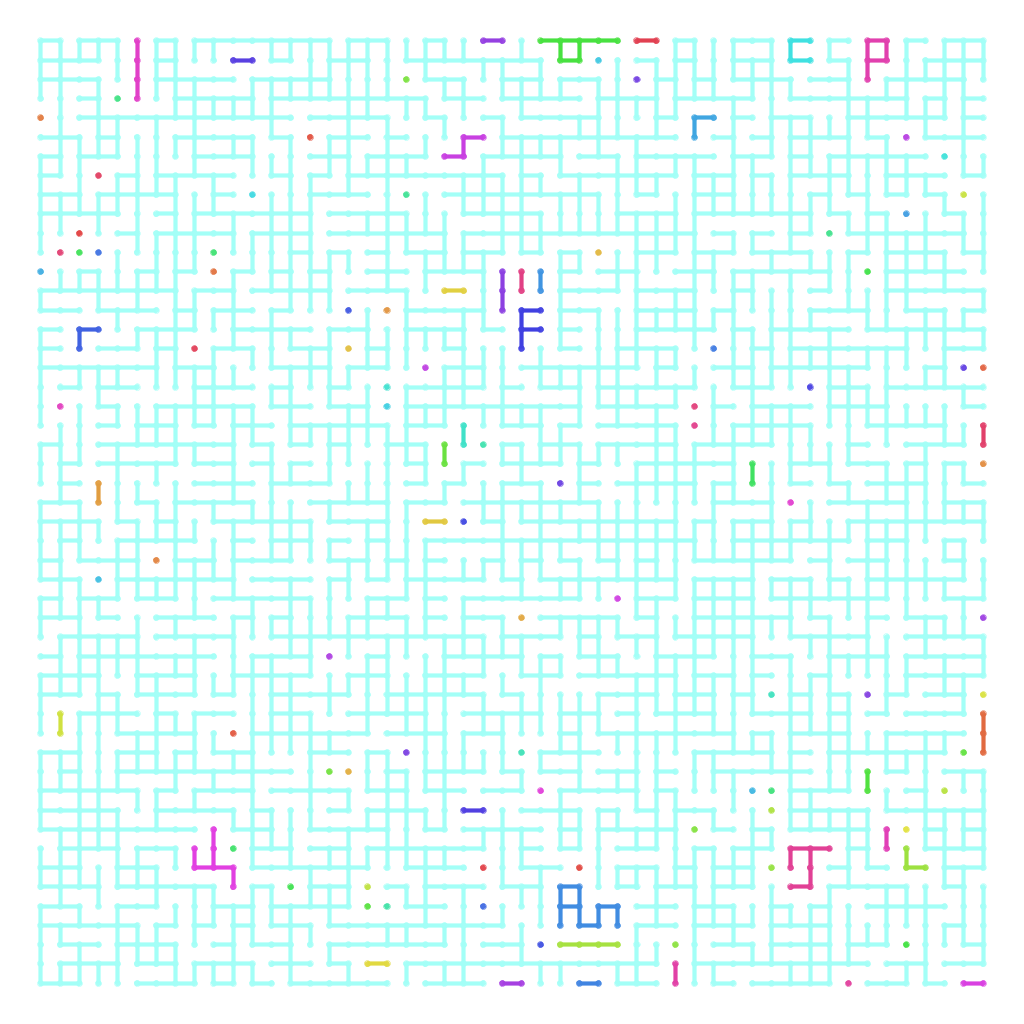} \\
	\includegraphics[scale=0.17]{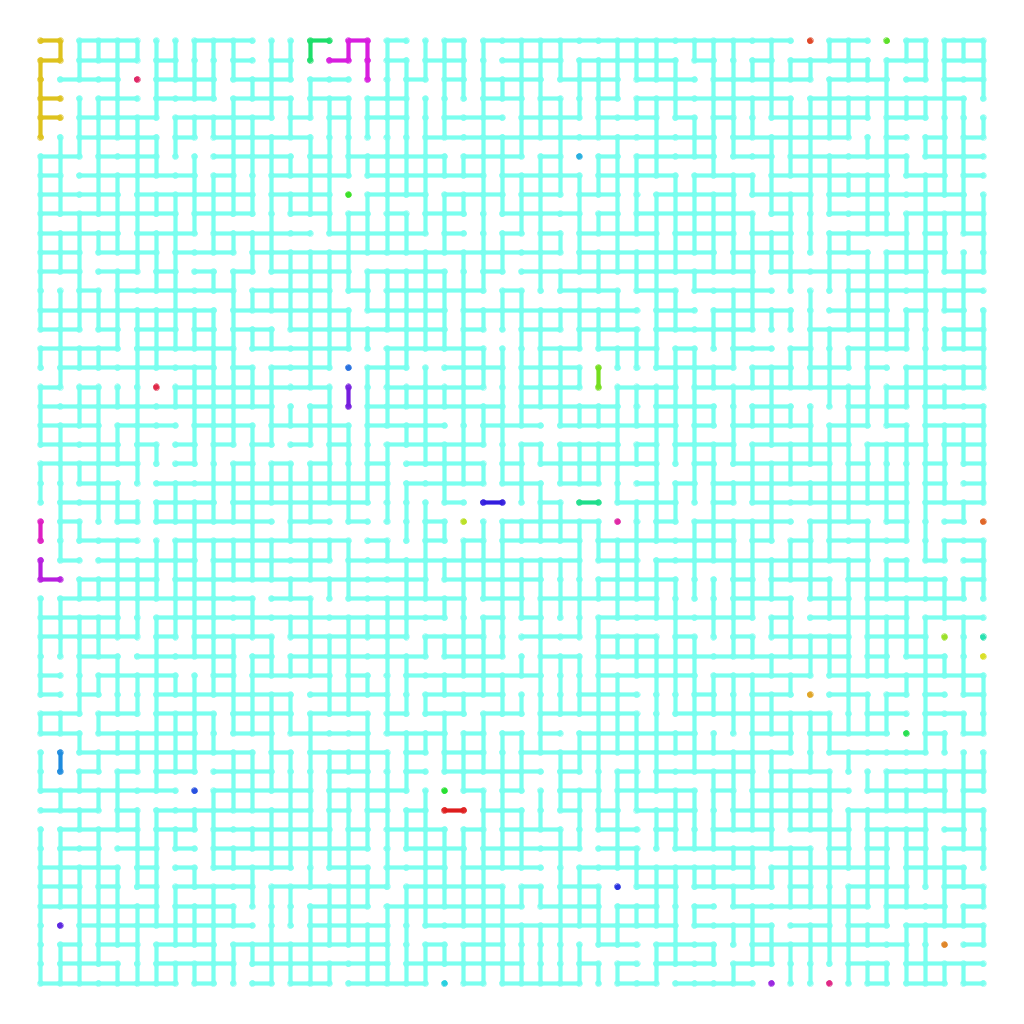}\includegraphics[scale=0.17]{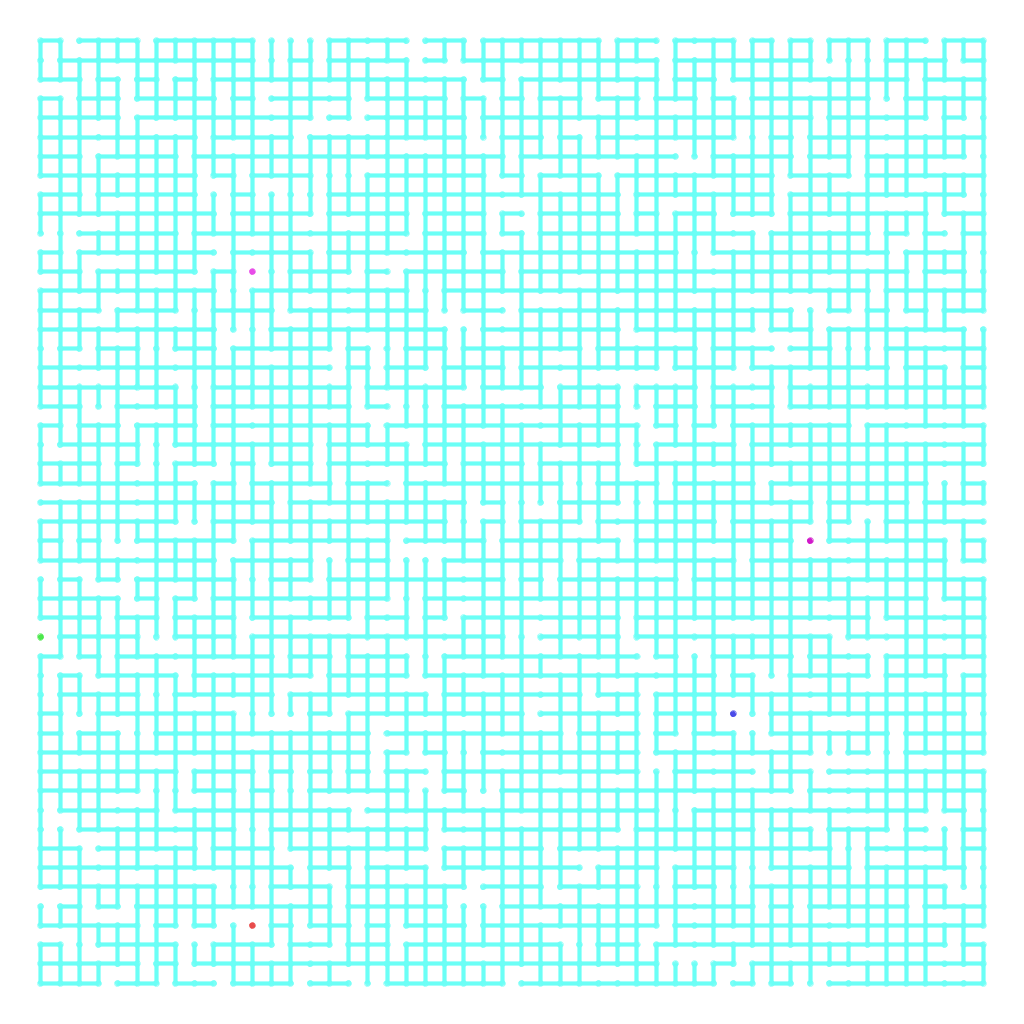} \\
\end{centering}
\caption{$50 \times 50$ bond percolation with $p$ increasing from 0.3 to 0.8 in steps of 0.1. Criticality occurs at $p=0.5$. The largest connected component is highlighted in sky blue throughout.}
\label{fig:colourpercolation}
\end{figure}
Percolation will now be more rigorously defined. Take the hypercubic lattice $\mathcal{D} = \mathbb{L}^d$. Look at all the edges of the lattice, and designate them either \textit{open} or \textit{closed} with probability $p$. This is called \textit{bond percolation}\footnote{Site percolation, as an alternative, turns the lattice squares themselves on or off with probability $p$.}. The question is, what proportion of the lattices display a path of open links from the origin $\mathcal{O}$ to the boundary of the lattice $\partial\mathcal{D}$, or to infinity if we consider an infinite lattice? There is a very well studied singularity at $p=p_c$ in both 
\begin{eqnarray}
	\theta(p) = \mathbb{P}(\singlenorm{C}=\infty),\qquad \chi(p) = \mathbb{E}_{p}\singlenorm{C}
\end{eqnarray}
where $C$ is the set of all vertices connected to the origin via open paths, and $\chi(p)$ is the mean cluster size. This singularity is due to either an undefined derivative at $\theta(p_c)$, or an undefined magnitude at $\chi(p_c)$, and is due to conjectured fractional power-law behaviour. 

\begin{theorem}[There exists a non-trivial critical point for $d\geq2$] We have that
	\begin{equation}\label{e:sGs}
		0 < p_c < 1 \qquad \text{if } \quad d \geq 2
	\end{equation}
\end{theorem}
Notice that $\theta\left(0\right)=0$, since there are no links. Also $\theta\left(1\right)=1$, since all vertices are linked. This theorem is the existence of a non-trivial critical parameter $p_c \in (0,1)$, and immediately introduces the idea of a self-avoiding random walk (SAW).

SAW is necessarily a sequence of adjacent vertices which does not repeat\footnote{Though this is not the exact definition, it is suitable here.}. It is thus a path through the lattice, of length $\nu$. Call $\Sigma_n$ the number of SAWs of length $n$ which leave from the origin (or, equivalently, finish at the origin). Let $\Sigma_n^{\star}$ be the number of these walks which only consist of open links. This is related to the percolation probability
\begin{equation}
	\theta(p) = \lim_{n \to \infty}P_{p}\left(\Sigma_n^{\star} \geq 1\right)
\end{equation}
Considering first $\Sigma_n$, we have that (at worst)
\begin{equation}
	\Sigma_n \leq 2d\left(2d - 1\right)^{n-1}
\end{equation}
given a symmetric random walk \cite{grimmettbook} from the origin which takes no steps back along its current path. Also
\begin{equation}
	P_{p}\left(\Sigma_n^{\star} \geq 1\right) \leq \mathbb{E}\left(\Sigma_n^{\star}\right) = p^n\Sigma_n
\end{equation}
where we note that the expected number of SAWs along open links is not a trivial object, but can be written in terms of the total number of SAWs of length $n$. Thus
\begin{equation}
	\theta(p) \leq \lim_{n \to \infty} 2d\left(2d - 1\right)^{n-1}p^n
\end{equation}
implying $\theta(p)$ can only be non-zero when
\begin{equation}
	p_c \geq \frac{1}{2d-1}
\end{equation}
demonstrating the lower bound for $p_c$, since $\theta(p)$ is non-decreasing in $p$, and equals unity when $p=1$.

For the upper bound, the idea is to create a dual lattice identical to the main lattice (in both geometry and open edge configuration), but spatially translated half a link-length south, and half a link-length to the east. Any finite cluster of open links will induce the dual lattice to display a finite loop of \textit{closed} links which encloses this open cluster. 

Event $M_n$ is a completed loop of length $n$. We have
\begin{equation}
	P_{p}\left(\sum_{n}M_n\geq1\right) \leq \mathbb{E}_{p}\left(\sum_{n}M_n\geq1\right)
\end{equation}
after considering the usual power series representation of the expectation. Also
\begin{equation}
	\mathbb{E}_{p}\left(\sum_{n}M_n\geq1\right) = \sum_{n=4}^{\infty}\mathbb{E}_{p}\left(M_n\right) \leq \sum_{n=4}^{\infty}\left(n 4^n\right)\left(1-p\right)^n
\end{equation}
after bounding the number of complete loops of length $n$, each of which is completely closed with probability $\left(1-p\right)^n$. This implies that
\begin{equation}
	1 - \theta(p) = P_{p}\left(\sum_{n}M_n\geq1\right) \leq \sum_{n=4}^{\infty}\left(n 4^n\right)\left(1-p\right)^n
\end{equation}
with the sum on the r.h.s strictly smaller than unity for some non-zero $\epsilon$ in $p=1-\epsilon$. Whatever this value of $\epsilon$, we have $\theta(p)>0$, and thus $0 < p_c < 1$ given $d=2$.
Also, since we can embed a lower dimensional space in a higher one
\begin{equation}
	1 \geq p_c\left(\mathbb{L}^2\right) \geq p_c\left(\mathbb{L}^3\right) \geq p_c\left(\mathbb{L}^4\right) \geq ... \geq 0
\end{equation}
so $0 < p_c < 1$ for $d \geq 2$ as required.


\chapter{Proof of the isolated vertices theorem} \label{appendix:isolatedvertices}
In this appendix we prove the following important theorem from Section \ref{sec:connect}:
\begin{theorem}[Connectivity is the same as isolated vertices]
	For almost all sets $\mathcal{Y}$
\begin{eqnarray}\label{e:connectivityisthesameasisolatednodes}
	\mathbb{P}\left(G\left(n, \pi r_0^2\right)\textit{ is connected}\right) = \mathbb{P}\left(\textit{no isolated nodes}\right)
\end{eqnarray}
\end{theorem}
\noindent This was first proved by Penrose \cite{penrose1997}. We follow some parts of \cite{waltersreview}. We first show that no two components in the graph are `large'.
\begin{lemma}[No two components in the graph are large]\label{l:nolargecomponents}
		Assuming $c > 0$, there exists a $C$ such that w.h.p the random graph $G\left(n, c \log n\right)$ does not consist of two or more connected components each with \textit{Euclidean diameter}\footnote{This is the largest Euclidean distance which can be found between any two vertices in a connected component.} at least $C\sqrt{\log n}$ . 
	\end{lemma}
		First, tile $S_n$ with tiles (i.e. like on a household wall) of side $r_0/\sqrt{20}$, which insists that any two vertices in $G\left(n, c \log n\right)$ found in any two adjacent squares are no more than $r_0/2$ apart. Then argue that
\begin{enumerate}
	\item A component $U$ of Euclidean diameter at least $C\sqrt{\log n}$ covers many tiles as $n \to \infty$.
	\item Since the tiles have side $r_0/\sqrt{20}$, all tiles adjacent to $U$ are empty.
	\item There are many empty boundary tiles, given the size of $U$.
	\item They cannot all be empty, and so any two large components will merge asymptotically.
\end{enumerate}
We do not explicitly prove parts $1$-$3$, but refer directly to Walter's review \cite{waltersreview} (or Penrose's 1997 paper \cite{penrose1997}). Hopefully these arguments are clear. The crux is to prove that boundary tiles in the limit are so numerous that at least one of them houses at least one vertex (almost surely).

How many boundary tiles surround $U$? At least as many as the square root of the number of tiles underlying $U$. Say $\singlenorm{U^{T}}$ tiles underly $U$. The edge isoperimetric inequality for the grid states\footnote{This is according to Bollob\'{a}s and Leader \cite{bollobas1991}.} that the number of boundary tiles $\singlenorm{\partial U^{T}}$ is given by
\begin{eqnarray}
	\min\left\{2\sqrt{ \singlenorm{U^{T}}},2\sqrt{ \singlenorm{{U^{T}}^{c}}}\right\}
\end{eqnarray}
where ${U^{T}}^{c}$ is the complement of the set of tiles underlying the component $U$. Since each component is of a diameter greater than $C\sqrt{\log n}$, it meets at least\footnote{It could be a diagonal strip-like component over $S_n$, which is as small as its Euclidean diameter allows.} $\left(C\sqrt{\log n}\right)/\left(r/\sqrt{20}\right)$ tiles, so
\begin{eqnarray}
	\singlenorm{\partial U^{T}} \geq 2\sqrt{\frac{C}{r}\sqrt{20\log n}}
\end{eqnarray}
Now we know the size $\singlenorm{\partial U^{T}}$ of the boundary, we show that at least one of its tiles is non-empty. Notice lemma 2.14 in \cite{waltersreview}:
\begin{lemma}\label{l:connectedsubgraphs}
Suppose that $G$ is a graph with maximum degree $\Delta$, and that $v$ is some vertex in $G$. The number of connected subgraphs of $G$ with $n$ vertices that contain $v$ is at most $(e \Delta)^n$.
\end{lemma}
This is not particularly difficult to demonstrate, and we refer to \cite{bollobasbook} for an in depth proof (Problem 45 ``Connected Subgraphs'' in the book of problems by Bollob\'{a}s).

To highlight the main concern, consider the rooted\footnote{A rooted graph is a graph with one marked vertex, called the root. If the graph $G$ has a root vertex we denote it $G^{\star}$.} graph $G^{\star}$, and the number $N(G^{\star},n)$ of subtrees of that graph consisting of $n$ vertices \textit{as well as} the root vertex. Then
\begin{eqnarray}
	N(G^{\star},n) \leq {\Delta n - n + 1 \choose n}
\end{eqnarray}
which requires some work which we exclude. Given
\begin{eqnarray}
	n!{\Delta n - n + 1 \choose n} &=& \left(\Delta-1\right)n\left(\left(\Delta-1\right)^2n^2-1\right)\left(n-3\right)!{(\Delta-1)n-2 \choose n-3} \nonumber \\
	&\leq& e^n \left(\Delta-1\right)^n
\end{eqnarray}
our bound in \ref{l:connectedsubgraphs} follows. The probability that almost all of set of $u$ tiles is empty is
\begin{eqnarray}
	\left(e^{-r^2/20}\right)^{u}
\end{eqnarray}
and (given $\Delta=8$ for the lattice over $S_n$ taking $d=2$), we have
\begin{eqnarray}
	\mathbb{P}\left(\text{two components do not merge}\right) \leq n\left(8e\right)^{u}e^{-ur^2/20}
\end{eqnarray}
after enumerating all possible boundary tile combinations. We consider the graph $G\left(n, c\log n \right)$, so $r = \sqrt{c \log n / \pi}$ and thus if $C$ satisfies
\begin{eqnarray}
	2\sqrt{\frac{C}{r}\sqrt{20\log n}} > 20\pi/c
\end{eqnarray}
then any two components of Euclidean diameter at least $C \sqrt{\log n}$ merge asymptotically. This occurs when
\begin{eqnarray}\label{e:C}
	C\geq\frac{100\pi^2}{c^2\left(\sqrt{20\pi/c}\right)}
\end{eqnarray}
and we have an idea of how large a component can get before it merges. 

Now, all disconnected components must have a Euclidean diameter strictly less than $C \sqrt{\log n}$. Notice that if all of these small disconnected components are a single vertex, then only isolated vertices can disconnect the graph $G\left(n, c\log n \right)$ as $n \to \infty$.
\begin{lemma}[All small components consist of a single vertex.]\label{l:theyareisolatedvertices}
	Suppose $C$ is as in \ref{l:nolargecomponents}, and take the case $c=9/10$. Then the graph $G\left(n, \log n - \frac{1}{2}\log \log n \right)$ contains no components $H$ of 1) more than one vertex and 2) Euclidean diameter strictly less than $C\sqrt{\log n}$.
\end{lemma}
The component $H$ has small Euclidean diameter, and no more than one vertex. Is this diameter smaller than $\eta = \left(\log \log n\right)^2 / \sqrt{\log n}$? Centre a ball of radius $\eta$ at some point $x$ in $H$. Call this ball $B(x,\eta)$. Consider another ball $B(x,r_0)$ with similar center. Then if $\text{diam}(H)<\eta$ is empty then there is at least one more point in $B(x,\eta)$, and $B(x,r_0) \setminus B(x,\eta)$ is empty. The probability of this dies away with $n$, since 
\begin{eqnarray}
	\left(1 - \exp{\left(-\pi\eta^2\right)}\right)\exp{\left(-\left(A-\pi\eta^2\right)\right)} = \mathcal{O}\left(\frac{\sqrt{\log n}\left(\log \log n\right)^{4}}{n \log n}\right) 
\end{eqnarray}
which is bounded by the growth of $1/n$, so even with a chance at every vertex $\mathbb{P}\left(\text{diam}(H)<\eta\right) \to 0$ asymptotically with $n$.
\begin{figure}[!]
\begin{centering}
\includegraphics[scale=0.21]{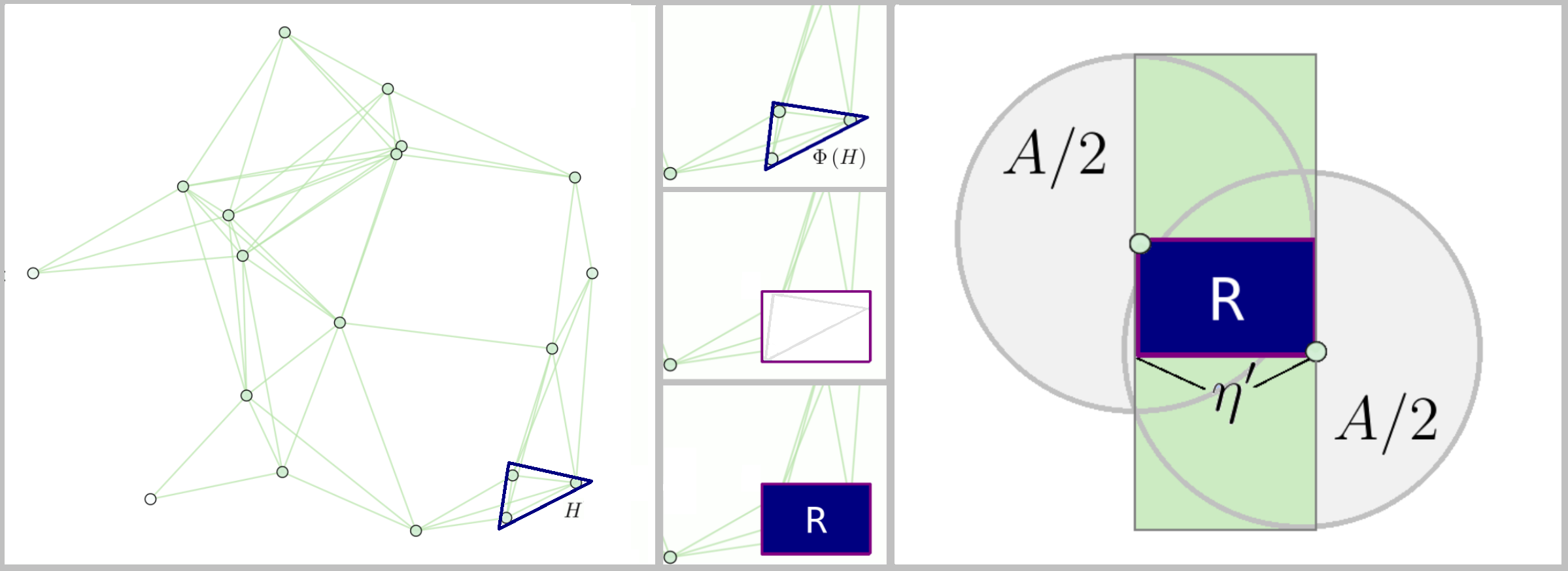}
\caption{The component $H$ is bounded by its convex hull $\Phi\left(H\right)$, itself encased in the upright rectangle $R$. In this case, the exclusion area is non-empty, and the vertex combination fails to satisfy the conditions of the component $H$ in Lemma \ref{l:theyareisolatedvertices}.}\label{fig:Walters}
\end{centering}
\end{figure}
Is $\text{diam}(H)>\eta$? Start by drawing the convex hull\footnote{This is the smallest convex set which encloses all points in $H$.} $\Phi\left(H\right)$ of $H$. Encase $\Phi\left(H\right)$ in an upright rectangle $R$. The diagonal length of $R$ must be at least $\eta$. Fig. \ref{fig:Walters} indicates in grey an area around $R$. $R$ cannot contain less than four points, given it is a rectangle (though $H$ could potentially contain fewer points). So the best situation has points on the north-east and south-west corners defining an area $A_{\text{exclusion}}$ at least
\begin{eqnarray}
	A_{\text{exclusion}} \geq A + \left(1+\text{o}\left(1\right)\right)2 \eta' r_0
\end{eqnarray}
with $\eta'$ the length of the longest side (which must be at least $\eta/\sqrt{2}$, at which point $R$ is a square). Our green rectangle in the right panel of Fig. \ref{fig:Walters} is slightly larger than required, hence the $\text{o}\left(1\right)$ term. Since we have
\begin{multline}
	A + \left(1 + \text{o}\left(1\right)\right)2 r_0 \eta' \geq A + \left(1 + \text{o}\left(1\right)\right)\sqrt{2}r_0 \eta > A + \eta \sqrt{\log n} = A + \left(\log \log n\right)^2 \nonumber
\end{multline}
and thus
\begin{eqnarray}\label{e:exclusion}
	\mathbb{P}\left(A_{\text{exclusion}} \text{ is empty}\right) &\leq& \exp{\left(A + r_0 \eta\right)} = \text{o}\left(\frac{1}{n \left(\log n\right)^{3}}\right) 
\end{eqnarray}
Given there are $n$ vertices available for isolation in $G\left(n, \log n - \frac{1}{2}\log \log n \right)$, the number of four-vertex combinations embedded in $S_n$ which satisfy the geometric constriction imposed\footnote{Note that we can pick any of the $n$ vertices in $G$ to be surrounded by $R$, but the other three must be located strictly within a distance $C\sqrt{\log n}$. There are $\mathcal{O}\left(\log n\right)$ choices per vertex, so there are $\mathcal{O}\left(n \times \log n \times \log n \times \log n\right)$ combinations to try.} by $R$ is order $\mathcal{O}\left(n\left(\log n\right)^{3}\right)$. Thus the number of combinations cannot grow fast enough for the decaying exclusion probability (Eq. \ref{e:exclusion}), and therefore no such rectangle $R$ exists w.h.p. 

Finally, notice the disk scaling is $\log n - \frac{1}{2}\log \log n$, so the graph disconnects w.h.p.
\begin{eqnarray}
	\mathbb{P}\left(G\left(n, \log n - \frac{1}{2}\log \log n\right) \text{ is connected}\right) = e^{-e^{\log \sqrt{\log n}}} \to 0
\end{eqnarray}
and theorem \ref{e:connectivityisthesameasisolatednodes} follows.
\end{appendices}
\newpage
\begin{center}$\dagger$\end{center}
\newpage


\newpage
\begin{center}$\dagger$\end{center}
\newpage


\end{document}